\documentclass[12pt]{article}
\usepackage{amssymb}
\usepackage{amsmath}
\usepackage{latexsym}

\newtheorem{theorem}{Theorem}[section]
\newtheorem{proposition}[theorem]{Proposition}
\newtheorem{corollary}[theorem]{Corollary}
\newtheorem{lemma}[theorem]{Lemma}

\newtheorem{remark}[theorem]{Remark}

\begin{document}

\title{The Genus Two Partition Function for Free Bosonic and Lattice Vertex
Operator Algebras}
\author{Geoffrey Mason\thanks{%
Supported by the NSF, NSA, and the Committee on Research at the University
of California, Santa Cruz} \\
Department of Mathematics, \\
University of California Santa Cruz, \\
CA 95064, U.S.A. \and Michael P. Tuite\thanks{%
Partially supported by Millenium Fund, National University of Ireland,
Galway} \\
Department of Mathematical Physics, \\
National University of Ireland, \\
Galway, Ireland.}
\maketitle

\begin{abstract}
We define the $n$-point function for a vertex operator algebra on a genus two
Riemann surface in two separate sewing schemes where either two tori are
sewn together or a handle is sewn to one torus. We explicitly obtain closed
formulas for the genus two partition function for the Heisenberg free
bosonic string and lattice vertex operator algebras in both sewing schemes.
We prove that the partition functions are holomorphic in the sewing
parameters on given suitable domains and describe their modular properties.
Finally, we show that the partition functions cannot be equal in the
neighborhood of a two-tori degeneration point where they can be explicitly
compared.
\end{abstract}
\newpage
\begin{center}
\textbf{Table of Contents }
\end{center}

\medskip \noindent 1. Introduction \newline
2. Genus two Riemann surfaces from sewn tori \newline
\indent2.1 Some elliptic functions \newline
\indent2.2 The $\epsilon $-formalism for sewing two tori \newline
\indent2.3 The $\rho $-formalism for self-sewing a torus \newline
3. Graphical expansions \newline
\indent3.1 Rotationless and chequered cycles \newline
\indent3.2 Graphical expansions for $\Omega $ in the $\epsilon $-formalism
\newline
\indent3.3 Graphical expansions for $\Omega $ in the $\rho $-formalism
\newline
4. Vertex operator algebras and the Li-Zamolodchikov metric \newline
\indent4.1 Vertex operator algebras \newline
\indent4.2 The Li-Zamolodchikov metric \newline
5. Partition and $n$-point functions for a vertex operator algebra on a
Riemann surface \newline
\indent5.1 Genus zero \newline
\indent5.2 Genus one \newline
\indent5.2.1 Self-sewing the Riemann sphere \newline
\indent5.2.2 An alternative self-sewing of the Riemann sphere and the
Catalan series \newline
\indent5.3 Genus two \newline
6. Genus two partition function for free bosonic theories in the $\epsilon $%
-formalism \newline
\indent6.1 The genus two partition function $Z_{M}^{(2)}(\tau _{1},\tau
_{2},\epsilon )$ \newline
\indent6.2 Holomorphic and modular-invariance properties \newline
7. Genus two partition function for lattice theories in the $\epsilon $%
-formalism \newline
8. Genus two partition function for free bosonic theories in the $\rho $%
-formalism \newline
\indent8.1 The genus two partition function $Z_{M}^{(2)}(\tau ,w,\rho )$
\newline
\indent8.2 Holomorphic and modular-invariance properties \newline
9. Genus two partition function for lattice theories in the $\rho $%
-formalism \newline
10. Comparison of partition functions in the two formalisms \newline
11. Final Remarks \newline
12. Appendix

\indent12.1 A product formula

\indent12.2 Invertibility about a two tori degeneration point

\indent12.3 Corrections

\section{Introduction}

One of the most striking features of conformal field theory (aka the theory
of vertex algebras) is the occurrence of \emph{elliptic functions} and \emph{%
modular forms}, manifested in the form of $n$-point correlation trace
functions. This phenomenon has been present in the physics of strings since
the earliest days e.g. \cite{GSW}, \cite{P}. In mathematics it dates from
the Conway-Norton conjectures \cite{CN} proved by Borcherds (\cite{B1}, \cite%
{B2}), and Zhu's important paper \cite{Z1}. Geometrically, we are dealing
with probability amplitudes corresponding to a complex torus (compact
Riemann surface of genus one) inflicted with $n$ punctures corresponding to
local fields (vertex operators). For a Vertex Operator Algebra (VOA) $%
V=\oplus V_{n}$, the most familiar correlation function is the $0$-point
function, also called the \emph{graded dimension} or \emph{partition function}
\begin{equation}
Z_{V}(q)=q^{-c/24}\sum_{n}\dim V_{n}q^{n}  \label{grdim}
\end{equation}%
($c$ is the central charge). An example which motivates much of the present
paper is that of a lattice theory $V_{L}$ associated to a positive-definite
even lattice $L$. Then $c$ is the rank of $L$ and
\begin{equation}
Z_{V_{L}}(q)=\frac{\theta _{L}(q)}{q^{c/24}\prod_{n}(1-q^{n})^{c}}
\label{latticeg=1pfunc}
\end{equation}%
where $\theta _{L}(q)$ is the usual theta function of $L$ and the
denominator is the $c$th power of the Dedekind eta function $\eta (q)$. Both
of these functions are (holomorphic) elliptic modular forms of weight $c/2$
on a certain congruence subgroup of $SL(2,\mathbb{Z})$, so that $Z_{V_{L}}$
is an elliptic modular function of weight zero on the same subgroup. It is
widely expected that an analogous result holds for \emph{any rational}
vertex operator algebra, namely that $Z_{V}(q)$ is a modular function of
weight zero on a congruence subgroup of $SL(2,\mathbb{Z})$.

\medskip There are natural physical and mathematical reasons for wanting to
extend this picture to Riemann surfaces of \emph{higher genus}. In
particular, we want to know if there are natural analogs of (\ref{grdim})
and (\ref{latticeg=1pfunc}) for arbitrary rational vertex operator algebras
and arbitrary genus, in which \emph{genus g Siegel modular forms} occur.
This is much more challenging than the case of genus one. Many, but not all,
of the new difficulties that arise are already present at genus two, and it
is this case that we are concerned with in the present paper. Our goal,
then, is this: given a vertex operator algebra $V$, to introduce the $n$%
-point correlation functions on a compact Riemann surface of genus two which
are associated to $V$, and study their automorphic properties. An overview
of aspects of this program is given in the Introduction to \cite{MT2}. Brief
discussions of some of our methods and results can also be found in \cite{T}
and \cite{MT3}.

\medskip Generally, our approach to genus two correlation functions is to
define them in terms of genus one data coming from $V$. Now there are two
rather different ways to obtain a compact Riemann surface of genus two from
surfaces of genus one - one may sew two separate tori together, or self-sew
a torus (i.e. attach a handle). We refer to these two schemes as the $%
\epsilon $- and $\rho $-formalism respectively. As a result, much of our
discussion is bifurcated as we are obliged to treat the two cases
separately. The present paper depends heavily on results obtained in \cite%
{MT2}, reviewed in Section 2 below, concerned with a pair of maps
\begin{equation}
\mathcal{D}^{\epsilon }\overset{F^{\epsilon }}{\longrightarrow }\mathbb{H}%
_{2}\overset{F^{\rho }}{\longleftarrow }\mathcal{D}^{\rho }  \label{Fmaps1}
\end{equation}%
For $g\geq 1,\mathbb{H}_{g}$ is the genus $g$ Siegel upper half-space. Then $%
\mathcal{D}^{\epsilon }\subseteq \mathbb{H}_{1}\times \mathbb{H}_{1}\times
\mathbb{C}$ is the domain consisting of triples $(\tau _{1},\tau
_{2},\epsilon )$ which correspond to a pair of complex tori of modulus $\tau
_{1},\tau _{2}$ sewn together according to the relation $z_{1}z_{2}=\epsilon
$. (For more details, see \cite{MT2} and Section 2 below.) $\mathcal{D}%
^{\rho }\subseteq \mathbb{H}_{1}\times \mathbb{C}^{2}$ is similarly defined
in terms of data $(\tau ,w,\rho )$ needed to self-sew a torus of modulus $%
\tau $. In each case, sewing produces a compact Riemann surface of genus
two, and the maps $F^{\bullet }$ are those which assign to a point of $%
\mathcal{D}^{\bullet }$ the period matrix $\Omega $ of the sewn surface.
Thus $F^{\bullet }$ takes values in $\mathbb{H}_{2},$ and in \cite{MT2} we
established that the $F^{\bullet }$ are holomorphic. (Here and below, it is
sometimes convenient to use a bullet in place of subscripts and superscripts
when comparing the $\epsilon $- and $\rho $-formalisms).

\medskip In Section 3 we introduce some graph-theoretic technology which
provides a convenient way of describing the period matrix $\Omega $ in the $%
\epsilon $- and $\rho $-formalisms. Similar graphical techniques are
employed later on as means of computing the genus two partition function for
free bosonic (Heisenberg)  and lattice VOAs. Section 4 is a brief review of
some necessary background on VOA theory  and the Li-Zamalodchikov or Li-Z
metric.  We assume throughout that the Li-Z  metric is unique and invertible
(which follows if $V$ is simple \cite{Li}).

\medskip Section 5 develops a theory of $n$-point functions for VOAs on
Riemann surfaces of genus 0, 1 and 2 motivated by ideas in conformal field
theory (\cite{FS}, \cite{So1}, \cite{P}). The Zhu theory \cite{Z1} of genus
one $n$-point functions is reformulated in this language in terms of the
self-sewing a Riemann sphere to obtain a torus. We consider two examples of
such a sewing procedure where an interesting application of the Catalan
series from combinatorics arises. We next consider two separate definitions
of genus two $n$-point functions based on the $\epsilon $- and $\rho $%
-formalisms.  The genus two partition functions involve extending (\ref%
{Fmaps1}) to a diagram
\begin{equation*}
\begin{array}{ccccc}
\mathcal{D}^{\epsilon } & \overset{F^{\epsilon }}{\longrightarrow } &
\hspace{0.3cm}\mathbb{H}_{2} & \overset{F^{\rho }}{\longleftarrow } &
\mathcal{D}^{\rho } \\
& \searrow  & \hspace{0.3cm}\downarrow  & \hspace{0.5cm}\swarrow  &  \\
&  & \hspace{0.2cm}\mathbb{C} &  &
\end{array}%
\end{equation*}%
where the partition functions are maps $\mathcal{D}^{\bullet }\rightarrow
\mathbb{C}$, defined purely in terms of genus one data coming from $V$.
Explicitly,  the genus two partition function of $V$ in the two formalisms
are as follows:
\begin{eqnarray}
Z_{V,\epsilon }^{(2)}(\tau _{1},\tau _{2},\epsilon ) &=&\sum_{n\geq
0}\epsilon ^{n}\sum_{u\in V_{[n]}}Z_{V}^{(1)}(u,\tau _{1})Z_{V}^{(1)}(\bar{u}%
,\tau _{2}),  \label{pfeps} \\
Z_{V,\rho }^{(2)}(\tau ,w,\rho ) &=&\sum_{n\geq 0}\rho ^{n}\sum_{u\in
V_{[n]}}Z_{V}^{(1)}(\bar{u},u,w,\tau ).  \label{pfrho}
\end{eqnarray}%
Here, $Z_{V}^{(1)}(u,\tau )$ and $Z_{V}^{(1)}(\bar{u},u,w,\tau )$ are (genus
one) $1$- and $2$-point functions respectively with $\bar{u}$ the Li-Z
metric dual of $u$. The precise meaning of (\ref{pfeps}) and (\ref{pfrho}),
together with similar definitions for $n$-point functions, is given in
Section 5. In any case, at genus two we get not one, but two, rather
differently defined partition functions of $V$.

\medskip In Sections 6-9 we investigate the case of free bosonic theories
(i.e. the Heisenberg vertex operator algebra $M^{c}$ corresponding to $c$
free bosons) and lattice theories in depth. The definitions of the genus two
partition functions are \emph{a priori} just formal power series in
variables $\epsilon ,q_{1},q_{2}$ and $\rho ,w,q$ respectively (as usual, $%
q=e^{2\pi i\tau }$, etc.). However, we will see, at least in the case of $%
M^{c}$ and $V_{L}$, that they are in fact holomorphic functions on $\mathcal{%
D}^{\bullet }$. It is natural to expect that this result holds in much wider
generality. Section 6 is devoted to the case of $M^{c}$ in the $\epsilon $%
-formalism. In this case, holomorphy depends on an interesting new formula
for the genus two partition function. Namely, we prove (Theorem \ref%
{Theorem_Z2_boson}) by reinterpreting (\ref{pfeps}) in terms of certain
graphical expansion, that
\begin{equation}
Z_{M^{c},\epsilon }^{(2)}(\tau _{1},\tau _{2},\epsilon )=\frac{%
Z_{M^{c}}^{(1)}(\tau _{1})Z_{M^{c}}^{(1)}(\tau _{2})}{\det
(I-A_{1}A_{2})^{c/2}}.  \label{Z2M}
\end{equation}%
Here, the $A_{i}$ are certain infinite matrices whose entries are
expressions involving quasi-modular forms, and $Z_{M^{c}}^{(1)}(\tau
_{i})=q_{i}^{c/24}/\eta (q_{i})$ the corresponding genus one partition
function\footnote{%
In (\ref{pfeps}), (\ref{pfrho}) there are no overall factors $q_{i}^{-c/24}$%
, whereas in (\ref{latticeg=1pfunc}) there is.}. The matrices $A_{i}$ and
the infinite determinant that occurs in (\ref{Z2M}) were introduced and
discussed at length in \cite{MT2}. The results obtained there are important
here, as are the explicit computations of $1$-point functions obtained in
\cite{MT1}. We also give in Section 6 a product formula for the infinite
determinant (Theorem \ref{Theorem_Z2_boson_prod}) which depends on the
graphical interpretation of the entries of the $A_{i}$.

\medskip The domain $\mathcal{D}^{\epsilon }$ admits the group $G_{0}=SL(2,%
\mathbb{Z})\times SL(2,\mathbb{Z})$ as automorphisms (in fact, there is a
larger automorphism group $G$ that contains $G_{0}$ with index $2$.) We show
(cf. Theorem \ref{Theorem_Z2_G}) that the modular partition function  $%
Z_{M^{c},\mathrm{mod}}^{(2)}(\tau _{1},\tau _{2},\epsilon )\equiv
Z_{M^{c},\epsilon }^{(2)}q_{1}^{-c/24}q_{2}^{-c/24}$ is an automorphic form
of weight $-c/2$ on $G$. This is a bit imprecise in several ways: we have
not explained here what the automorphy factor is, and in fact this is an
interesting point because it depends on the map $F^{\epsilon }$. And
similarly to the eta-function, there is a $24$th root of unity,
corresponding to a character of $G$, that intervenes in the functional
equation. In any case if, for example, the central charge is $26$ (the case
of relevance to the critical bosonic string), $Z_{M^{26},\mathrm{mod}}^{(2)}$
is an automorphic form of weight $-13$ on $G$ transforming according to an
explicit character of $G$ of order $12$. These properties of $Z_{M^{c},%
\mathrm{mod}}^{(2)}(\tau _{1},\tau _{2},\epsilon )$ justify the idea that it
should be thought of as the genus two analog (in the $\epsilon $-formalism)
of $\eta (q)^{-c}$.

\medskip Section 8 is devoted to a reprise of the contents of Section 6 in
the $\rho $-formalism. This is not merely a simple translation, however, and
there are several issues in the $\rho $-formalism that need additional
attention and which generally make the $\rho $-formalism more complicated
than its $\epsilon$-counterpart. This circumstance was already
evident in \cite{MT2},
and stems in part from the fact that the map $F^{\rho }$ involves a
logarithmic term that is absent in the $\epsilon $-formalism. The $A$%
-matrices are also more unwieldy, as they must be considered as having
entries which themselves are $2\times 2$ block-matrices with entries which
are elliptic-type functions. Be that as it may, we establish in Theorem \ref%
{Theorem_Z2_boson_rho} the analog of (\ref{Z2M}) for $Z_{M^{c},\rho
}^{(2)}(\tau ,w,\rho )$ in the $\rho $-formalism, a product formula (Theorem %
\ref{Theorem_Z2_boson_prod_rho}) for the intervening infinite determinant,
and automorphy of $Z_{M^{c},\mathrm{mod}}^{(2)}(\tau ,w,\rho )\equiv $ $%
Z_{M^{c},\rho }^{(2)}q^{-c/24}$ with respect to a group $\Gamma _{1}\cong SL(%
\mathbb{Z})$ (Theorem \ref{Theorem_Z2_rho_G1}). In particular we obtain the
holomorphy of the partition function.

\medskip
In view of these strong formal similarities between the two partition
functions $Z_{M^{c},\epsilon }^{(2)}(\tau _{1},\tau _{2},\epsilon )$ and $%
Z_{M^{c},\rho }^{(2)}(\tau ,w,\rho )$ associated to $2$ free bosons, it is
natural to ask if they are equal in some sense. One could ask the same
question for rational theories $V$ too. In the very special case in which $V$
is holomorphic (i.e. it has a \emph{unique} irreducible module) one knows
(e.g. \cite{TUY}) that the genus $2$ conformal block is $1$-dimensional, in
which case an identification of the two partition functions might seem
inevitable. Of course, the partition functions are defined on different
domains, so there is no question of them being literally equal. What one
could (and should) ask is the following: is there a holomorphic map $f:%
\mathbb{H}_{2}\rightarrow \mathbb{C}$ which makes the following diagram
commute?

\begin{equation}
\begin{array}{lllll}
\ \hspace{0cm}D^{\epsilon } & \overset{F^{\epsilon }}{\longrightarrow } & \
\ \mathbb{H}_{2} & \overset{F^{\rho }}{\longleftarrow } & D^{\rho } \\
& \overset{Z_{\epsilon }^{(2)}}{\searrow } & \ \ \ \downarrow f & \overset{%
Z_{\rho }^{(2)}}{\swarrow } &  \\
&  &  &  &  \\
&  & \hspace{0.4cm}\mathbb{C} &  &
\end{array}
\label{fexists?}
\end{equation}%
One should also expect $f$ to be automorphic in a way that is consistent
with the automorphic properties of the partition functions. Prior comments
notwithstanding, we will explain in Section 10 that there is \emph{no} such $%
f$, even for holomorphic theories such as $V_{L}$ (where the lattice $L$ is
self-dual). How can we understand this situation? In the holomorphic case,
say, how can it be reconciled with the uniqueness of the conformal block?

\medskip In Sections 7 and 9 we elucidate the expressions (\ref{pfeps}) and (%
\ref{pfrho}), thereby obtaining some of the main results of the present
paper. Taking $L$ of rank $c$, we prove that
\begin{equation}
\frac{Z_{V_{L},\epsilon }^{(2)}(\tau _{1},\tau _{2},\epsilon )}{%
Z_{M^{c},\epsilon }^{(2)}(\tau _{1},\tau _{2},\epsilon )}=\frac{%
Z_{V_{L},\rho }^{(2)}(\tau ,w,\rho )}{Z_{M^{c},\rho }^{(2)}(\tau ,w,\rho )}%
=\theta _{L}^{(2)}(\Omega ).  \label{normpfs}
\end{equation}%
Here, $\theta _{L}^{(2)}(\Omega )$ is the genus two Siegel theta function
attached to $L$ (\cite{F}), and the equalities hold (even as formal power
series) when $\Omega$ is interpreted as $F^{\epsilon}(\tau_1, \tau_2,
\epsilon)$ or $F^{\rho}(\tau, w, \rho)$ according to the sewing scheme.
Inasmuch as the genus one partition function for $c$ free bosons is $\eta
(q)^{-c}$, we see that (\ref{latticeg=1pfunc}) is just the genus one analog
of (\ref{normpfs}).

\medskip Call the displayed quotients of partition functions the \emph{%
normalized} partition function. Then (\ref{normpfs}) says exactly that (\ref%
{fexists?}) holds (with $f$ being the Siegel theta function) if we replace
the partition functions with their normalized versions. Thus the following
diagram commutes

\begin{equation*}
\begin{array}{lllll}
\ \hspace{0cm}D^{\epsilon } & \overset{F^{\epsilon }}{\longrightarrow } & \
\ \mathbb{H}_2 & \overset{F^{\rho }}{\longleftarrow } & D^{\rho } \\
& \overset{\hat{Z}_{\epsilon }^{(2)}}{\searrow } & \ \ \ \downarrow \theta
_{L}^{(2)} & \overset{\hat{Z}_{\rho }^{(2)}}{\swarrow } &  \\
&  &  &  &  \\
&  & \hspace{0.4cm}\mathbb{C} &  &
\end{array}%
\end{equation*}

\medskip We can put it this way: the \emph{normalized} partition functions
are \emph{independent of the sewing scheme}. They can be identified, via the
sewing maps $F^{\bullet }$, with \emph{a holomorphic genus two Siegel
modular form of weight $c/2$ - the Siegel theta function in this case}. It
is the normalized partition function which can be appropriately identified
with an element of the conformal block. We conjecture that analogous results
hold for \emph{any} rational vertex operator algebra and \emph{any} genus.
Section 11 contains a brief further discussion of these issues in the light
of related ideas in string theory and algebraic geometry. That the two
normalized partition functions turn out to be equal is, from our current
perspective, something of a miracle. It would be very useful to have
available an \emph{a priori} proof of this circumstance.

\section{Genus Two Riemann Surfaces from Sewn Tori}

In this section we review the main results of  \cite{MT2} relevant to
the present work. We review two separate constructions of a genus two
Riemann surface based on a general sewing formalism due to Yamada \cite{Y}.
\ In the first construction, which we refer to as the $\epsilon $-formalism,
we parameterize a genus two Riemann surface by sewing together two
once-punctured tori. In the second construction, which we refer to as the $%
\rho $-formalism, we parameterize a genus two Riemann surface by self-sewing
a twice-punctured torus. In both cases, the period matrix is described by an
explicit formula which defines a holomorphic map from a specified domain
into the genus two Siegel upper half plane $\mathbb{H}_{2}$. In each case,
this map is equivariant under a suitable subgroup of $Sp(4,\mathbb{Z})$. We
also review the convergence and holomorphy of some infinite determinants
that naturally arise and which play a dominant r\^{o}le later on.

\subsection{Some Elliptic Function Theory}

We begin with the definition of various modular and elliptic functions that
permeate this work \cite{MT1}, \cite{MT2}. We define
\begin{eqnarray}
P_{2}(\tau ,z) &=&\wp (\tau ,z)+E_{2}(\tau )  \notag \\
&=&\frac{1}{z^{2}}+\sum_{k=2}^{\infty }(k-1)E_{k}(\tau )z^{k-2},  \label{P2}
\end{eqnarray}%
where $\tau $ $\in \mathbb{H}_{1}$, the complex upper half-plane and where $%
\wp (\tau ,z)$ is the Weierstrass function and $E_{k}(\tau )$ is equal to $0$
for $k$ odd, and for $k$ even is the Eisenstein series
\begin{equation}
E_{k}(\tau )=E_{k}(q)=-\frac{B_{k}}{k!}+\frac{2}{(k-1)!}\sum_{n\geq 1}\sigma
_{k-1}(n)q^{n}.  \label{Eisenk}
\end{equation}%
Here and below, we take $q=\exp (2\pi i\tau )$; $\sigma
_{k-1}(n)=\sum_{d\mid n}d^{k-1}$, and $B_{k}$ is a $k$th Bernoulli number
e.g. \cite{Se}. If $k\geq 4$ then $E_{k}(\tau )$ is a holomorphic modular
form of weight $k$ on $SL(2,\mathbb{Z})$ whereas $E_{2}(\tau )$ is a
quasi-modular form \cite{KZ},\cite{MT2}. We define $P_{0}(\tau ,z)$, up to a
choice of the logarithmic branch, and $P_{1}(\tau ,z)$ by
\begin{eqnarray}
P_{0}(\tau ,z) &=&-\log (z)+\sum_{k\geq 2}\frac{1}{k}E_{k}(\tau )z^{k},
\label{P0} \\
P_{1}(\tau ,z) &=&\frac{1}{z}-\sum_{k\geq 2}E_{k}(\tau )z^{k-1},  \label{P1}
\end{eqnarray}%
where $P_{2}=-\frac{d}{dz}P_{1}$ and $P_{1}=-\frac{d}{dz}P_{0}.$ $P_{0}$ is
related to the elliptic prime form $K(\tau ,z)\,$\ by \cite{Mu1}
\begin{equation}
K(\tau ,z)=\exp (-P_{0}(\tau ,z)).  \label{Primeform}
\end{equation}%
Define elliptic functions $P_{k}(\tau ,z)\,$ for $k\geq 3$
\begin{equation}
P_{k}(\tau ,z)=\frac{(-1)^{k-1}}{(k-1)!}\frac{d^{k-1}}{dz^{k-1}}P_{1}(\tau
,z).  \label{Pkdef}
\end{equation}

\noindent Define for $k,l\geq 1$
\begin{eqnarray}
C(k,l) &=&C(k,l,\tau )=(-1)^{k+1}\frac{(k+l-1)!}{(k-1)!(l-1)!}E_{k+l}(\tau ),
\label{Ckldef} \\
D(k,l,z) &=&D(k,l,\tau ,z)=(-1)^{k+1}\frac{(k+l-1)!}{(k-1)!(l-1)!}%
P_{k+l}(\tau ,z).  \label{Dkldef}
\end{eqnarray}%
$\,$Note that $C(k,l)=C(l,k)$ and $D(k,l,z)=(-1)^{k+l}D(l,k,z)$. We also
define for $k\geq 1$,

\begin{eqnarray}
C(k,0) &=&C(k,0,\tau )=(-1)^{k+1}E_{k}(\tau ),  \label{C(k,0)} \\
D(k,0,z,\tau ) &=&(-1)^{k+1}P_{k}(z,\tau ).  \label{D(k,0)}
\end{eqnarray}%
The Dedekind eta-function is defined by

\begin{equation}
\eta (\tau )=q^{1/24}\prod_{n=1}^{\infty }(1-q^{n}).  \label{etafun}
\end{equation}

\subsection{The $\protect\epsilon $-Formalism for Sewing Two Tori}

We now review a general method due to Yamada \cite{Y} and discussed at
length in  \cite{MT2} for calculating the period matrix $\Omega $ of the
genus two Riemann surface formed by sewing together two tori $\mathcal{S}%
_{a} $ for $a=1,2$. We shall sometimes refer to $\mathcal{S}_{1}$ and $%
\mathcal{S}_{2}$ as the left and right torus respectively. Consider an
oriented torus $\mathcal{S}_{a}=\mathbb{C}/\Lambda _{a}$ with lattice $%
\Lambda _{a}=2\pi i(\mathbb{Z}\tau _{a}\oplus \mathbb{Z})$ for $\tau _{a}\in
\mathbb{H}_{1}$. For local coordinate $z_{a}\in \mathbb{C}/\Lambda _{a}$
consider the closed disk $\left\vert z_{a}\right\vert \leq r_{a}$ which is
contained in $\mathcal{S}_{a}$ provided\ $r_{a}<\frac{1}{2}D(q_{a})$ where
\begin{equation*}
D(q_{a})=\min_{\lambda \in \Lambda _{a},\lambda \neq 0}|\lambda |,
\end{equation*}%
is the minimal lattice distance. Introduce a complex sewing parameter $%
\epsilon $ where $|\epsilon |\leq r_{1}r_{2}<\frac{1}{4}D(q_{1})D(q_{2})$
and excise the disk $\{z_{a},\left\vert z_{a}\right\vert \leq |\epsilon |/r_{%
\bar{a}}\}$ centred at $z_{a}=0$ where we use the convention
\begin{equation}
\overline{1}=2,\quad \overline{2}=1.  \label{abar}
\end{equation}%
Defining the annulus

\begin{equation}
\mathcal{A}_{a}=\{z_{a},|\epsilon |/r_{\bar{a}}\leq \left\vert
z_{a}\right\vert \leq r_{a}\},  \label{annuli_eps}
\end{equation}%
we identify $\mathcal{A}_{1}$ with $\mathcal{A}_{2}$ via the sewing relation
\begin{equation}
z_{1}z_{2}=\epsilon .  \label{pinch}
\end{equation}

\begin{center}
\bigskip

\begin{picture}(300,100)

\put(50,50){\qbezier(10,-20)(-30,0)(10,20)}
\put(50,52){\qbezier(10,18)(50,35)(90,18)}

\put(50,48){\qbezier(10,-18)(50,-35)(90,-18)}

\put(45,50){\qbezier(25,0)(45,17)(60,0)}
\put(45,50){\qbezier(20,2)(45,-17)(65,2)}


\put(175,52){\qbezier(10,18)(50,35)(90,18)}
\put(175,50){\qbezier(90,20)(130,0)(90,-20)}
\put(175,48){\qbezier(10,-18)(50,-35)(90,-18)}

\put(200,50){\qbezier(25,0)(45,17)(60,0)}
\put(200,50){\qbezier(20,2)(45,-17)(65,2)}

\put(140,50){\circle{16}}
\put(140,50){\circle{40}}

\put(140,50){\vector(-1,-2){0}}
\put(50,50){\qbezier(90,0)(100,15)(90,30)}%
\put(140,90){\makebox(0,0){$z_1=0$}}

\put(140,50){\line(-1,1){14.1}}
\put(127,55){\makebox(0,0){$r_1$}}

\put(140,50){\line(1,0){8}}
\put(145,50){\vector(1,4){0}}
\put(55,20){\qbezier(90,4)(85,17)(90,30)}%
\put(150,15){\makebox(0,0){$|\epsilon|/r_2$}}

\put(25,50){\makebox(0,0){$\mathcal{S}_1$}}


\put(185,50){\circle{16}}
\put(185,50){\circle{40}}

\put(185,50){\vector(1,-2){0}}
\put(95,50){\qbezier(90,0)(80,15)(90,30)}%
\put(185,90){\makebox(0,0){$z_2=0$}}

\put(185,50){\line(-1,-1){14.1}}
\put(171,45){\makebox(0,0){$r_2$}}

\put(185,50){\line(1,0){8}}
\put(190,50){\vector(1,4){0}}
\put(90,20){\qbezier(100,4)(95,17)(100,30)}%
\put(190,15){\makebox(0,0){$ |\epsilon|/r_1$}}

\put(300,50){\makebox(0,0){$\mathcal{S}_2$}}

\end{picture}

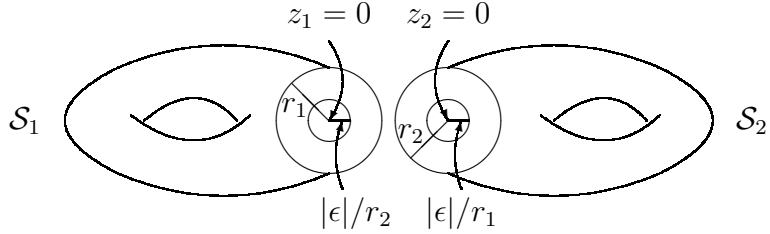
{\small Fig. 1 Sewing Two Tori}
\end{center}

The genus two Riemann surface is then parameterized by the domain
\begin{equation}
\mathcal{D}^{\epsilon }=\{(\tau _{1},\tau _{2},\epsilon )\in \mathbb{H}_{1}%
\mathbb{\times H}_{1}\mathbb{\times C}\ |\ |\epsilon |<\frac{1}{4}%
D(q_{1})D(q_{2})\}.  \label{Deps}
\end{equation}

The genus two period matrix $\Omega \in \mathbb{H}_{2}$, the Siegel upper
half plane, may be determined as a function of\ $(\tau _{1},\tau
_{2},\epsilon )\in \mathcal{D}^{\epsilon }$ in terms of certain infinite
dimensional moment matrices $A_{a}(\tau _{a},\epsilon )=(A_{a}(k,l,\tau
_{a},\epsilon ))$ for $k,l\geq 1$ where
\begin{equation}
A_{a}(k,l,\tau _{a},\epsilon )=\frac{\epsilon ^{(k+l)/2}}{\sqrt{kl}}%
C(k,l,\tau _{a}).  \label{Akldef}
\end{equation}%
These matrices play a dominant role both in the description of $\Omega $ and
in our later discussion of the free bosonic and lattice VOA genus two
partition functions. In particular, the matrix $I-A_{1}A_{2}$ and $\det
(I-A_{1}A_{2})$ (where $I$ denotes the infinite identity matrix) play an
important role where $\det (I-A_{1}A_{2})$ is defined by
\begin{eqnarray}
\log \det (I-A_{1}A_{2}) &=&tr\log (I-A_{1}A_{2})  \notag \\
&=&-\sum_{n\geq 1}\frac{1}{n}tr((A_{1}A_{2})^{n}).  \label{logdet}
\end{eqnarray}%
One finds

\begin{theorem}
\label{Theorem_A1A2} \ \

(a) (op. cite., Proposition 1) The infinite matrix
\begin{equation}
(I-A_{1}A_{2})^{-1}=\sum_{n\geq 0}(A_{1}A_{2})^{n},  \label{I_minus_A1A2}
\end{equation}%
is convergent for $(\tau _{1},\tau _{2},\epsilon )\in \mathcal{D}^{\epsilon
} $.

(b) (op. cite., Theorem 2 \& Proposition 3) $\det (I-A_{1}A_{2})$ is
non-vanishing and holomorphic for $(\tau _{1},\tau _{2},\epsilon )\in
\mathcal{D}^{\epsilon }$. $\square $
\end{theorem}

Furthermore, the genus two period matrix $\Omega $ is given by:

\begin{theorem}
\label{Theorem_period_eps}(op. cite., Theorem 4) The $\epsilon $-formalism
determines a holomorphic map
\begin{eqnarray}
F^{\epsilon }:\mathcal{D}^{\epsilon } &\rightarrow &\mathbb{H}_{2},  \notag
\\
(\tau _{1},\tau _{2},\epsilon ) &\mapsto &\Omega (\tau _{1},\tau
_{2},\epsilon ),  \label{Fepsmap}
\end{eqnarray}%
where $\Omega =\Omega (\tau _{1},\tau _{2},\epsilon )$ is given by
\begin{eqnarray}
2\pi i\Omega _{11} &=&2\pi i\tau _{1}+\epsilon
(A_{2}(I-A_{1}A_{2})^{-1})(1,1),  \label{Om11eps} \\
2\pi i\Omega _{22} &=&2\pi i\tau _{2}+\epsilon
(A_{1}(I-A_{2}A_{1})^{-1})(1,1),  \label{Om22eps} \\
2\pi i\Omega _{12} &=&-\epsilon (I-A_{1}A_{2})^{-1}(1,1).  \label{Om12eps}
\end{eqnarray}%
Here $(1,1)$ refers to the $(1,1)$-entry of a matrix. $\square $
\end{theorem}

$\mathcal{D}^{\epsilon }$ is preserved under the action of $G\simeq (SL(2,%
\mathbb{Z})$ $\times SL(2,\mathbb{Z}))\rtimes \mathbb{Z}_{2}$, the direct
product of two copies of $SL(2,\mathbb{Z})$ (the left and right torus
modular groups) which are interchanged upon conjugation by an involution $%
\beta $ as follows%
\begin{eqnarray}
\gamma _{1}.(\tau _{1},\tau _{2},\epsilon ) &=&(\frac{a_{1}\tau _{1}+b_{1}}{%
c_{1}\tau _{1}+d_{1}},\tau _{2},\frac{\epsilon }{c_{1}\tau _{1}+d_{1}}),
\notag \\
\gamma _{2}.(\tau _{1},\tau _{2},\epsilon ) &=&(\tau _{1},\frac{a_{2}\tau
_{2}+b_{2}}{c_{2}\tau _{2}+d_{2}},\frac{\epsilon }{c_{2}\tau _{2}+d_{2}}),
\notag \\
\beta .(\tau _{1},\tau _{2},\epsilon ) &=&(\tau _{2},\tau _{1},\epsilon ),
\label{GDeps}
\end{eqnarray}%
for $(\gamma _{1},\gamma _{2})\in SL(2,\mathbb{Z})\times SL(2,\mathbb{Z})$
with $\gamma _{i}=\left(
\begin{array}{cc}
a_{i} & b_{i} \\
c_{i} & d_{i}%
\end{array}%
\right) $. There is a natural injection $G\rightarrow Sp(4,\mathbb{Z})$ in
which the two $SL(2,\mathbb{Z})$ subgroups are mapped to
\begin{equation}
\Gamma _{1}=\left\{ \left[
\begin{array}{cccc}
a_{1} & 0 & b_{1} & 0 \\
0 & 1 & 0 & 0 \\
c_{1} & 0 & d_{1} & 0 \\
0 & 0 & 0 & 1%
\end{array}%
\right] \right\} ,\;\Gamma _{2}=\left\{ \left[
\begin{array}{cccc}
1 & 0 & 0 & 0 \\
0 & a_{2} & 0 & b_{2} \\
0 & 0 & 1 & 0 \\
0 & c_{2} & 0 & d_{2}%
\end{array}%
\right] \right\} ,  \label{Gamma1Gamma2}
\end{equation}%
and the involution is mapped to
\begin{equation}
\beta =\left[
\begin{array}{cccc}
0 & 1 & 0 & 0 \\
1 & 0 & 0 & 0 \\
0 & 0 & 0 & 1 \\
0 & 0 & 1 & 0%
\end{array}%
\right] .  \label{betagen}
\end{equation}%
Thus as a subgroup of $Sp(4,\mathbb{Z})$, $G$ also has a natural action on
the Siegel upper half plane $\mathbb{H}_{2}$ where for $\gamma =\left(
\begin{array}{ll}
A & B \\
C & D%
\end{array}%
\right) \in Sp(4,\mathbb{Z})$
\begin{equation}
\gamma .\Omega {=(A\Omega +B)(C\Omega +D)^{-1},}  \label{eq: modtrans}
\end{equation}%
One then finds

\begin{theorem}
\label{TheoremGequiv} (op. cit., Theorem 5) $F^{\epsilon }$ is equivariant
with respect to the action of $G$ i.e. there is a commutative diagram for $%
\gamma \in G$,
\begin{equation*}
\begin{array}{ccc}
\mathcal{D}^{\epsilon } & \overset{F^{\epsilon }}{\rightarrow } & \mathbb{H}%
_{2} \\
\gamma \downarrow &  & \downarrow \gamma \\
\mathcal{D}^{\epsilon } & \overset{F^{\epsilon }}{\rightarrow } & \mathbb{H}%
_{2}%
\end{array}%
\end{equation*}
$\square$
\end{theorem}

The limit $\epsilon \rightarrow 0$ corresponds (by construction) to a two
tori degeneration point of the genus two surface where $\Omega \rightarrow
\mathrm{diag}(\Omega _{11}=\tau _{1},\Omega _{22}=\tau _{2})$. One can then
show that there is a $G$-invariant neighborhood of each degeneration point
on which $F^{\epsilon }$ is invertible (op.cite., Proposition 4.10).

\subsection{The $\protect\rho $-Formalism for Self-Sewing a Torus\label%
{SubSect_rhosewing}}

We may alternatively construct a genus two Riemann surface by self-sewing a
twice-punctured torus. Consider an oriented torus $\mathcal{S}=\mathbb{C}%
/\Lambda $ with local coordinate $z$ for lattice $\Lambda =2\pi i(\mathbb{Z}%
\tau \oplus \mathbb{Z})$ with $\tau \in \mathbb{H}_{1}$. Consider two disks
centred at $z=0$ and $z=w$ with local coordinates $z_{1}=z$ and $z_{2}=z-w$
of radius $r_{a}<\frac{1}{2}D(q)$ for $a=1,2$. Note that $r_{1},r_{2}$ must
also be sufficiently small to ensure that the disks do not intersect on $%
\mathcal{S}$. Introduce a complex parameter ${\rho }$ where $|{\rho }|\leq
r_{1}r_{2}$ and define annular regions $\mathcal{A}_{a}=\{z_{a},|{\rho }|r_{%
\bar{a}}^{-1}\leq \left\vert z_{a}\right\vert \leq r_{a}\}$. We identify $%
\mathcal{A}_{1}$ with $\mathcal{A}_{2}$ as a single region via the sewing
relation
\begin{equation}
z_{1}z_{2}=\rho .  \label{rhosew}
\end{equation}%
The genus two Riemann surface (excluding the degeneration point $\rho =0$)
so constructed is then parameterized by the domain
\begin{equation}
\mathcal{D}^{\rho }=\{(\tau ,w,\rho )\in \mathbb{H}_{1}\times \mathbb{%
C\times C}|\quad |w-\lambda |>2|\rho |^{1/2}>0\text{ for all }\lambda \in
\Lambda \},  \label{Drho}
\end{equation}%
where the first inequality follows from the requirement that the annuli do
not intersect.

In the $\rho $-formalism the genus two period matrix is expressed as a
function of\ $(\tau ,w,\rho )\in \mathcal{D}^{\rho }$ in terms of a
doubly-indexed infinite matrix $R(\tau ,w,\rho )=(R_{ab}(k,l,\tau ,w,\rho ))$
for $k,l\geq 1$ and $a,b\in \{1,2\}$ where \cite{MT2}
\begin{equation}
R(k,l,\tau ,w,\rho )=-\frac{\rho ^{(k+l)/2}}{\sqrt{kl}}\left[
\begin{array}{cc}
D(k,l,\tau ,w) & C(k,l,\tau ) \\
C(k,l,\tau ) & D(l,k,\tau ,w)%
\end{array}%
\right] .  \label{Rdef}
\end{equation}%
Note that $R_{ab}(k,l)=R_{\bar{a}\bar{b}}(l,k)$. Similarly to Theorem \ref%
{Theorem_A1A2} \ we find $I-R$ and $\det (I-R)$ play a central role (where
now $I$ denotes the doubly-indexed identity matrix)

\begin{theorem}
\label{Theorem_R}

(a) (op. cite., Proposition 6) The infinite matrix
\begin{equation}
(I-R)^{-1}=\sum_{n\geq 0}R^{n},  \label{1_minus_R}
\end{equation}%
is convergent for $(\tau ,w,\rho )\in \mathcal{D}^{\rho }$.

(b) (op. cite., Theorem 7) $\det (I-R)$ is non-vanishing and holomorphic for
$(\tau ,w,\rho )\in \mathcal{D}^{\rho }$. $\square $
\end{theorem}

The period matrix $\Omega $ is determined as follows:

\begin{theorem}
\label{Theorem_period_rho} (op. cite., Proposition 11) The $\rho $-formalism
determines a holomorphic map%
\begin{eqnarray}
F^{\rho }:\mathcal{D}^{\rho } &\rightarrow &\mathbb{H}_{2},  \notag \\
(\tau ,w,\rho ) &\mapsto &\Omega (\tau ,w,\rho ),  \label{Frhomap}
\end{eqnarray}%
where $\Omega =\Omega (\tau ,w,\rho )$ is given by
\begin{eqnarray}
2\pi i\Omega _{11} &=&2\pi i\tau -\rho \sigma ((I-R)^{-1}(1,1)),
\label{Om11rho} \\
2\pi i\Omega _{12} &=&w-\rho ^{1/2}\sigma ((b(I-R)^{-1}(1)),  \label{Om12rho}
\\
2\pi i\Omega _{22} &=&\log (-\frac{\rho }{K(\tau ,w)^{2}})-b(I-R)^{-1}\bar{b}%
^{T}.  \label{Om22rho}
\end{eqnarray}%
$K$ is the elliptic prime form, $b=(b_{a}(k,\tau ,w,{\rho }))$ is a
doubly-indexed infinite row vector\footnote{%
Note that $b$ is denoted by $\beta $ in ref. \cite{MT2}.}%
\begin{equation}
b(k,\tau ,w,{\rho })=\frac{\rho ^{k/2}}{\sqrt{k}}(P_{k}(\tau ,w)-E_{k}(\tau
))[-1,(-1)^{k}],  \label{bk}
\end{equation}%
with $\bar{b}_{a}=b_{\bar{a}}$. $(1,1)$ and $(1)$ refer to the $(k,l)=(1,1)$%
, respectively, $(k)=(1)$ entries of an infinite matrix and row vector
respectively. $\sigma (M)$ denotes the sum over the finite indices for a
given $2\times 2$ or $1\times 2$ matrix $M$. $\square $
\end{theorem}

The domain $\mathcal{D}^{\rho }$ is preserved by the Jacobi group $J=SL(2,%
\mathbb{Z})\ltimes \mathbb{Z}^{2}$ where
\begin{eqnarray}
(a,b).(\tau ,w,\rho ) &=&(\tau ,w+2\pi ia\tau +2\pi ib,\rho ),\quad (a,b)\in
\mathbb{Z}^{2}  \label{ab_rho} \\
\gamma _{1}.(\tau ,w,\rho ) &=&(\frac{a_{1}\tau +b_{1}}{c_{1}\tau +d_{1}},%
\frac{w}{c_{1}\tau +d_{1}},\frac{\rho }{(c_{1}\tau +d_{1})^{2}}),\quad
\gamma _{1}\in \Gamma _{1},  \label{gam1_rho}
\end{eqnarray}%
with $\Gamma _{1}=SL(2,\mathbb{Z})$. However, due to the branch structure of
the logarithmic term in (\ref{Om22rho}) the map $F^{\rho }$ is not
equivariant with respect to $J$. (Instead one must pass to a
simply-connected covering space $\mathcal{\hat{D}}^{\rho }$ on which $L=\hat{%
H}\Gamma _{1}$, a split extension of $SL(2,\mathbb{Z})$ by an integer
Heisenberg group $\hat{H}$, acts - see Section 6.3 of ref. \cite{MT2} for
details). However, on restricting to the modular subgroup $\Gamma _{1}$ we
find

\begin{theorem}
\label{TheoremG1equiv} (op. cit., Theorem 11, Corollary 2) $F^{\rho }$ is
equivariant with respect to the action of $\Gamma _{1}$ i.e. there is a
commutative diagram for $\gamma _{1}\in \Gamma _{1}$,
\begin{equation*}
\begin{array}{ccc}
\mathcal{D}^{\rho } & \overset{F^{\rho }}{\rightarrow } & \mathbb{H}_{2} \\
\gamma _{1}\downarrow &  & \downarrow \gamma _{1} \\
\mathcal{D}^{\rho } & \overset{F^{\rho }}{\rightarrow } & \mathbb{H}_{2}%
\end{array}%
,
\end{equation*}%
where the action of $\Gamma _{1}$ on $\mathbb{H}_{2}$ is that of (\ref%
{Gamma1Gamma2}). $\square $
\end{theorem}

The two tori degeneration limit can also be considered in the $\rho $%
-formalism as follows. Define the $\Gamma _{1}$-invariant
\begin{equation}
\chi =-\frac{\rho }{w^{2}}.  \label{chidef}
\end{equation}%
Then one finds that two tori degeneration limit is given by $\rho
,w\rightarrow 0$ for fixed $\chi $ where
\begin{equation}
\Omega \rightarrow \left(
\begin{array}{cc}
\tau & 0 \\
0 & \frac{1}{2\pi i}\log (f(\chi ))%
\end{array}%
\right)  \label{Omrho_degen}
\end{equation}%
with $f(\chi )$ the Catalan series
\begin{equation}
f(\chi )=\frac{1-\sqrt{1-4\chi }}{2\chi }-1=\sum_{n\geq 1}\frac{1}{n}\binom{%
2n}{n+1}\chi ^{n}.  \label{Catalan}
\end{equation}%
(The coefficients $\frac{1}{n}\binom{2n}{n+1}$ are the ubiquitous Catalan
numbers of combinatorics e.g. \cite{St}.)

\medskip
In order to describe the limit (\ref{Omrho_degen}) more precisely, we
introduce the domain
\begin{equation}
\mathcal{D}^{\chi }=\{(\tau ,w,\chi )\in \mathbb{H}_{1}\times \mathbb{%
C\times C}\ |\ (\tau ,w,-w^{2}\chi )\in \mathcal{D}^{\rho },0<|\chi |<\frac{1%
}{4}\},  \label{Dchidomain}
\end{equation}%
and a $\Gamma _{1}$-equivariant holomorphic map
\begin{eqnarray}
F^{\chi }:\mathcal{D}^{\chi } &\rightarrow &\mathbb{H}_{2},  \notag \\
(\tau ,w,\chi ) &\mapsto &\Omega ^{(2)}(\tau ,w,-w^{2}\chi ).  \label{Fchi}
\end{eqnarray}%
Then
\begin{equation*}
\mathcal{D}_{0}^{\chi }=\{(\tau ,0,\chi )\in \mathbb{H}_{1}\times \mathbb{%
C\times C}|0<|\chi |<\frac{1}{4}\},
\end{equation*}%
is the space\ of two-tori degeneration limit points of the domain $\mathcal{D%
}^{\chi }$. We may then compare the two parameterizations on certain $\Gamma
_{1}$-invariant neighborhoods of a two tori degeneration point in both
parameterizations to obtain

\begin{theorem}
\label{Theorem_epsrho_11map} (op. cit., Theorem 12) There exists a 1-1
holomorphic mapping between $\Gamma _{1}$-invariant open domains $\mathcal{I}%
^{\chi }\subset (\mathcal{D}^{\chi }\cup \mathcal{D}_{0}^{\chi })$ and $%
\mathcal{I}^{\epsilon }\subset \mathcal{D}^{\epsilon }$ where $\mathcal{I}%
^{\chi }$ and $\mathcal{I}^{\epsilon }$ are open neighborhoods of a two tori
degeneration point. $\square $
\end{theorem}

\section{Graphical expansions}

\subsection{Rotationless and chequered cycles}

We set up some notation and discuss certain types of labelled graphs. These
arise directly from consideration of the terms that appear in the
expressions for $\Omega _{ij}$\ reviewed in the last Section, and will later
play an important r\^{o}le in the analysis of genus two partition functions
for vertex operator algebras.

Consider a set of independent (non-commuting) variables $x_{i}$ indexed by
the elements of a finite set $I=\{1,...,N\}$. The set of all distinct
monomials $x_{i_{1}}...x_{i_{n}}(n\geq 0)$ may be considered as a basis for
the tensor algebra associated with an $N$ dimensional vector space. Call $n$
the degree of the monomial $x_{i_{1}}...x_{i_{n}}$.

Let $\rho =\rho _{n}$ be the standard cyclic permutation which acts on
monomials of degree $n$ via $\rho :x_{i_{1}}...x_{i_{n}}\mapsto
x_{i_{n}}x_{i_{1}}...x_{i_{n-1}}$. The \emph{rotation group} of a given
monomial $x=x_{i_{1}}...x_{i_{n}}$ is the subgroup of $\langle \rho
_{n}\rangle $ that leaves $x$ invariant. Call $x$ \emph{rotationless} in
case its rotation group is trivial. Let us say that two monomials $x,y$ of
degree $n$ are \emph{equivalent} in case $y=\rho _{n}^{r}(x)$ for some $r\in
Z$, and denote the corresponding equivalence class by $(x)$. We call these
\emph{cycles}. Note that equivalent monomials have the same rotation group,
so we may meaningfully refer to the rotation group of a cycle. In
particular, a \emph{rotationless cycle} is a cycle whose representative
monomials are themselves rotationless. Let $C_{n}$ be the set of
inequivalent cycles of degree $n$.

It is convenient to identify a cycle $(x_{i_{1}}...x_{i_{n}})$ with a \emph{%
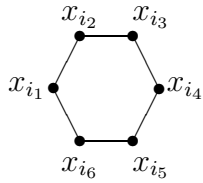
cyclic labelled graph} or \emph{labelled polygon}, that is, a graph with $n$
vertices labelled $x_{i_{1}},...,x_{i_{n}}$ and with edges $%
x_{i_{1}}x_{i_{2}},...,x_{i_{n-1}}x_{i_{n}},x_{i_{n}}x_{i_{1}}$. We will
sometimes afflict the graph with one of the two canonical orientations.

\bigskip

\begin{center}
\begin{picture}(250,80)


\put(100,50){\line(1,2){10}}
\put(90,50){\makebox(0,0){$x_{i_1 }$}}
\put(100,50){\circle*{4}}

\put(110,70){\line(1,0){20}}
\put(110,77){\makebox(0,0){$x_{i_2 }$}}
\put(110,70){\circle*{4}}

\put(130,70){\line(1,-2){10}}
\put(137,77){\makebox(0,0){$x_{i_3}$}}
\put(130,70){\circle*{4}}

\put(140,50){\line(-1,-2){10}}
\put(150,50){\makebox(0,0){$x_{i_4}$}}
\put(140,50){\circle*{4}}

\put(100,50){\line(1,-2){10}}
\put(110,20){\makebox(0,0){$x_{i_6 }$}}
\put(110,30){\circle*{4}}

\put(110,30){\line(1,0){20}}
\put(137,20){\makebox(0,0){$x_{i_5 }$}}
\put(130,30){\circle*{4}}

\end{picture}

{\small Fig. 2 A Cyclic Labelled Graph}
\end{center}

A cycle is rotationless precisely when its graph admits no non-trivial
rotations (a rotation now being an orientation-preserving automorphism of
the graph which preserves labels of nodes).

\bigskip Next we introduce the notion of a \emph{chequered cycle }as a
(clockwise) oriented, labelled polygon $L$ with $2n$ nodes for some integer $%
n\geq 0$, and nodes labelled by arbitrary positive integers. Moreover, edges
carry a label $1$ or $2$ which alternate as one moves around the polygon.

\bigskip

\begin{center}
\begin{picture}(250,80)


\put(100,50){\line(1,2){10}}
\put(90,50){\makebox(0,0){$i_1$}}
\put(100,50){\circle*{4}}
\put(102,63){\makebox(0,0){\scriptsize 1}}
\put(107,63){\vector(1,2){0}}

\put(110,70){\line(1,0){20}}
\put(105,78){\makebox(0,0){$i_2$}}
\put(110,70){\circle*{4}}
\put(120,75){\makebox(0,0){\scriptsize 2}}
\put(122,70){\vector(1,0){0}}

\put(130,70){\line(1,-2){10}}
\put(139,78){\makebox(0,0){$i_3$}}
\put(130,70){\circle*{4}}
\put(139,63){\makebox(0,0){\scriptsize 1}}
\put(137,57){\vector(1,-2){0}}

\put(140,50){\line(-1,-2){10}}
\put(150,50){\makebox(0,0){$i_4$}}
\put(140,50){\circle*{4}}
\put(138,36){\makebox(0,0){\scriptsize 2}}
\put(134,38){\vector(-1,-2){0}}

\put(110,30){\line(1,0){20}}
\put(135,20){\makebox(0,0){$i_5$}}
\put(130,30){\circle*{4}}
\put(120,25){\makebox(0,0){\scriptsize 1}}
\put(117,30){\vector(-1,0){0}}

\put(100,50){\line(1,-2){10}}
\put(105,20){\makebox(0,0){$i_6$}}
\put(110,30){\circle*{4}}
\put(102,36){\makebox(0,0){\scriptsize 2}}
\put(104,42){\vector(-1,2){0}}

\end{picture}

{\small Fig. 3 Chequered Cycle}
\end{center}

We call a node with label $1$ \emph{distinguished} if its abutting edges are
of type $\overset{2}{\longrightarrow }\overset{1}{\bullet }\overset{1}{%
\longrightarrow }$. Set
\begin{eqnarray}  \label{mathcalRdef}
\mathcal{R} &=& \{%
\mbox{isomorphism classes of rotationless chequered cycles
} \}, \\
\mathcal{R}_{21} &=&\{
\mbox{isomorphism classes of rotationless  chequered
cycles }  \notag \\
&&\mbox{ with a  distinguished node}\},  \notag \\
\mathcal{L}_{21} &=&\{\mbox{isomorphism classes of chequered cycles  with a}
\notag \\
&&\mbox{ \emph{unique} distinguished node}\},  \notag
\end{eqnarray}

\medskip

Let $S$ be a commutative ring and $S[t]$ the polynomial ring with
coefficients in $S$. Let $M_{1}$ and $M_{2}$ be infinite matrices with $%
(k,l) $-entries
\begin{equation}
M_{a}(k,l)=t^{k+l}s_{a}(k,l)  \label{eq: abstractAmatrix}
\end{equation}%
for $a=1,2$ and $k,l\geq 1$, where $s_{a}(k,l)\in S$. Given this data, we
define a map, or \emph{weight function},
\begin{equation*}
\omega :\{\mbox{chequered cycles}\}\longrightarrow S[t]
\end{equation*}%
as follows: if $L$ is a chequered cycle then $L$ has edges $E$ labelled as \
$\overset{k}{\bullet }\overset{a}{\longrightarrow }\overset{l}{\bullet }$.
Then set $\omega (E)=M_{a}(k,l)$ and
\begin{equation}
\omega (L)=\prod \omega (E)  \label{eq: eqomega}
\end{equation}%
where the product is taken over all edges of $L$.

\bigskip It is useful to also introduce a variation on the theme of
chequered polygons, namely \emph{\ oriented chequered necklaces}. These are
connected graphs with $n\geq 3$ nodes, $(n-2)$ of which have valency $2$ and
two of which have valency $1$ (these latter are the \emph{end nodes})
together with an orientation, say from left to right. There is also a
degenerate necklace $N_{0}$ with a single node and no edges. As before,
nodes are labelled with arbitrary positive integers and edges are labelled
with an index $1$ or $2$ which alternate along the necklace. For such a
necklace $N$, we define the weight function $\omega (N)$ as a product of
edge weights as in (\ref{eq: eqomega}), with $\omega (N_{0})=1$.

\medskip Among all chequered necklaces there is a distinguished set for
which both end nodes are labelled by $1$. There are four types of such
chequered necklaces, which may be further distinguished by the labels of the
two edges at the extreme left and right. Using the convention (\ref{abar})
we say that the chequered necklace
\begin{equation*}
\underset{1}{\bullet }\overset{\bar{a}}{\longrightarrow }\underset{i}{%
\bullet }\ldots \underset{j}{\bullet }\overset{\bar{b}}{\longrightarrow }%
\underset{1}{\bullet }
\end{equation*}

\begin{center}
{\small {Fig. 4 }}
\end{center}

\noindent is of \emph{type $ab$ }for $a,b\in \{1,2\}$, and set
\begin{eqnarray}  \label{Nabdef}
\mathcal{N}_{ab} &=&\{\mbox{isomorphism classes of oriented chequered}
\notag \\
&&\mbox{ \ necklaces of type}\ ab\},
\end{eqnarray}
\begin{eqnarray}  \label{omegaabdef}
\omega _{ab} &=&\sum_{N\in \mathcal{N}_{ab}}\omega (N).
\end{eqnarray}

\medskip There are other types of graphs that we will need, adapted to the $%
\rho$-formalism. Rather than discussing them here, we delay their
introduction until Subsection 3.3 below.

\subsection{\protect\bigskip Graphical expansions for $\Omega$ in the $%
\protect\epsilon$-formalism}

We now apply the formalism of the previous Subsection to the expressions for
$\Omega _{ij}$ in the $\epsilon $-formalism reviewed in Section 2. The ring $%
S$ is taken to be the product $S_{1}\times S_{2}$ where for $a=1,2$, $S_{a}$
is the ring of quasi-modular forms $\mathbb{C}[E_{2}(\tau _{a}),E_{4}(\tau
_{a}),E_{6}(\tau _{a})]$, and $t = \epsilon^{1/2}$. The matrices $M_a$ are
taken to be the $A_a$ defined in (\ref{Akldef}). Thus
\begin{equation}
s_{a}(k,l)=\frac{C(k,l,\tau _{a})}{\sqrt{kl}},  \label{eq: subst}
\end{equation}%
and for the edge $E$ labelled as \ $\overset{k}{\bullet }\overset{a}{%
\longrightarrow }\overset{l}{\bullet}$ we have
\begin{eqnarray}  \label{omegaval}
\omega (E)=A_{a}(k,l).
\end{eqnarray}

\medskip We can now state

\begin{proposition}
\label{Prop_Om12_R21expansion}In the $\epsilon -$formalism we have
\begin{equation}
\Omega _{12}=-\frac{\epsilon }{2\pi i}\prod_{L\in \mathcal{R}_{21}}(1-\omega
(L))^{-1}.  \label{eq: omega_12prod}
\end{equation}
\end{proposition}
\noindent
\textbf{Proof.} Beyond the intrinsic interest of this product formula, our
main use of it will be to provide an alternate proof of Theorem \ref%
{Theorem_Z2_G} below. We therefore relegate the proof to Proposition \ref%
{Prop_R21expansion} in an Appendix. \ \ \ \ \ $\Box$

\bigskip We will need some further identities of this nature. Recalling the
notation (\ref{omegaabdef}),

\begin{proposition}
\label{Propepsperiodgraph} (\cite{MT2}, Proposition 4) In the $\epsilon $%
-formalism we have for $a=1,2$ that
\begin{eqnarray*}
\Omega _{aa} &=&\tau _{a}+\frac{\epsilon }{2\pi i}\omega _{aa}, \\
\Omega _{a\bar{a}} &=&-\frac{\epsilon }{2\pi i}\omega _{a\bar{a}}.\ \ \ \ \
\ \ \square \quad
\end{eqnarray*}
\end{proposition}

\begin{remark}
\label{Om12_graphs}Proposition \ref{Prop_Om12_R21expansion} implies that
\begin{eqnarray*}
\omega _{12}=\prod_{L\in \mathcal{R}_{21}}(1-\omega (L))^{-1}.
\end{eqnarray*}
\end{remark}

\subsection{\protect\bigskip Graphical expansions for $\Omega$ in the $%
\protect\rho$-formalism}

We turn next to the expressions for $\Omega $ in the $\rho $-formalism
reviewed in Section 2. In this case it is natural to introduce \emph{%
doubly-indexed cycles.} These are (clockwise) oriented, labelled polygons $L$
with $n$ nodes for some integer $n\geq 1$, nodes being labelled by a pair of
integers $k,a$ where $k\geq 1$ and $a\in \{1,2\}$. Thus, a typical
doubly-indexed cycle looks as follows:

\bigskip

\begin{center}
\begin{picture}(250,80)


\put(100,50){\line(1,2){10}}
\put(82,50){\makebox(0,0){$k_1,a_1$}}
\put(100,50){\circle*{4}}
\put(107,63){\vector(1,2){0}}

\put(110,70){\line(1,0){20}}
\put(100,78){\makebox(0,0){$k_2,a_2$}}
\put(110,70){\circle*{4}}
\put(122,70){\vector(1,0){0}}

\put(130,70){\line(1,-2){10}}
\put(145,78){\makebox(0,0){$k_3,a_3$}}
\put(130,70){\circle*{4}}
\put(137,57){\vector(1,-2){0}}

\put(140,50){\line(-1,-2){10}}
\put(160,50){\makebox(0,0){$k_4,a_4$}}
\put(140,50){\circle*{4}}
\put(134,38){\vector(-1,-2){0}}

\put(110,30){\line(1,0){20}}
\put(145,20){\makebox(0,0){$k_5,a_5$}}
\put(130,30){\circle*{4}}
\put(117,30){\vector(-1,0){0}}

\put(100,50){\line(1,-2){10}}
\put(100,20){\makebox(0,0){$k_6,a_6$}}
\put(110,30){\circle*{4}}
\put(104,42){\vector(-1,2){0}}

\end{picture}

{\small Fig. 5 Doubly-Indexed Cycle}
\end{center}

Next we define a weight function $\omega $ with values in the ring of
elliptic functions and quasi-modular forms $\mathbb{C}[P_{2}(\tau
,w),P_{3}(\tau ,w),E_{2}(\tau ),E_{4}(\tau ),E_{6}(\tau )]$ as follows: if $%
L $ is a doubly-indexed cycle then $L$ has edges $E$ labelled as \ $\overset{%
k,a}{\bullet }\rightarrow \overset{l,b}{\bullet }$, and we set
\begin{equation*}
\omega (E)=R_{ab}(k,l,\tau ,w,{\rho }),
\end{equation*}%
with $R_{ab}(k,l)$ as in (\ref{Rdef}) and
\begin{equation*}
\omega (L)=\prod \omega (E).
\end{equation*}%
As usual, the product is taken over all edges of $L$. It is straightforward
(but of less value for our purposes) to obtain an analog of Proposition \ref%
{Prop_R21expansion}.

\medskip We also introduce \emph{doubly-indexed necklaces }$\mathcal{N}%
=\{N\} $. These are connected graphs with $n\geq 2$ nodes, $(n-2)$ of which
have valency $2$ and two of which have valency $1$ together with an
orientation, say from left to right, on the edges. In this case, each vertex
carries two integer labels $k,a$ with $k\geq 1$ and $a\in \{1,2\}$. We
define the degenerate necklace $N_{0}$ to be a single node with no edges,
and set $\omega (N_{0})=1$.

\medskip We define necklaces with distinguished end nodes labelled $k,a;l,b$
as follows:
\begin{equation*}
\underset{k,a}{\bullet }\longrightarrow \underset{k_{1},a_{1}}{\bullet }%
\ldots \underset{k_{2},a_{2}}{\bullet }\longrightarrow \underset{l,b}{%
\bullet }\hspace{10mm}\mbox{(type $k,a;l,b$)}
\end{equation*}%
and set
\begin{equation*}
\mathcal{N}(k,a;l,b)=\{\mbox{isomorphism classes of
necklaces of type}\ k,a;l,b\}.
\end{equation*}%
It is convenient to define
\begin{eqnarray}
\omega _{11} &=&\sum_{a_{1},a_{2}=1,2}\sum_{N\in \mathcal{N}%
(1,a_{1};1,a_{2})}\omega (N),  \notag \\
\omega _{b1} &=&\sum_{a_{1},a_{2}=1,2}\sum_{k\geq 1}b_{a_{1}}(k)\sum_{N\in
\mathcal{N}(k,a_{1};1,a_{2})}\omega (N),  \label{om_rho_weights} \\
\omega _{b\bar{b}} &=&\sum_{a_{1},a_{2}=1,2}\sum_{k,l\geq 1}b_{a_{1}}(k)\bar{%
b}_{a_{2}}(l)\sum_{N\in \mathcal{N}(k,a_{1};l,a_{2})}\omega (N).  \notag
\end{eqnarray}%
Then we find

\begin{proposition}
\label{Prop_rhoperiod_graph} (\cite{MT2}, Proposition 12) In the $\rho $%
-formalism we have
\begin{eqnarray*}
2\pi i\Omega _{11} &=&2\pi i\tau -\rho \omega _{11}, \\
2\pi i\Omega _{12} &=&w-\rho ^{1/2}\omega _{b1}, \\
2\pi i\Omega _{22} &=&\log (-\frac{\rho }{K(\tau ,w)^{2}})-\omega _{b\bar{b}%
}.\ \ \ \ \ \ \ \ \square
\end{eqnarray*}
\end{proposition}

\section{Vertex operator algebras and the Li-Zamalodchikov metric}

\subsection{Vertex operator algebras}

We review some relevant aspects of vertex operator algebras (\cite{FHL},\cite%
{FLM}, \cite{Ka}, \cite{LL}, \cite{MN}). A vertex operator algebra (VOA) is
a quadruple $(V,Y,\mathbf{1},\omega )$ consisting of a $\mathbb{Z}$-graded
complex vector space $V=\bigoplus_{n\in \mathbb{Z}}V_{n}$, a linear map $%
Y:V\rightarrow (\mathrm{End}V)[[z,z^{-1}]]$, for formal parameter $z$, and a
pair of distinguished vectors (states): the vacuum $\mathbf{1}\in V_{0}$ ,
and the conformal vector $\omega \in V_{2}$. For each state $v\in V$ the
image under the $Y$ map is the vertex operator

\begin{equation}
Y(v,z)=\sum_{n\in \mathbb{Z}}v(n)z^{-n-1},  \label{Ydefn}
\end{equation}%
with modes $v(n)\in \mathrm{End}V$ where $\mathrm{Res}_{z=0}z^{-1}Y(v,z)%
\mathbf{1}=v(-1)\mathbf{1}=v$. Vertex operators satisfy the Jacobi identity
or equivalently, operator locality or Borcherds's identity for the modes
(loc. cit.).

The vertex operator for the conformal vector $\omega $ is defined as

\begin{equation*}
Y(w,z)=\sum_{n\in \mathbb{Z}}L(n)z^{-n-2}.
\end{equation*}%
The modes $L(n)$ satisfy the Virasoro algebra of central charge $c$:

\begin{equation*}
\lbrack L(m),L(n)]=(m-n)L(m+n)+(m^{3}-m)\frac{c}{12}\delta _{m,-n}.
\end{equation*}
We define the homogeneous space of weight $k$ to be $V_{k}=\{v\in
V|L(0)v=kv\}$ where we write $wt(v)=k$ for $v$ in $V_{k}$ . Then as an
operator on $V$ we have

\begin{equation*}
v(n):V_{m}\rightarrow V_{m+k-n-1}.
\end{equation*}%
In particular, the \textit{zero mode} $o(v)=v(wt(v)-1)$ is a linear operator
on $V_{m}$. A state $v$ is said to be \textit{quasi-primary} if $L(1)v=0$
and \textit{primary} if additionally $L(2)v=0$.

The subalgebra $\{L(-1),L(0),L(1)\}$ generates a natural action on vertex
operators associated with $SL(2,\mathbb{C})$ M\"{o}bius transformations on $%
z $ ( \cite{B1}, \cite{DGM}, \cite{FHL}, \cite{Ka}). In particular, we note
the inversion $z\mapsto 1/z$ for which
\begin{equation}
Y(v,z)\mapsto Y^{\dagger }(v,z)=Y(e^{zL(1)}(-\frac{1}{z^{2}})^{L(0)}v,\frac{1%
}{z}).  \label{eq: adj op}
\end{equation}%
$Y^{\dagger }(v,z)$ is the \emph{adjoint} vertex operator \cite{FHL}. Under
the dilatation $z\mapsto az$ we have
\begin{equation}
Y(v,z)\mapsto a^{L(0)}Y(v,z)a^{-L(0)}=Y(a^{L(0)}v,az).  \label{Y_D}
\end{equation}%
We also note (\cite{BPZ},\cite{Z2}) that under a general origin-preserving
conformal map $z\mapsto w=\phi (z)$,
\begin{equation}
Y(v,z)\mapsto Y((\phi ^{\prime }(z))^{L(0)}v,w),  \label{Y_phi}
\end{equation}%
for any primary vector $v$.

\bigskip We consider some particular VOAs, namely Heisenberg free boson and
lattice VOAs. Consider an $l$-dimensional complex vector space (i.e.,
abelian Lie algebra) $\mathfrak{H}$ equipped with a non-degenerate,
symmetric, bilinear form $(\ ,)$ and a distinguished orthonormal basis $%
a_{1},a_{2},...a_{l}$. The corresponding affine Lie algebra is the
Heisenberg Lie algebra $\mathfrak{\hat{H}}=\mathfrak{H}\otimes \mathbb{C}%
[t,t^{-1}]\oplus \mathbb{C}k$ with brackets $[k,\mathfrak{\hat{H}}]=0$ and

\begin{equation}
\lbrack a_{i}\otimes t^{m},a_{j}\otimes t^{n}]=m\delta _{i,j}\delta _{m,-n}k.
\label{Fockbracket}
\end{equation}%
Corresponding to an element $\lambda $ in the dual space $\mathfrak{H}^{\ast
}$ we consider the Fock space defined by the induced (Verma) module
\begin{equation*}
M^{(\lambda )}=U(\mathfrak{\hat{H}})\otimes _{U(\mathfrak{H}\otimes \mathbb{C%
}[t]\oplus \mathbb{C}k)}\mathbb{C},
\end{equation*}%
where $\mathbb{C}$ is the $1$-dimensional space annihilated by $\mathfrak{H}%
\otimes t\mathbb{C}[t]$ and on which $k$ acts as the identity and $\mathfrak{%
H}\otimes t^{0}$ via the character $\lambda $; $U$ denotes the universal
enveloping algebra. There is a canonical identification of linear spaces

\begin{equation*}
M^{(\lambda )}=S(\mathfrak{H}\otimes t^{-1}\mathbb{C}[t^{-1}]),
\end{equation*}%
where $S$ denotes the (graded) symmetric algebra. The Heisenberg free boson
VOA $M^{l}$ corresponds to the case $\lambda =0$. The Fock states

\begin{equation}
v=a_{1}(-1)^{e_{1}}.a_{1}(-2)^{e_{2}}....a_{1}(-n)^{e_{n}}....a_{l}(-1)^{f_{1}}.a_{l}(-2)^{f_{2}}...a_{l}(-p)^{f_{p}}.%
\mathbf{1,}  \label{Fockstate}
\end{equation}%
for non-negative integers $e_{i},...,f_{j}$ form a\textit{\ }basis of $M^{l}$%
. The vacuum $\mathbf{1}$ is canonically identified with the identity of $%
M_{0}=\mathbb{C}$, while the weight 1 subspace $M_{1}$ may be naturally
identified with $\mathfrak{H}$. $M^{l}$ is a simple VOA of central charge $l$%
.

\bigskip Next we consider the case of a lattice vertex operator algebra $%
V_{L}$ associated to a positive-definite even lattice $L$ (cf. \cite{B1},
\cite{FLM}). Thus $L$ is a free abelian group of rank $l$ equipped with a
positive definite, integral bilinear form $(\ ,):L\otimes L\rightarrow
\mathbb{Z}$ such that $(\alpha ,\alpha )$ is even for $\alpha \in L$. Let $%
\mathfrak{H}$ be the space $\mathbb{C}\otimes _{\mathbb{Z}}L$ equipped with
the $\mathbb{C}$-linear extension of $(\ ,)$ to $\mathfrak{H}\otimes
\mathfrak{H}$ and let $M^{l}$ be the corresponding Heisenberg VOA. The Fock
space of the lattice theory may be described by the linear space
\begin{equation}
V_{L}=M^{l}\otimes \mathbb{C}[L]=\sum_{\alpha \in L}M^{l}\otimes e^{\alpha },
\label{VLdefn}
\end{equation}%
where $\mathbb{C}[L]$ denotes the group algebra of $L$ with canonical basis $%
e^{\alpha }$, $\alpha \in L$. $M^{l}$ may be identified with the subspace $%
M^{l}\otimes e^{0}$ of $V_{L}$, in which case $M^{l}$ is a subVOA of $V_{L}$
and the rightmost equation of (\ref{VLdefn}) then displays the decomposition
of $V_{L}$ into irreducible $M^{l}$-modules. $V_{L}$ is a simple VOA of
central charge $l$. Each $\mathbf{1}\otimes e^{\alpha }\in V_{L}$ is a
primary state of weight $\frac{1}{2}(\alpha ,\alpha )$ with vertex operator
(loc. cit.)
\begin{eqnarray}
Y(\mathbf{1}\otimes e^{\alpha },z) &=&Y_{-}(\mathbf{1}\otimes e^{\alpha
},z)Y_{+}(\mathbf{1}\otimes e^{\alpha },z)e^{\alpha }z^{\alpha },  \notag \\
Y_{\pm }(\mathbf{1}\otimes e^{\alpha },z) &=&\exp (\mp \sum_{n>0}\frac{%
\alpha (\pm n)}{n}z^{\mp n}).  \label{Yealpha}
\end{eqnarray}%
The operators $e^{\alpha }\in \mathbb{C}[L]$ obey
\begin{equation}
e^{\alpha }e^{\beta }=\epsilon (\alpha ,\beta )e^{\alpha +\beta }
\label{eps_cocycle}
\end{equation}%
for $2$-cocycle $\epsilon (\alpha ,\beta )$ satisfying $\epsilon (\alpha
,\beta )\epsilon (\beta ,\alpha )=(-1)^{(\alpha ,\beta )}$.

\subsection{The Li-Zamolodchikov metric}

A bilinear form $\langle \ ,\rangle :V\times V{\longrightarrow }\mathbb{C}$
is called \emph{invariant} in case the following identity holds for all $%
a,b,c\in V$ (\cite{FHL}):
\begin{equation}
\langle Y(a,z)b,c\rangle =\langle b,Y^{\dagger }(a,z)c\rangle,
\label{eq: inv bil form}
\end{equation}%
with $Y^{\dagger }(a,z)$ the adjoint operator (\ref{eq: adj op}).

\begin{remark}
\label{Rem_Zam}Note that
\begin{eqnarray}
\langle u,v\rangle &=&\mathrm{Res}_{w=0}w^{-1}\mathrm{Res}%
_{z=0}z^{-1}\langle Y(u,w)\mathbf{1},Y(v,z)\mathbf{1}\rangle  \notag \\
&=&\mathrm{Res}_{w=0}w^{-1}\mathrm{Res}_{z=0}z^{-1}\langle \mathbf{1}%
,Y^{\dagger }(u,w)Y(v,z)\mathbf{1}\rangle  \notag \\
&=&``\langle \mathbf{1},Y(u,z=\infty )Y(v,z=0)\mathbf{1}\rangle
\textquotedblright ,  \label{Zam}
\end{eqnarray}%
with $w=1/z$, following (\ref{eq: adj op}). Thus the invariant bilinear form
is equivalent to what is known as the (chiral) Zamolodchikov metric in
Conformal Field Theory (\cite{BPZ}, \cite{P}).
\end{remark}

Generally a VOA may have no non-zero invariant bilinear form, but a result
of Li \cite{Li} guarantees that if $V_{0}$ is spanned by the vacuum vector $%
\mathbf{1}$, and $V$ is self-dual in the sense that $V$ is isomorphic to the
contragredient module $V^{\prime }$ as a $V$-module, then $V$ has a unique
non-zero invariant bilinear form up to scalar. Note that $\langle \ ,\rangle
$ is necessarily symmetric by a theorem of \cite{FHL}. Furthermore if $V$ is
simple then such a form is necessarily non-degenerate. All of the VOAs that
occur in this paper satisfy these conditions, so that if we normalize so
that $\langle \mathbf{1},\mathbf{1}\rangle =1$ then $\langle \ ,\rangle $ is
unique. We refer to this particular bilinear form as the \emph{%
Li-Zamolodchikov metric} on $V$, or LiZ-metric for short.

\begin{remark}
\label{Rem_tensorVOA}Uniqueness entails that the LiZ-metric on the tensor
product $V_{1}\otimes V_{2}$ of a pair of simple VOAs satisfying the
appropriate conditions is just the tensor product of the LiZ metrics on $%
V_{1}$ and $V_{2}$.
\end{remark}

If $a$ is a homogeneous, quasi-primary state, the component form of (\ref%
{eq: inv bil form}) reads
\begin{equation}
\langle a(n)u,v\rangle =(-1)^{wt(a)}\langle u,a(2wt(a)-n-2)v\rangle .
\label{eq: qp bil form}
\end{equation}%
In particular, since the conformal vector $\omega $ is quasi-primary of
weight $2$ we may take $\omega $ in place of $a$ in (\ref{eq: qp bil form})
and obtain
\begin{equation}
\langle L(n)u,v\rangle =\langle u,L(-n)v\rangle .
\label{eq: conf vec bil form}
\end{equation}%
The case $n=0$ of (\ref{eq: conf vec bil form}) shows that the homogeneous
spaces $V_{n},V_{m}$ are orthogonal if $n\not=m$. Taking $u=\mathbf{1}$ and
using $a=a(-1)\mathbf{1}$ in (\ref{eq: qp bil form}) yields
\begin{equation}
\langle a,v\rangle =(-1)^{wt(a)}\langle \mathbf{1},a(2wt(a)-1)v\rangle ,
\label{eq: liz form}
\end{equation}%
for $a$ quasi-primary, and this affords a practical way to compute the
LiZ-metric.

\bigskip Consider the rank one Heisenberg (free boson) VOA $M = M^1$
generated by a weight one state $a$ with $(a,a)=1$. Then $\langle a,a\rangle
=-\langle \mathbf{1},a(1)a(-1)\mathbf{1}\rangle =-1$. Using (\ref%
{Fockbracket}), it is straightforward to verify that in general the Fock
basis consisting of vectors of the form
\begin{equation}
v=a(-1)^{e_{1}}...a(-p)^{e_{p}}.\mathbf{1,}  \label{eq: fock basis}
\end{equation}%
for non-negative integers $\{e_{i}\}$ is orthogonal with respect to the
LiZ-metric, and that
\begin{equation}
\langle v,v\rangle =\prod_{1\leq i\leq p}(-i)^{e_{i}}e_{i}!.
\label{eq: inner prod}
\end{equation}%
This result generalizes in an obvious way for a rank $l$ free boson VOA $%
M^{l}$ with Fock basis (\ref{Fockstate}) following Remark \ref{Rem_tensorVOA}%
.

\medskip
We consider next the lattice vertex operator algebra $V_{L}$ for a
positive-definite even lattice $L$. We take as our Fock basis the states $%
\{v\otimes e^{\alpha }\}$ where $v$ is as in (\ref{Fockstate}) and $\alpha $
ranges over the elements of $L$.

\begin{lemma}
\label{Lemma_LiZ_lattice} If $u,v\in M^{l}$ and $\alpha ,\beta \in L$, then
\begin{eqnarray*}
\langle u\otimes e^{\alpha },v\otimes e^{\beta }\rangle &=&\langle
u,v\rangle \langle \mathbf{1}\otimes e^{\alpha },\mathbf{1}\otimes e^{\beta
}\rangle \\
&=&(-1)^{\frac{1}{2}(\alpha ,\alpha )}\epsilon (\alpha ,-\alpha )\langle
u,v\rangle \delta _{\alpha ,-\beta }.
\end{eqnarray*}
\end{lemma}

\noindent \textbf{Proof.} It follows by successive applications of (\ref{eq:
qp bil form}) that the first equality in the Lemma is true, and that it is
therefore enough to prove it in the case that $u=v=\mathbf{1}$. We identify
the primary vector $\mathbf{1}\otimes e^{\alpha }$ with $e^{\alpha }$ in the
following. Apply (\ref{eq: liz form}) in this case to see that $\langle
e^{\alpha },e^{\beta }\rangle $ is given by$\ $
\begin{align*}
& (-1)^{\frac{1}{2}(\alpha ,\alpha )}\langle \mathbf{1},e^{\alpha
}((a,a)-1)e^{\beta }\rangle \\
& =(-1)^{\frac{1}{2}(\alpha ,\alpha )}\mathrm{Res}_{z=0}z^{(a,a)-1}\langle
\mathbf{1,}Y(e^{\alpha },z)e^{\beta }\rangle \\
& =(-1)^{\frac{1}{2}(\alpha ,\alpha )}\epsilon (\alpha ,\beta )\mathrm{Res}%
_{z=0}z^{(\alpha ,\beta )+(a,a)-1}\langle \mathbf{1,}Y_{-}(\mathbf{1}\otimes
e^{\alpha },z).\mathbf{1}\otimes e^{\alpha +\beta }\rangle .
\end{align*}%
Unless $\alpha +\beta =0$, all states to the left inside the bracket $%
\langle \ ,\rangle $ on the previous line have positive weight, hence are
orthogonal to $\mathbf{1}$. So $\langle e^{\alpha },e^{\beta }\rangle =0$ if
$\alpha +\beta \not=0$. In the contrary case, the exponential operator
acting on the vacuum yields just the vacuum itself among weight zero states,
and we get $\langle e^{\alpha },e^{-\alpha }\rangle =(-1)^{\frac{1}{2}%
(\alpha ,\alpha )}\epsilon (\alpha ,-\alpha )$ in this case. This completes
the proof of the Lemma. \ \ \ \ \ \ \ $\square $

\begin{corollary}
\label{Corollary_cocycle_choice}We may choose the cocycle so that $\epsilon
(\alpha ,-\alpha )=(-1)^{\frac{1}{2}(\alpha ,\alpha )}$ (see Appendix). In
this case, we have
\begin{equation}
\langle u\otimes e^{\alpha },v\otimes e^{\beta }\rangle =\langle u,v\rangle
\delta _{\alpha ,-\beta }.  \label{LiZ_lattice}
\end{equation}
\end{corollary}

\section{Partition and $n$-point functions for vertex operator algebras on a
Riemann Surface}

In this section we consider the partition and $n$-point functions for a VOA
on a Riemann surface of genus zero, one or two. Our definitions are based on
sewing schemes for the given Riemann surface in terms of one or more
surfaces of lower genus and are motivated by ideas in conformal field theory
especially   \cite{FS}, \cite{So1} and \cite{P}. We assume throughout
that $V$ has a non-degenerate LiZ metric $\langle \ ,\rangle $. Then for any
$V$ basis $\{u^{(a)}\}$, we may define the \emph{dual basis }$\{\bar{u}%
^{(a)}\}$ with respect to the LiZ metric where
\begin{equation}
\langle u^{(a)},\bar{u}^{(b)}\rangle =\delta _{ab}\text{.}  \label{LiZdual}
\end{equation}

\subsection{Genus zero}

We begin with the definition of the genus zero $n$-point function given by:%
\begin{equation}
Z_{V}^{(0)}(v_{1},z_{1};\ldots v_{n},z_{n})=\langle \mathbf{1}%
,Y(v_{1},z_{1})\ldots Y(v_{n},z_{n})\mathbf{1}\rangle ,  \label{Zzero}
\end{equation}%
for $v_{1},\ldots v_{n}\in V$. In particular, the genus zero partition (or $%
0 $-point) function is $Z_{V}^{(0)}=\langle \mathbf{1},\mathbf{1}\rangle =1$%
. \ The genus zero $n$-point function is a rational function of $%
z_{1},\ldots z_{n}$ (with possible poles at $z_{i}=0$ and $z_{i}=z_{j},i\neq
j$). Thus we may consider $z_{1},\ldots z_{n}\in \mathbb{C}\cup \{\infty \}$%
, the Riemann sphere, with $Z_{V}^{(0)}(v_{1},z_{1};\ldots ;v_{n},z_{n})$
evaluated for $\left\vert z_{1}\right\vert >\left\vert z_{2}\right\vert
>\ldots >\left\vert z_{n}\right\vert $ (e.g. \cite{FHL}, \cite{Z2}, \cite{GG}%
). \ The $n$-point function has a canonical geometric interpretation for
primary vectors $v_{i}$ of $L(0)$ weight $wt(v_{i})$. Then $%
Z_{V}^{(0)}(v_{1},z_{1};\ldots ;v_{n},z_{n})$ parameterizes a global
meromorphic differential form on the Riemann sphere,
\begin{equation}
\mathcal{F}_{V}^{(0)}(v_{1},\ldots v_{n})=Z_{V}^{(0)}(v_{1},z_{1};\ldots
;v_{n},z_{n})\prod_{1\leq i\leq n}(dz_{i})^{wt(v_{i})}.  \label{Zzeroform}
\end{equation}%
It follows from (\ref{Y_phi}) that $\mathcal{F}_{V}^{(0)}$ is conformally
invariant.\ This interpretation is the starting point of various
algebraic-geometric approaches to $n$-point functions (\emph{apart }from the
partition or $0$-point function) at higher genera (e.g. \cite{TUY}, \cite{Z2}%
).

\medskip It is instructive to consider $\mathcal{F}_{V}^{(0)}$ in the
context of a trivial sewing of two Riemann spheres parameterized by $z_{1}$
and $z_{2}$ to form another Riemann sphere as follows. For $r_{a}>0,a=1,2$,
and a complex parameter $\epsilon $ satisfying $|\epsilon |\leq r_{1}r_{2}$,
excise the open disks $\left\vert z_{a}\right\vert <|\epsilon |r_{\bar{a}%
}^{-1}$ (recall convention (\ref{abar})) and identify the annular regions $%
r_{a}\geq \left\vert z_{a}\right\vert \geq |\epsilon |r_{\bar{a}}^{-1}$ via
the sewing relation
\begin{equation}
z_{1}z_{2}=\epsilon .  \label{spheresew}
\end{equation}%
Consider $Z_{V}^{(0)}(v_{1},x_{1};\ldots v_{n},x_{n})$ for quasi-primary $%
v_{i}$ with $r_{1}\geq |x_{i}|\geq |\epsilon |r_{2}^{-1}$ and let $%
y_{i}=\epsilon /x_{i}$. Then for $0\leq k\leq n-1$ we find from (\ref%
{LiZdual}) that
\begin{eqnarray}
&&\hspace{1cm}Y(v_{k+1},x_{k+1})\ldots Y(v_{n},x_{n})\mathbf{1}=  \notag \\
&&\sum_{r\geq 0}\sum_{u\in V_{r}}\langle \bar{u},Y(v_{k+1},x_{k+1})\ldots
Y(v_{n},x_{n})\mathbf{1}\rangle u,  \label{Yxkplus_u_expansion}
\end{eqnarray}%
where the inner sum is taken over any basis for $V_{r}$. Thus
\begin{eqnarray*}
&&\hspace{3.5cm}Z_{V}^{(0)}(v_{1},x_{1};\ldots v_{n},x_{n})= \\
&&\sum_{r\geq 0}\sum_{u\in V_{r}}\langle \mathbf{1},Y(v_{1},x_{1})\ldots
Y(v_{k},x_{k})u\rangle \langle \bar{u},Y(v_{k+1},x_{k+1})\ldots
Y(v_{n},x_{n})\mathbf{1}\rangle .
\end{eqnarray*}%
But
\begin{equation*}
\langle \mathbf{1},Y(v_{1},x_{1})\ldots Y(v_{k},x_{k})u\rangle =\mathrm{Res}%
_{z_{1}=0}z_{1}^{-1}Z_{V}^{(0)}(v_{1},x_{1};\ldots v_{k},x_{k};u,z_{1}),
\end{equation*}%
and
\begin{eqnarray*}
&&\langle \bar{u},Y(v_{k+1},x_{k+1})\ldots Y(v_{n},x_{n})\mathbf{1}\rangle \\
&=&\langle \mathbf{1},Y^{\dagger }(v_{n},x_{n})\ldots Y^{\dagger
}(v_{k+1},x_{k+1})\bar{u}\rangle \\
&=&\langle \mathbf{1},\epsilon ^{L(0)}Y^{\dagger }(v_{n},x_{n})\epsilon
^{-L(0)}\ldots \epsilon ^{L(0)}Y^{\dagger }(v_{k+1},x_{k+1})\epsilon
^{-L(0)}\epsilon ^{L(0)}\bar{u}\rangle \\
&=&\epsilon ^{r}\mathrm{Res}_{z_{2}=0}z_{2}^{-1}Z_{V}^{(0)}(v_{n},y_{n};%
\ldots v_{k+1},y_{k+1};\bar{u},z_{2})\prod_{k+1\leq i\leq n}(-\frac{\epsilon
}{x_{i}^{2}})^{wt(v_{i})}.
\end{eqnarray*}%
The last equation holds since for quasiprimary states $v_{i}$, the M\"{o}%
bius transformation $x\mapsto y=\epsilon /x$ induces
\begin{equation*}
Y(v_{i},x_{i})\mapsto \epsilon ^{L(0)}Y^{\dagger }(v_{i},x_{i})\epsilon
^{-L(0)}=(-\frac{\epsilon }{x_{i}^{2}})^{wt(v_{i})}Y(v_{i},y_{i}).
\end{equation*}%
Thus we find

\begin{proposition}
\label{Prop_Z0sew} For quasiprimary states $v_{i}$ with the sewing scheme (%
\ref{spheresew}), we have
\begin{gather*}
\mathcal{F}_{V}^{(0)}(v_{1},\ldots ,v_{n})= \\
\sum_{r\geq 0}\epsilon ^{r}\sum_{u\in V_{r}}\mathrm{Res}%
_{z_{1}=0}z_{1}^{-1}Z_{V}^{(0)}(v_{1},x_{1};\ldots v_{k},x_{k};u,z_{1}) \\
\mathrm{Res}_{z_{2}=0}z_{2}^{-1}Z_{V}^{(0)}(v_{n},y_{n};\ldots
v_{k+1},y_{k+1};\bar{u},z_{2})\prod_{1\leq i\leq
k}(dx_{i})^{wt(v_{i})}\prod_{k+1\leq i\leq n}(dy_{i})^{wt(v_{i})}.\quad
\square
\end{gather*}
\end{proposition}

\subsection{Genus one}

\subsubsection{Self-sewing the Riemann sphere}

We now consider genus one $n$-point functions defined in terms of
self-sewings of a Riemann sphere. We first consider the case where punctures
are located at the origin and the point at infinity. Choose local
coordinates $z_{1}=z$ in the neighborhood of the origin and $%
z_{2}=1/z^{\prime }$ for $z^{\prime }$ in the neighborhood of the point at
infinity. For $a=1,2$ and $r_{a}>0$, identify the annular regions $|q|r_{%
\bar{a}}^{-1}\leq \left\vert z_{a}\right\vert \leq r_{a}$ for complex $q$
satisfying $|q|\leq r_{1}r_{2}$ via the sewing relation $z_{1}z_{2}=q$ i.e. $%
z=qz^{\prime }$. Then it is straightforward to show that the annuli do not
intersect for $|q|<1$, and that $q=\exp (2\pi i\tau )$ where $\tau $ is the
torus modular parameter (e.g. \cite{MT2}, Proposition 8).

\medskip We define the genus one partition function by
\begin{eqnarray}
&&\hspace{3cm}Z_{V}^{(1)}(q)=Z_{V}^{(1)}(\tau )=  \notag \\
&&\sum_{n\geq 0}q^{n}\sum_{u\in V_{n}}\mathrm{Res}_{z_{2}=0}z_{2}^{-1}%
\mathrm{Res}_{z_{1}=0}z_{1}^{-1}\langle \mathbf{1},Y^{\dag }(u,z_{2})Y(\bar{u%
},z_{1})\mathbf{1}\rangle ,  \label{Z1def}
\end{eqnarray}%
where the inner sum is taken over any basis for $V_{n}$. Note that the genus
zero partition function is recovered in the Riemann sphere degeneration
limit $q\rightarrow 0$ where $Z_{V}^{(1)}(\tau )\rightarrow Z_{V}^{(0)}=1$.
\ From (\ref{Zam}) and (\ref{LiZdual}) it follows that
\begin{equation}
Z_{V}^{(1)}(\tau )=\sum_{n\geq 0}\dim V_{n}q^{n}=Tr_{V}(q^{L(0)}),
\label{Z_tau}
\end{equation}%
the standard graded trace definition (excluding, for the present, the extra $%
q^{-c/24}$ factor). The genus one $n$-point function is similarly given by%
\begin{gather*}
\sum_{r\geq 0}q^{r}\sum_{u\in V_{r}}\mathrm{Res}_{z_{2}=0}z_{2}^{-1}\mathrm{%
Res}_{z_{1}=0}z_{1}^{-1}\langle \mathbf{1},Y^{\dag
}(u,z_{2})Y(v_{1},x_{1})\ldots Y(v_{n},x_{n})Y(\bar{u},z_{1})\mathbf{1}%
\rangle \\
=Tr_{V}(Y(v_{1},x_{1})\ldots Y(v_{n},x_{n})q^{L(0)}).
\end{gather*}%
\qquad It is natural to consider the conformal map $x=q_{z}\equiv \exp (z)$
in order to describe the elliptic properties of the $n$-point function \cite%
{Z1}. Since from (\ref{Y_phi}), for a primary state $v$, $Y(v,w)\rightarrow
Y(q_{z}^{L(0)}v,q_{z})$ under this conformal map, we are led to the
following definition of the genus one $n$-point function (op. cite.):
\begin{eqnarray}
&&\hspace{1.5cm}Z_{V}^{(1)}(v_{1},z_{1};\ldots v_{n},z_{n};\tau )=  \notag \\
&&Tr_{V}(Y(q_{z_{1}}^{L(0)}v_{1},q_{z_{1}})\ldots
Y(q_{z_{n}}^{L(0)}v_{n},q_{z_{n}})q^{L(0)}).  \label{Z1npt}
\end{eqnarray}%
Note again that $Z_{V}^{(0)}(v_{1},z_{1};\ldots v_{n},z_{n})$ is recovered
in the degeneration limit $q\rightarrow 0$. For homogeneous primary states $%
v_{i}$ of weight $wt(v_{i})$, $Z_{V}^{(1)}$ parameterizes a global
meromorphic differential form on the torus
\begin{equation*}
\mathcal{F}_{V}^{(1)}(v_{1},\ldots v_{n};\tau
)=Z_{V}^{(1)}(v_{1},z_{1};\ldots v_{n},z_{n};\tau )\prod_{1\leq i\leq
n}(dz_{i})^{wt(v_{i})}.
\end{equation*}

\bigskip Zhu introduced (\cite{Z1}) a second VOA $(V,Y[,],\mathbf{1},\tilde{%
\omega})$ which is isomorphic to $(V,Y(,),\mathbf{1},\omega ).$ It has
vertex operators
\begin{equation}
Y[v,z]=\sum_{n\in \mathbb{Z}}v[n]z^{-n-1}=Y(q_{z}^{L(0)}v,q_{z}-1),
\label{Ysquare}
\end{equation}%
and conformal vector $\tilde{\omega}=$ $\omega -\frac{c}{24}\mathbf{1}$. Let
\begin{equation}
Y[\tilde{\omega},z]=\sum_{n\in \mathbb{Z}}L[n]z^{-n-2} ,  \label{Ywtilde}
\end{equation}%
and write $wt[v]=k$ if $L[0]v=kv$, $V_{[k]}=\{v\in V|wt[v]=k\}$. Only
primary vectors are homogeneous with respect to both $L(0)$ and $L[0]$, in
which case $wt(v)=wt[v]$. Similarly, we define the square bracket LiZ metric
$\langle \ ,\rangle _{\mathrm{sq}}$ which is invariant with respect to the
square bracket adjoint.

\medskip An explicit description of the $n$-point functions for Heisenberg
and lattice VOAs is given in \cite{MT1}, using a Fock basis of $L[0]$%
-homogeneous states. We will make extensive use of the $1$-point and $2$%
-point functions for these VOAs below. We denote any $1$-point function by
\begin{equation}
Z_{V}^{(1)}(v,\tau )=Z_{V}^{(1)}(v,z;\tau )=Tr_{V}(o(v)q^{L(0)}).
\label{Z1_1_pt}
\end{equation}%
$Z^{(1)}(v,\tau )$ is independent of $z$ since only the zero mode $o(v)$
contributes to the trace. Any $n$-point function can be expressed in terms
of $1$-point functions (\cite{MT1}, Lemma 3.1) as follows:
\begin{eqnarray}
&&Z_{V}^{(1)}(v_{1},z_{1};\ldots v_{n},z_{n};\tau )  \notag \\
&=&Z_{V}^{(1)}(Y[v_{1},z_{1}]\ldots Y[v_{n-1},z_{n-1}]Y[v_{n},z_{n}]\mathbf{1%
},\tau )  \label{Z1Ysq1} \\
&=&Z_{V}^{(1)}(Y[v_{1},z_{1n}]\ldots Y[v_{n-1},z_{n-1n}]v_{n},\tau ),
\label{Z1Ysq2}
\end{eqnarray}%
where $z_{in}=z_{i}-z_{n}$. Thus $Z_{V}^{(1)}(v_{1},z_{1};v_{2},z_{2};\tau )$
is a function of $z_{12}$ only, which we denote by
\begin{eqnarray}
Z_{V}^{(1)}(v_{1},v_{2},z_{12},\tau )
&=&Z_{V}^{(1)}(v_{1},z_{1};v_{2},z_{2};\tau )  \notag \\
&=&Tr_{V}(o(Y[v_{1},z_{12}]v_{2})q^{L(0)}).  \label{Z1_2pt}
\end{eqnarray}

We may consider a trivial sewing of a torus with local coordinate $z_{1}$ to
a Riemann sphere with local coordinate $z_{2}$ by identifying the annuli $%
r_{a}\geq |z_{a}|\geq |\epsilon |r_{\bar{a}}^{-1}$ via the sewing relation $%
z_{1}z_{2}=\epsilon .$ Consider $Z_{V}^{(1)}(v_{1},x_{1};\ldots v_{n},x_{n})$
for quasi-primary $v_{i}$ of $L[0]$ weight $wt[v_{i}]$, with $r_{1}\geq
|x_{i}|\geq |\epsilon |r_{2}^{-1}$, and let $y_{i}=\epsilon /x_{i}$. Using (%
\ref{Z1Ysq1}), and employing the square bracket version of (\ref%
{Yxkplus_u_expansion}) with square bracket LiZ metric $\langle \ ,\rangle _{%
\mathrm{sq}}$, we have
\begin{eqnarray*}
&&\hspace{3.5cm}Z_{V}^{(1)}(v_{1},x_{1};\ldots v_{n},x_{n};\tau )= \\
&&\sum_{r\geq 0}\sum_{u\in V_{[r]}}Z_{V}^{(1)}(Y[v_{1},x_{1}]\ldots
Y[v_{k},x_{k}]u;\tau )\langle \bar{u},Y[v_{k+1},x_{k+1}]\ldots Y[v_{n},x_{n}]%
\mathbf{1}\rangle _{\mathrm{sq}},
\end{eqnarray*}%
where the inner sum is taken over any basis $\{u\}$ of $V_{[r]}$, and $\{%
\bar{u}\}$ is the dual basis with respect to $\langle \ ,\rangle _{\mathrm{sq%
}}$. Now
\begin{equation*}
Z_{V}^{(1)}(Y[v_{1},x_{1}]\ldots Y[v_{k},x_{k}]u;\tau )=\mathrm{Res}%
_{z_{1}=0}z_{1}^{-1}Z_{V}^{(1)}(v_{1},x_{1};\ldots v_{k},x_{k};u,z_{1};\tau
).
\end{equation*}%
Using the isomorphism between the round and square bracket formalisms, we
find as before that
\begin{eqnarray*}
&&\hspace{2cm}\langle \bar{u},Y[v_{k+1},x_{k+1}]\ldots Y[v_{n},x_{n}]\mathbf{%
1}\rangle _{\mathrm{sq}} \\
&=&\epsilon ^{r}\mathrm{Res}_{z_{2}=0}z_{2}^{-1}Z_{V}^{(0)}(v_{n},y_{n};%
\ldots v_{k+1},y_{k+1};\bar{u},z_{2})\prod_{k+1\leq i\leq n}(-\frac{\epsilon
}{x_{i}^{2}})^{wt[v_{i}]}.
\end{eqnarray*}%
We thus obtain a natural analogue of Proposition \ref{Prop_Z0sew}:

\begin{proposition}
\label{Prop_Z1sewZ0} For quasiprimary states $v_{i}$ with the above sewing
scheme, we have
\begin{gather*}
\mathcal{F}_{V}^{(1)}(v_{1},\ldots ,v_{n};\tau )= \\
\sum_{r\geq 0}\epsilon ^{r}\sum_{u\in V_{[r]}}\mathrm{Res}%
_{z_{1}=0}z_{1}^{-1}Z_{V}^{(1)}(v_{1},x_{1};\ldots v_{k},x_{k};u,z_{1};\tau
). \\
\mathrm{Res}_{z_{2}=0}z_{2}^{-1}Z_{V}^{(0)}(v_{n},y_{n};\ldots
v_{k+1},y_{k+1};\bar{u},z_{2}).\prod_{1\leq i\leq
k}(dx_{i})^{wt[v_{i}]}\prod_{k+1\leq i\leq n}(dy_{i})^{wt[v_{i}]},\quad
\end{gather*}%
where the inner sum is taken over any basis $\{u\}$ for $V_{[r]}$ and $\{%
\bar{u}\}$ is the dual basis with respect to $\langle \ ,\rangle _{\mathrm{sq%
}}$. \ \ \ \ \ \ \ \ $\square $
\end{proposition}

\subsubsection{An alternative self-sewing of the Riemann sphere and the
Catalan series}

We now consider an alternative construction of a torus by self-sewing a
Riemann sphere with punctures located at the origin and an arbitrary point $%
w $. We will show that the resulting partition function is again (\ref{Z_tau}%
). Choose local coordinates $z_{1}$ in the neighborhood of the origin and $%
z_{2}=z-w$ for $z$ in the neighborhood of $w$. For a complex sewing
parameter $\rho $, identify the annuli $|{\rho }|r_{\bar{a}}^{-1}\leq
\left\vert z_{a}\right\vert \leq r_{a}$ for $a=1,2$ and $|{\rho }|\leq
r_{1}r_{2}$ via the sewing relation
\begin{equation}
z_{1}z_{2}=\rho .  \label{spheresew2}
\end{equation}%
With $\chi $ as in (\ref{chidef}), the annuli do not intersect provided $%
|\chi |<\frac{1}{4}$ (\cite{MT2}) and the torus modular parameter is
(Proposition 9, op.cite.)
\begin{equation}
q=f(\chi ),  \label{qCat}
\end{equation}%
where $f(\chi )$ satisfies $f=\chi (1+f)^{2}$. Thus $f(\chi )$ is the
Catalan series (\ref{Catalan}) familiar from combinatorics (cf. \cite{St}).
We note the following identity (which can be proved by induction)

\begin{lemma}
\label{Lemma_catm}$f(\chi )$ satisfies
\begin{equation*}
f(\chi )^{m}=\sum_{n\geq m}\frac{m}{n}\binom{2n}{n+m}\chi ^{n}.\quad \square
\end{equation*}
\end{lemma}

\bigskip We now define the genus one partition function in the sewing scheme
(\ref{spheresew2}) by
\begin{eqnarray}
&&\hspace{4.5cm}Z_{V,\rho }^{(1)}(\rho ,w)=  \notag \\
&&\sum_{n\geq 0}\rho ^{n}\sum_{u\in V_{n}}\mathrm{Res}_{z_{2}=0}z_{2}^{-1}%
\mathrm{Res}_{z_{1}=0}z_{1}^{-1}\langle \mathbf{1},Y(u,w+z_{2})Y(\bar{u}%
,z_{1})\mathbf{1}\rangle ,  \label{Z1def_rho}
\end{eqnarray}%
where the $\rho \ $subscript refers to the sewing scheme we are currently
using. In fact, this partition function is equivalent to $Z_{V}^{(1)}(q)$:

\begin{theorem}
\label{Theorem_Z1_q_rho} In the sewing scheme (\ref{spheresew2}), we have
\begin{equation}
Z_{V,\rho }^{(1)}(\rho ,w)=Z_{V}^{(1)}(q).  \label{Z1_q_chi}
\end{equation}%
where $q=f(\chi )$ is given by (\ref{qCat}).
\end{theorem}

\noindent \textbf{Proof. }The summand in (\ref{Z1def_rho}) is
\begin{eqnarray*}
\langle \mathbf{1,}Y(u,w)\bar{u}\rangle &=&\langle Y^{\dagger }(u,w)\mathbf{1%
},\bar{u}\rangle \\
&=&(-w^{-2})^{n}\langle Y(e^{wL(1)}u,w^{-1})\mathbf{1},\bar{u}\rangle \\
&=&(-w^{-2})^{n}\langle e^{w^{-1}L(-1)}e^{wL(1)}u,\bar{u}\rangle ,
\end{eqnarray*}%
where we have used (\ref{eq: adj op}) and also $Y(v,z)\mathbf{1}=\exp
(zL(-1))v.$ (See \cite{Ka}\ or \cite{MN} for the latter equality.) Hence we
find that
\begin{eqnarray*}
Z_{V,\rho }^{(1)}(\rho ,w) &=&\sum_{n\geq 0}(-\frac{\rho }{w^{2}}%
)^{n}\sum_{u\in V_{n}}\langle e^{w^{-1}L(-1)}e^{wL(1)}u,\bar{u}\rangle \\
&=&\sum_{n\geq 0}\chi ^{n}Tr_{V_{n}}(e^{w^{-1}L(-1)}e^{wL(1)}).
\end{eqnarray*}%
Expanding the exponentials yields
\begin{equation}
Z_{V,\rho }^{(1)}(\rho ,w)=Tr_{V}(\sum_{r\geq 0}\frac{L(-1)^{r}L(1)^{r}}{%
(r!)^{2}}\chi ^{L(0)}),  \label{Z1rhotrace}
\end{equation}%
an expression which depends only on $\chi $.

\medskip In order to compute (\ref{Z1rhotrace}) we consider the
quasi-primary decomposition of $V$. Let $Q_{m}=\{v\in V_{m}|L(1)v=0\}$
denote the space of quasiprimary states of weight $m\geq 1$. Then $\dim
Q_{m}=p_{m}-p_{m-1}$ with $p_{m}=\dim V_{m}$. Consider the decomposition of $%
V$ into $L(-1)$-descendents of quasi-primaries
\begin{equation}
V_{n}=\bigoplus_{m=1}^{n}L(-1)^{n-m}Q_{m}.  \label{V_qprim}
\end{equation}

\begin{lemma}
\label{Lemma_Quasip}Let $v\in Q_{m}$ for $m\geq 1$. For an integer $n\geq m,$
\begin{equation*}
\sum_{r\geq 0}\frac{L(-1)^{r}L(1)^{r}}{(r!)^{2}}L(-1)^{n-m}v=\binom{2n-1}{n-m%
}L(-1)^{n-m}v.
\end{equation*}
\end{lemma}

\noindent \textbf{Proof.} First use induction on $t\geq 0$ to show that
\begin{eqnarray*}
L(1)L(-1)^{t}v=t(2m+t-1)L(-1)^{t-1}v.
\end{eqnarray*}
Then by induction in $r$ it follows that
\begin{equation*}
\frac{L(-1)^{r}L(1)^{r}}{(r!)^{2}}L(-1)^{n-m}v=\binom{n-m}{r}\binom{n+m-1}{r}%
L(-1)^{n-m}v.
\end{equation*}%
Hence
\begin{eqnarray*}
\sum_{r\geq 0}\frac{L(-1)^{r}L(1)^{r}}{(r!)^{2}}L(-1)^{n-m}v &=&\sum_{r\geq
0}^{n-m}\binom{n-m}{r}\binom{n+m-1}{r}L(-1)^{n-m}v, \\
&=&\binom{2n-1}{n-m}L(-1)^{n-m}v,
\end{eqnarray*}%
where the last combinatorial identity follows from a comparison of the $%
x^{n-m}\,$coefficient of both sides of $%
(1+x)^{n-m}(1+x)^{n+m-1}=(1+x)^{2n-1} $. \ \ \ \ \ \ $\square $

\medskip Lemma \ref{Lemma_Quasip} \ and (\ref{V_qprim}) imply that for $%
n\geq 1,$
\begin{eqnarray*}
Tr_{V_{n}}(\sum_{r\geq 0}\frac{L(-1)^{r}L(1)^{r}}{(r!)^{2}})
&=&\sum_{m=1}^{n}Tr_{Q_{m}}(\sum_{r\geq 0}\frac{L(-1)^{r}L(1)^{r}}{(r!)^{2}}%
L(-1)^{n-m}) \\
&=&\sum_{m=1}^{n}(p_{m}-p_{m-1})\binom{2n-1}{n-m}.
\end{eqnarray*}%
The coefficient of $p_{m}$ is
\begin{equation*}
\binom{2n-1}{n-m}-\binom{2n-1}{n-m-1}=\frac{m}{n}\binom{2n}{m+n},
\end{equation*}%
and hence
\begin{equation*}
Tr_{V_{n}}(\sum_{r\geq 0}\frac{L(-1)^{r}L(1)^{r}}{(r!)^{2}})=\sum_{m=1}^{n}%
\frac{m}{n}\binom{2n}{m+n}p_{m}.
\end{equation*}%
Using Lemma \ref{Lemma_catm}, we find that
\begin{eqnarray*}
Z_{V,\rho }^{(1)}(\rho ,w) &=&1+\sum_{n\geq 1}\chi ^{n}\sum_{m=1}^{n}\frac{m%
}{n}\binom{2n}{m+n}p_{m}, \\
&=&1+\sum_{m\geq 1}p_{m}\sum_{n\geq m}\frac{m}{n}\binom{2n}{m+n}\chi ^{n} \\
&=&1+\sum_{m\geq 1}p_{m}(f(\chi ))^{m},
\end{eqnarray*}%
and Theorem \ref{Theorem_Z1_q_rho} follows. $\ \ \ \ \ \ \square $

\subsection{Genus two}

We now discuss two separate definitions of the genus two $n$-point function
associated respectively with the $\epsilon $- and $\rho $-sewing schemes
reviewed in Sections 2.3 and 2.4.

\medskip In the $\epsilon $-sewing scheme we consider a pair of tori $%
\mathcal{S}_{1},\mathcal{S}_{2}$ with modular parameters\ $\tau _{1},\tau
_{2}$ respectively. We sew them together via the sewing relation (\ref{pinch}%
) of Section 2.3. By definition, in the limit $\epsilon \rightarrow 0$ the
genus two surface degenerates to two tori. For $x_{1},\ldots ,x_{k}\in
\mathcal{S}_{1}$ with $\left\vert x_{i}\right\vert \geq |\epsilon |/r_{2}$
and $y_{k+1},\ldots ,y_{n}\in \mathcal{S}_{2}$ with $\left\vert
y_{i}\right\vert \geq |\epsilon |/r_{1}$ we define the genus two $n$-point
function in the $\epsilon $-formalism by
\begin{gather}
Z_{V,\epsilon }^{(2)}(v_{1},x_{1};\ldots v_{k},x_{k};v_{k+1},y_{k+1};\ldots
v_{n},y_{n};\tau _{1},\tau _{2},\epsilon )=  \notag \\
\sum_{r\geq 0}\epsilon ^{r}\sum_{u\in V_{[r]}}\mathrm{Res}%
_{z_{1}=0}z_{1}^{-1}Z_{V}^{(1)}(v_{1},x_{1};\ldots v_{k},x_{k};u,z_{1};\tau
_{1}).  \notag \\
\mathrm{Res}_{z_{2}=0}z_{2}^{-1}Z_{V}^{(1)}(v_{k+1},y_{k+1};\ldots
v_{n},y_{n};\bar{u},z_{2};\tau _{1}),.  \label{Z2n_pt_eps}
\end{gather}%
As before, the inner sum is taken over any basis $\{u\}$ for $V_{[r]}$ and $%
\{\bar{u}\}$ is the dual basis with respect to $\langle \ ,\rangle _{\mathrm{%
sq}}$. This definition is motivated by Proposition \ref{Prop_Z1sewZ0}. In
this paper we will concentrate on the genus two partition function, i.e. the
$0$-point function. (A discussion of genus two $n$-point functions will
appear elsewhere \cite{MT4}). In the notation of (\ref{Z1_1_pt}), this is
given by
\begin{equation}
Z_{V,\epsilon }^{(2)}(\tau _{1},\tau _{2},\epsilon )=\sum_{n\geq 0}\epsilon
^{n}\sum_{u\in V_{[n]}}Z_{V}^{(1)}(u,\tau _{1})Z_{V}^{(1)}(\bar{u},\tau
_{2}).  \label{Z2_def_eps}
\end{equation}%
Note again that $Z_{V,\epsilon }^{(2)}(\tau _{1},\tau _{2},\epsilon
)\rightarrow Z_{V}^{(1)}(\tau _{1})Z_{V}^{(1)}(\tau _{2})$ in the two tori
degeneration limit $\epsilon \rightarrow 0$.

\medskip In the $\rho $-sewing scheme we self-sew a torus $\mathcal{S}$ with
modular parameter\ $\tau $ via the sewing relation (\ref{rhosew}). For $%
x_{1},\ldots ,x_{n}\in \mathcal{S}$ with $\left\vert x_{i}\right\vert \geq
|\epsilon |/r_{2}$ and $\left\vert x_{i}-w\right\vert \geq |\epsilon |/r_{1}$%
, we define the genus two $n$-point function in the $\rho $-formalism by
\begin{gather}
Z_{V,\rho }^{(2)}(v_{1},x_{1};\ldots v_{n},x_{n};\tau ,w,\rho )=  \notag \\
\sum_{r\geq 0}\rho ^{r}\sum_{u\in V_{[r]}}\mathrm{Res}_{z_{1}=0}z_{1}^{-1}%
\mathrm{Res}_{z_{2}=0}z_{2}^{-1}Z_{V}^{(1)}(\bar{u},w+z_{2};v_{1},x_{1};%
\ldots v_{n},x_{n};u,z_{1};\tau ).  \label{Z2n_pt_rho}
\end{gather}%
This definition is motivated by Proposition \ref{Prop_Z1sewZ0} and (\ref%
{Z1def_rho}). With the notation (\ref{Z1_2pt}), the genus two partition
function is then
\begin{equation}
Z_{V,\rho }^{(2)}(\tau ,w,\rho )=\sum_{n\geq 0}\rho ^{n}\sum_{u\in
V_{[n]}}Z_{V}^{(1)}(\bar{u},u,w,\tau ).  \label{Z2_def_rho}
\end{equation}

\medskip We next consider $Z_{V,\rho }^{(2)}(\tau ,w,\rho )$ in the two-tori
degeneration limit $w,\rho \rightarrow 0$ for fixed $\left\vert \chi
\right\vert <\frac{1}{4}$ of (\ref{chidef}) as reviewed in Section 2.3. We
then find that $Z_{V,\epsilon }^{(2)}(\tau _{1},\tau _{2},\epsilon )$ and $%
Z_{V,\rho }^{(2)}(\tau ,w,\rho )$ agree in this limit:

\begin{theorem}
\label{Theorem_Degen_pt}For fixed $|\chi |<\frac{1}{4}$ we have
\begin{equation*}
\lim_{w,\rho \rightarrow 0}Z_{V,\rho }^{(2)}(\tau ,w,\rho
)=Z_{V}^{(1)}(q)Z_{V}^{(1)}(f(\chi )),
\end{equation*}%
where $f(\chi )$ is the Catalan series (\ref{Catalan}).
\end{theorem}

\noindent \textbf{Proof.} Using (\ref{Z1_2pt}), the inner sum in (\ref%
{Z2_def_rho}) is
\begin{equation*}
Z_{V}^{(1)}(\bar{u},u,w,\tau )=Tr_{V}(o(Y[\bar{u},w]u)q^{L(0)}).
\end{equation*}%
Using the square bracket version of (\ref{Yxkplus_u_expansion}), we have
\begin{equation*}
Y[\bar{u},w]u=\sum_{m\geq 0}\sum_{v\in V_{[m]}}\langle \bar{v},Y[\bar{u}%
,w]u\rangle _{\mathrm{sq}}v.
\end{equation*}%
Similarly to the first part of Theorem \ref{Theorem_Z1_q_rho}, we also have
\begin{eqnarray*}
\langle \bar{v},Y[\bar{u},w]u\rangle _{\mathrm{sq}} &=&(-w^{-2})^{n}\langle
Y[e^{wL[1]}\bar{u},w^{-1}]\bar{v},u\rangle _{\mathrm{sq}} \\
&=&(-w^{-2})^{n}\langle e^{w^{-1}L[-1]}Y[\bar{v},-w^{-1}]e^{wL[1]}\bar{u}%
,u\rangle _{\mathrm{sq}} \\
&=&(-w^{-2})^{n}\langle E[\bar{v},w]\bar{u},u\rangle _{\mathrm{sq}},
\end{eqnarray*}%
where
\begin{equation*}
E[\bar{v},w]=\exp (w^{-1}L[-1])Y[\bar{v},-w^{-1}]\exp (wL[1]).
\end{equation*}%
Hence we find
\begin{eqnarray*}
Z_{V,\rho }^{(2)}(\tau ,w,\rho ) &=&\sum_{m\geq 0}\sum_{v\in
V_{[m]}}\sum_{n\geq 0}\chi ^{n}\sum_{u\in V[n]}\langle E[\bar{v},w]\bar{u}%
,u\rangle Z_{V}^{(1)}(v,q) \\
&=&\sum_{m\geq 0}\sum_{v\in V_{[m]}}Tr_{V}(E[\bar{v},w]\chi
^{L[0]})Z_{V}^{(1)}(v,q).
\end{eqnarray*}%
Now consider
\begin{gather*}
Tr_{V}(E[\bar{v},w]\chi ^{L[0]})= \\
w^{m}\sum\limits_{r,s\geq 0}(-1)^{r+m}\frac{1}{r!s!}Tr_{V}(L[-1]^{r}\bar{v}%
[r-s-m-1]L[1]^{s}\chi ^{L[0]}).
\end{gather*}%
The leading term in $w$ thus arises from $\bar{v}=\mathbf{1}$ of weight $m=0$%
, and is equal to
\begin{equation*}
tr_{V}(E[\mathbf{1},w]\chi ^{L[0]})=Z_{V}^{(1)}(f(\chi )).
\end{equation*}%
This follows from (\ref{Z1rhotrace}) and the isomorphism between the
original and square bracket formalisms. Taking $w\rightarrow 0$ for fixed $%
\chi $, we find $Z_{V,\rho }^{(2)}(\tau ,w,\rho )\rightarrow
Z_{V}^{(1)}(q)Z_{V}^{(1)}(f(\chi )).$ \ \ \ \ \ \ $\square $

\section{Genus two partition function for free bosons in the $\protect%
\epsilon $-formalism}

\subsection{\protect\bigskip The genus two partition function $Z_{M}^{(2)}(%
\protect\tau _{1},\protect\tau _{2},\protect\epsilon )$}

We begin by recalling the general definition (\ref{Z2_def_eps}) of the genus
two partition function associated to a VOA $V$ in the $\epsilon $ formalism:

\begin{equation}
Z_{V,\epsilon }^{(2)}(\tau _{1},\tau _{2},\epsilon )=Z_{V}^{(2)}(\tau
_{1},\tau _{2},\epsilon )=\sum_{n\geq 0}\epsilon ^{n}\sum_{u\in
V_{[n]}}Z_{V}^{(1)}(u,\tau _{1})Z_{V}^{(1)}(\bar{u},\tau _{2}).
\label{eq: part func}
\end{equation}%
The main ingredients in (\ref{eq: part func}) are as follows: $\{u\}$ is a
basis for $V_{[n]},n\geq 0$; $\{\bar{u}\}$ is the dual basis with respect to
the square-bracket Li-Z metric $\langle \ ,\rangle _{\mathrm{sq}}$, and $%
Z_{V}^{(1)}(u,\tau )$ is the genus one graded trace of a state $u$ with
respect to $(V,Y(\ ,))$. In this section we obtain closed formulas for $%
Z_{V}^{(2)}(\tau _{1},\tau _{2},\epsilon )$ in terms of an infinite
determinant and also as an infinite product. We also discuss its
modular properties.

\medskip The partition function (\ref{eq: part func}) is independent of the
choice of basis, in particular if we can choose a diagonal basis $\{u\}$
then $\bar{u}=u/\langle u,u\rangle _{\mathrm{sq}}$ and we obtain
\begin{equation}
Z_{V}^{(2)}(\tau _{1},\tau _{2},\epsilon )=\sum_{n\geq 0}\epsilon
^{n}\sum_{u\in V_{[n]}}\frac{Z_{V}^{(1)}(u,\tau _{1})Z_{V}^{(1)}(u,\tau _{2})%
}{\langle u,u\rangle _{\mathrm{sq}}}.  \label{eq: part func diag}
\end{equation}

Let us consider the case of the rank one Heisenberg (free boson) VOA $M$.
Recalling the definition (\ref{logdet}) we wish to establish the following
closed formula:

\begin{theorem}
\label{Theorem_Z2_boson} Let $M$ be the vertex operator algebra of one free
boson. Then
\begin{equation}
Z_{M}^{(2)}(\tau _{1},\tau _{2},\epsilon )=Z_{M}^{(1)}(\tau
_{1})Z_{M}^{(1)}(\tau _{2})(\det (1-A_{1}A_{2}))^{-1/2},  \label{Z2_1bos}
\end{equation}%
where $A_{1}$ and $A_{2}$ are as in (\ref{Akldef}) and $Z_{M}^{(1)}(\tau
)=q^{1/24}/\eta (\tau )$.
\end{theorem}

\begin{remark}
\label{Rem_Z2mult}\textsc{\ }From Remark \ref{Rem_tensorVOA} it follows that
(\ref{eq: part func}) is multiplicative over tensor products of vertex
operator algebras. Thus the genus two partition function for $l$ free bosons
$M^{l}$ is just the $l^{th}$ power of (\ref{Z2_1bos}).
\end{remark}

\noindent \textbf{Proof of Theorem.} In the following we use the notation
and results of \cite{MT1} concerning the $1$-point functions $%
Z_{M}^{(1)}(u,\tau )$, noting the absence of an overall $q^{-c/24}$ factor
in (\ref{Z_tau}) and (\ref{Z1_1_pt}). Thus $Z_{M}^{(1)}(\tau )=q^{1/24}/\eta
(\tau )$. We take as our diagonal basis of $(V,Y[,])$ the standard Fock
vectors (in the square bracket formulation)
\begin{equation}
v=a[-1]^{e_{1}}...a[-p]^{e_{p}}\mathbf{1.}  \label{eq: sq fock vec}
\end{equation}%
Of course, these Fock vectors correspond in a natural 1-1 manner with
unrestricted partitions, the state $v$ (\ref{eq: sq fock vec}) corresponding
to a partition $\lambda =\{1^{e_{1}}...p^{e_{p}}\}$ of $n=\sum_{1\leq i\leq
p}ie_{i}$. We sometimes write $v=v(\lambda )$ to indicate this
correspondence. As discussed at length in \cite{MT1}, the partition $\lambda
$ may be thought of as a labelled set $\Phi =\Phi _{\lambda }$ with $e_{i}$
elements labelled $i$. One of the main results of \cite{MT1} (loc.cit.
Corollary 1 and eqn.(53)) is then that
\begin{equation}
Z_{M}^{(1)}(v(\lambda ),\tau )=Z_{M}^{(1)}(\tau )\sum_{\phi \in F(\Phi
_{\lambda })}\Gamma (\phi ),  \label{eq: genus1 ptn func}
\end{equation}%
with
\begin{equation}
\Gamma (\phi ,\tau )=\Gamma (\phi )=\prod_{(r,s)}C(r,s,\tau ),
\label{eq: Gamma}
\end{equation}%
for $C$ of (\ref{Ckldef}), where $\phi $ ranges over the elements of $F(\Phi
_{\lambda })$ (the fixed-point-free involutions in $\Sigma (\Phi _{\lambda
}) $) and $(r,s)$ ranges over the orbits of $\phi $ on $\Phi _{\lambda }$.

\medskip With this notation, and using the Fock basis (\ref{eq: sq fock vec}%
) as well as (\ref{eq: inner prod}), (\ref{eq: part func diag}) reads

\begin{equation}
Z_{M}^{(2)}(\tau _{1},\tau _{2},\epsilon )=Z_{M}^{(1)}(\tau
_{1})Z_{M}^{(1)}(\tau _{1})\sum_{\lambda =\{i^{e_{i}}\}}\frac{E(\lambda )}{%
\prod_{i}(-i)^{e_{i}}e_{i}!}\epsilon ^{\sum ie_{i}},
\label{eq: part func lambda}
\end{equation}%
where $\lambda $ ranges over all unrestricted partitions and where we have
set
\begin{eqnarray}
E(\lambda ) &=&\sum_{\phi ,\psi \in F(\Phi _{\lambda })}\Gamma _{1}(\phi
)\Gamma _{2}(\psi ),  \label{eq: expr E} \\
\Gamma _{i}(\phi ) &=&\Gamma (\phi ,\tau _{i}).  \label{eq: expr Gamma}
\end{eqnarray}

We now analyze the nature of the expression $E(\lambda )$ more closely. This
will lead us to the connection between $Z^{(2)}(\tau _{1},\tau _{2},\epsilon
)$ and the chequered cycles discussed  in Section 3.1. The idea is to
use the technique employed in the proof of Proposition 3.10 of \cite{MT1}.
If we fix for a moment a partition $\lambda $ then a pair of
fixed-point-free involutions $\phi ,\psi $ correspond (loc.cit.) to a pair
of complete matchings $\mu _{\phi },\mu _{\psi }$ on the labelled set $\Phi
_{\lambda }$ which we may represent pictorially as

\begin{equation*}
\begin{array}{ccccc}
\overset{r_{1}}{\bullet } & \overset{1}{\longrightarrow } & \overset{s_{1}}{%
\bullet } & \overset{2}{\longrightarrow } & \overset{t_{1}}{\bullet } \\
\overset{r_{2}}{\bullet } & \overset{1}{\longrightarrow } & \overset{s_{2}}{%
\bullet } & \overset{2}{\longrightarrow } & \overset{t_{2}}{\bullet } \\
\vdots &  & \vdots & \vdots & \vdots \\
\overset{r_{b}}{\bullet } & \overset{1}{\longrightarrow } & \overset{s_{b}}{%
\bullet } & \overset{2}{\longrightarrow } & \overset{t_{b}}{\bullet }%
\end{array}%
\end{equation*}

\begin{center}
{\small Fig. 5 \ Two complete matchings }
\end{center}

\noindent Here, $\mu _{\phi }$ is the matching with edges labelled $1$, $\mu
_{\psi }$ the matching with edges labelled $2$, and where we have denoted
the (labelled) elements of $\Phi _{\lambda }$ by $%
\{r_{1},s_{1},...,r_{b},s_{b}\}=\{s_{1},t_{1},...,s_{b},t_{b}\}$. From this
data we may create a chequered cycle in a natural way: starting with some
node of $\Phi _{\lambda }$, apply the involutions $\phi ,\psi $ successively
and repeatedly until the initial node is reached, using the complete
matchings to generate a chequered cycle. The resulting chequered cycle
corresponds to an orbit of $\langle \psi \phi \rangle $ considered as a
cyclic subgroup of $\Sigma (\Phi _{\lambda })$. Repeat this process for each
such orbit to obtain a \emph{chequered diagram} $D$ consisting of the union
of the chequered cycles corresponding to all of the orbits of $\langle \psi
\phi \rangle $ on $\Phi _{\lambda }$. To illustrate, for the partition $%
\lambda =\{1^{2}.2.3^{2}.5\}$ with matchings $\mu _{\phi }=(13)(15)(23)$ and
$\mu _{\psi }=(11)(35)(23)$, the corresponding chequered diagram is

\begin{center}
\begin{picture}(300,60)

\put(100,50){\line(1,0){20}}
\put(95,53){\makebox(0,0){$1$}}
\put(100,50){\circle*{4}}
\put(110,54){\makebox(0,0){\scriptsize $1$}}

\put(120,50){\line(0,-1){20}}
\put(125,53){\makebox(0,0){$3$}}
\put(120,50){\circle*{4}}
\put(123,40){\makebox(0,0){\scriptsize $2$}}

\put(100,30){\line(1,0){20}}
\put(125,27){\makebox(0,0){$5$}}
\put(120,30){\circle*{4}}
\put(110,26){\makebox(0,0){\scriptsize $1$}}

\put(100,50){\line(0,-1){20}}
\put(95,27){\makebox(0,0){$1$}}
\put(100,30){\circle*{4}}
\put(97,40){\makebox(0,0){\scriptsize $2$}}

\put(160,40){\qbezier(0,0)(10,10)(20,0)}
\put(160,40){\circle*{4}}
\put(153,40){\makebox(0,0){$2$}}
\put(170,50){\makebox(0,0){\scriptsize $1$}}

\put(180,40){\circle*{4}}
\put(185,40){\makebox(0,0){$3$}}
\put(160,40){\qbezier(0,0)(10,-10)(20,0)}
\put(170,31){\makebox(0,0){\scriptsize $2$}}

\end{picture}

{\small Fig. 6 Chequered Diagram}
\end{center}

As usual, two chequered diagrams are isomorphic if there is a bijection on
the nodes which preserves edges and labels of nodes and edges. If $\lambda
=\{1^{e_{1}}...p^{e_{p}}\}$ then $\Sigma (\Phi _{\lambda })$ acts on the
chequered diagrams which have $\Phi _{\lambda }$ as underlying set of
labelled nodes. The \emph{label subgroup} $\Lambda $, consisting of the
elements of $\Sigma (\Phi _{\lambda })$ which preserves node labels, is
isomorphic to $\Sigma _{e_{1}}\times \ldots \times \Sigma _{e_{p}}$. It
induces all isomorphisms among these chequered diagrams. Of course $|\Lambda
|=\prod_{1\leq i\leq p}e_{i}!$ . We have almost established the first step
in the proof of Theorem \ref{Theorem_Z2_boson}, namely

\begin{proposition}
\label{Prop_Z2boson_cheq}We have
\begin{equation}
Z_{M}^{(2)}(\tau _{1},\tau _{2},\epsilon )=Z_{M}^{(1)}(\tau
_{1})Z_{M}^{(1)}(\tau _{1})\sum_{D}\frac{\gamma (D)}{|\mathrm{Aut}(D)|},
\label{eq: part func D}
\end{equation}%
where $D$ ranges over isomorphism classes of chequered configurations and
\begin{equation}
\gamma (D)=\frac{E(\lambda )}{\prod_{i}(-i)^{e_{i}}}\epsilon ^{\sum ie_{i}}.
\label{gammaD}
\end{equation}
\end{proposition}

Proposition \ref{Prop_Z2boson_cheq} follows from what we have said together
with (\ref{eq: part func lambda}). It is only necessary to point out that
because the label subgroup induces all isomorphisms of chequered diagrams,
when we sum over isomorphism classes of such diagrams in (\ref{eq: part func
lambda}) the term $\prod_{i}e_{i}!$ must be replaced by $|\mathrm{Aut}(D)|$.
$\square $

\bigskip Recalling the weights (\ref{eq: eqomega}),$\ $we define
\begin{equation*}
\omega (D)=\Pi _{E}\omega (E),
\end{equation*}%
where the product is taken over the edges $E$ of $D$ and $\omega(E)$ is as
in (\ref{omegaval}).

\begin{lemma}
\label{Lemma_om_gamma} For all $D$ we have%
\begin{equation}
\omega (D)=\gamma (D).  \label{omgamma}
\end{equation}
\end{lemma}

\noindent \textbf{Proof. }Let $D$ be determined by a partition $\lambda
=\{1^{e_{1}}...p^{e_{p}}\}$ and a pair of involutions $\phi ,\psi \in F(\Phi
_{\lambda })$, and let $(a,b),(r,s)$ range over the orbits of $\phi $ resp. $%
\psi $ on $\Phi _{\lambda }$. Then we find
\begin{eqnarray*}
\frac{E(\lambda )}{\prod_{i}(-i)^{e_{i}}}\epsilon ^{\sum ie_{i}} &=&\frac{%
\prod_{(a,b)}C(a,b,\tau _{1})\prod_{(r,s)}C(r,s,\tau _{2})}{%
\prod_{i}(-i)^{e_{i}}}\epsilon ^{\sum ie_{i}} \\
&=&(-1)^{\sum e_{i}}\prod_{(ab)}\frac{\epsilon ^{(a+b)/2}}{\sqrt{ab}}%
C(a,b,\tau _{1})\prod_{(rs)}\frac{\epsilon ^{(r+s)/2}}{\sqrt{rs}}C(r,s,\tau
_{1}) \\
&=&(-1)^{\sum e_{i}}\prod_{(ab)}A_{1}(a,b)\prod_{(rs)}A_{2}(r,s) \\
&=&(-1)^{\sum e_{i}}\omega (D).
\end{eqnarray*}%
For $\sum e_{i}$ even this is (\ref{omgamma}), whereas if $\sum e_{i}$ is
odd then $\omega (D)$ vanishes since some Eisenstein series $E_{a+b}(\tau )$
vanishes. \ \ \ \ \ \ $\square $

\bigskip We may represent a chequered diagram formally as a product
\begin{equation}
D=\prod_{i}L_{i}^{m_{i}}  \label{DLprod}
\end{equation}%
in case $D$ is the disjoint union of unoriented chequered cycles $L_{i}$
with multiplicity $m_{i}$. Then $\mathrm{Aut}(D)$ is isomorphic to the
direct product of the groups $\mathrm{Aut}(L_{i}^{m_{i}})$, and
\begin{equation*}
|\mathrm{Aut}(D)|=\prod_{i}|\mathrm{Aut}(L_{i}^{m_{i}})|m_{i}!
\end{equation*}%
Noting that the expression $\omega (D)$ is multiplicative over disjoint
unions of diagrams, we calculate

\begin{eqnarray*}
\sum_{D}\frac{\omega (D)}{|\mathrm{Aut}(D)|} &=&\prod_{L}\sum_{k\geq 0}\frac{%
\omega (L^{k})}{|\mathrm{Aut}(L^{k})|} \\
&=&\prod_{L}\sum_{k\geq 0}\frac{\omega (L)^{k}}{|\mathrm{Aut}(L)|^{k}k!} \\
&=&\prod_{L}\exp \left( \frac{\omega (L)}{|\mathrm{Aut}(L)|}\right) \\
&=&\exp \left( \sum_{L}\frac{\omega (L)}{|\mathrm{Aut}(L)|}\right) ,
\end{eqnarray*}%
where $L$ ranges over isomorphism classes of unoriented chequered cycles.
Now $\mathrm{Aut}(L)$ is either a dihedral group of order $2r$ or a cyclic
group of order $r$ for some $r\geq 1$, depending on whether $L$ admits a
reflection symmetry or not. If we now \emph{orient} our cycles, say in a
clockwise direction, then we can replace the previous sum over $L$ by a sum
over the set of (isomorphism classes of) \emph{oriented} chequered cycles $%
\mathcal{O}$ to obtain
\begin{equation}
\sum_{D}\frac{\omega (D)}{|\mathrm{Aut}(D)|}=\exp \left( \frac{1}{2}%
\sum_{M\in \mathcal{O}}\frac{\omega (M)}{|\mathrm{Aut}(M)|}\right) .
\label{omDM}
\end{equation}

Let $\mathcal{O}_{2n}\subset $ $\mathcal{O}$ denoted the set of oriented
chequered cycles with $2n$ nodes. Then we have

\begin{lemma}
\label{Lemma_TrA1A2}%
\begin{equation}
Tr((A_{1}A_{2})^{n})=\sum_{M\in \mathcal{O}_{2n}}\frac{n}{|\mathrm{Aut}(M)|}%
\omega (M).  \label{trace}
\end{equation}
\end{lemma}

\noindent \textbf{Proof. }The contribution $%
A_{1}(i_{1},i_{2})A_{2}(i_{2},i_{3})\ldots A_{2}(i_{2n},i_{1})$ to the
left-hand-side of (\ref{trace}) is equal to the weight $\omega (M)$ for some
$M\in \mathcal{O}_{2n}$ with vertices $i_{1},i_{2},\ldots i_{2n}$. Let $%
\sigma =\left(
\begin{array}{ccccc}
i_{1} & \ldots & i_{k} & \ldots & i_{2n} \\
i_{3} & \ldots & i_{k+2} & \ldots & i_{2}%
\end{array}%
\right) $ denote the order $n$ permutation of the indices which generates
rotations of $M$. Then $\mathrm{Aut}(M)=\langle \sigma ^{m}\rangle $ for
some $m=n/|\mathrm{Aut}(M)|$. Now sum over all $i_{k}$ to compute $%
Tr((A_{1}A_{2})^{n})$, noting that for inequivalent $M$ the weight $\omega
(M)$ occurs with multiplicity $m$. The Lemma follows. \ \ \ \ \ \ $\square $

\medskip
We may now complete the proof of Theorem \ref{Theorem_Z2_boson}. From (\ref%
{omDM}) and (\ref{trace}) we obtain
\begin{eqnarray*}
\sum_{D}\frac{\omega (D)}{|\mathrm{Aut}(D)|} &=&\exp \left( \frac{1}{2}%
tr(\sum_{n}\frac{1}{n}(A_{1}A_{2})^{n})\right) \\
&=&\exp (-\frac{1}{2}tr(\log (1-A_{1}A_{2}))) \\
&=&\det (\exp (-\frac{1}{2}(\log (1-A_{1}A_{2})))) \\
&=&(\det (1-A_{1}A_{2}))^{-1/2}.\quad \square
\end{eqnarray*}

\bigskip We may also obtain a product formula for $Z_{M}^{(2)}(\tau
_{1},\tau _{2},\epsilon )$ as follows. Recalling the notation (\ref%
{mathcalRdef}), for each oriented chequered cycle $M$, $\mathrm{Aut}(M)$ is
a cyclic group of order $r$ for some $r\geq 1$. Furthermore it is evident
that there is a rotationless chequered cycle $N$ with $\omega (M)=\omega
(N)^{r}$. Indeed, $N$ may be obtained by taking a suitable consecutive
sequence of $n/r$ nodes of $M$, where $n$ is the total number of nodes of $M$%
. We thus see that
\begin{eqnarray*}
\sum_{M\in \mathcal{O}}\frac{\omega (M)}{|\mathrm{Aut}(M)|} &=&\sum_{N\in
\mathcal{R}}\sum_{r\geq 1}\frac{\omega (N)^{r}}{r} \\
&=&-\sum_{N\in \mathcal{R}}\log (1-\omega (N)).
\end{eqnarray*}%
Then (\ref{omDM}) implies%
\begin{equation}
\det (1-A_{1}A_{2})=\prod_{N\in \mathcal{R}}(1-\omega (N)),  \label{detprod}
\end{equation}%
\ and thus we obtain

\begin{theorem}
\label{Theorem_Z2_boson_prod} Let $M$ be the vertex operator algebra of one
free boson. Then
\begin{equation}
Z_{M}^{(2)}(\tau _{1},\tau _{2},\epsilon )= \frac{Z_{M}^{(1)}(\tau
_{1})Z_{M}^{(1)}(\tau _{2})}{\prod_{N\in \mathcal{R}}(1-\omega (N))^{1/2}}.
\label{eq: part func free bos}
\end{equation}
\end{theorem}

\subsection{Holomorphic and modular-invariance properties}

In Section 2.2 we reviewed the genus two $\epsilon $-sewing formalism and
introduced the domain $\mathcal{D}^{\epsilon }$ parameterizing the genus two
surface. An immediate consequence of Theorem \ref{Theorem_Z2_boson} and
Theorem \ref{Theorem_A1A2}b) is the following:

\begin{theorem}
\label{Theorem_Z2_boson_eps_hol}$Z_{M}^{(2)}(\tau _{1},\tau _{2},\epsilon )$
is holomorphic on\ the domain $\mathcal{D}^{\epsilon }$. \ \ \ \ \ $\square $
\end{theorem}

We next consider the automorphic properties of the genus two partition
function with respect to the group $G$ reviewed in Section 2.2. For a vertex
operator algebra of central charge $c$ at genus one, we usually consider the
\emph{modular partition function }
\begin{equation*}
Z_{V,\mathrm{mod}}^{(1)}(\tau )=q^{-c/24}Z_{V}^{(1)}(\tau ),
\end{equation*}%
because of its enhanced $SL(2,\mathbb{Z})$-invariance properties. In
particular, for \emph{two} free bosons this is
\begin{equation}
Z_{M^{2},\mathrm{mod}}^{(1)}(\tau )=\frac{1}{\eta (\tau )^{2}}.
\label{etafunceqn}
\end{equation}%
Let\footnote{%
There should be no confusion between the character $\chi $ introduced here
and the variable $\chi $ used in (\ref{chidef})} $\chi $ be the
character of $SL(2,\mathbb{Z})$ defined by its action on $\eta (\tau )^{-2}$%
, i.e.
\begin{equation}
\eta (\gamma \tau )^{-2}=\chi (\gamma )\eta (\tau )^{-2}(c\tau +d)^{-1},
\label{eq: eta}
\end{equation}%
where $\gamma =\left(
\begin{array}{cc}
a & b \\
c & d%
\end{array}%
\right) \in SL(2,\mathbb{Z})$. Recall (e.g. \cite{Se}) that $\chi (\gamma )$
is a twelfth root of unity. For a function $f(\tau )$ on $\mathbb{H}%
_{1},k\in \mathbb{Z}$ and $\gamma \in SL(2,\mathbb{Z})$, we define
\begin{equation}
f(\tau )|_{k}\gamma =f(\gamma \tau )\ (c\tau +d)^{-k},  \label{slashaction}
\end{equation}%
so that
\begin{equation}
Z_{M^{2},\mathrm{mod}}^{(1)}(\tau )|_{-1}\gamma =\chi (\gamma )Z_{M^{2},%
\mathrm{mod}}^{(1)}(\tau ).  \label{Z1modgam}
\end{equation}

Similarly, at genus two we define the \emph{genus two modular partition
function} for two free bosons
\begin{eqnarray}
Z_{M^{2},\mathrm{mod}}^{(2)}(\tau _{1},\tau _{2},\epsilon )
&=&Z_{M^{2}}^{(2)}(\tau _{1},\tau _{2},\epsilon )q_{1}^{-1/12}q_{2}^{-1/12}
\notag \\
&=&\frac{1}{\eta (\tau _{1})^{2}\eta (\tau _{2})^{2}\det (1-A_{1}A_{2})}.
\label{Z2mod_eps}
\end{eqnarray}%
Analogously to (\ref{slashaction}), we define
\begin{equation}
f(\tau _{1},\tau _{2},\epsilon )|_{k}\gamma =f(\gamma (\tau _{1},\tau
_{2},\epsilon ))\det (C\Omega +D)^{-k}.  \label{eq: Gaction}
\end{equation}%
Here, the action of $\gamma $ on the right-hand-side is as in (\ref{GDeps}).
We have abused notation by adopting the following conventions in (\ref{eq:
Gaction}), which we continue to use below:
\begin{equation}
\Omega =F^{\epsilon }(\tau _{1},\tau _{2},\epsilon ),\ \gamma =\left(
\begin{array}{cc}
A & B \\
C & D%
\end{array}%
\right) \in Sp(4,\mathbb{Z})  \label{Omegaconvention}
\end{equation}%
where $F^{\epsilon }$ is as in Theorem \ref{Theorem_period_eps}, and $\gamma
$ is identified with an element of $Sp(4,\mathbb{Z})$ via (\ref{Gamma1Gamma2}%
)-(\ref{betagen}). Note that (\ref{eq: Gaction}) defines a right action of $%
G $ on functions $f(\tau _{1},\tau _{2},\epsilon )$. We will establish the
natural extension of (\ref{Z1modgam}). To describe this, introduce the
character $\chi ^{(2)}$ of $G$ defined by
\begin{equation*}
\chi ^{(2)}(\gamma _{1}\gamma _{2}\beta ^{m})=(-1)^{m}\chi (\gamma
_{1}\gamma _{2}),\quad \ \gamma _{i}\in \Gamma _{a},a=1,2.
\end{equation*}%
($\Gamma _{i}$ is as in (\ref{Gamma1Gamma2})). Thus $\chi ^{(2)}$ takes
values which are twelfth roots of unity, and we have

\begin{theorem}
\label{Theorem_Z2_G}If $\gamma \in G$ then
\begin{equation*}
Z_{M^{2},\mathrm{mod}}^{(2)}(\tau _{1},\tau _{2},\epsilon )|_{-1}\gamma
=\chi ^{(2)}(\gamma )Z_{M^{2},\mathrm{mod}}^{(2)}(\tau _{1},\tau
_{2},\epsilon ).
\end{equation*}
\end{theorem}

\begin{corollary}
\label{Cor_Z2_24G} Set $Z_{M^{24},\mathrm{mod}}^{(2)}=(Z_{M^{2},\mathrm{mod}%
}^{(2)})^{12}$. Then for $\gamma \in G$ we have
\begin{equation*}
Z_{M^{24},\mathrm{mod}}^{(2)}(\tau _{1},\tau _{2},\epsilon )|_{-12}\gamma
=Z_{M^{24},\mathrm{mod}}^{(2)}(\tau _{1},\tau _{2},\epsilon ).
\end{equation*}
\end{corollary}

\noindent \textbf{Proof. } We will give two different proofs of this result.
Using the convention (\ref{Omegaconvention}), we have to show that
\begin{equation}
Z_{M^{2},\mathrm{mod}}^{(2)}(\gamma (\tau _{1},\tau _{2},\epsilon ))\det
(C\Omega +D)=\chi ^{(2)}(\gamma )Z_{M^{2},\mathrm{mod}}^{(2)}(\tau _{1},\tau
_{2},\epsilon )  \label{eq: Z^2identity}
\end{equation}%
for $\gamma \in G$, and it is enough to do this for a generating set of $G$.
If $\gamma =\beta $ then the result is clear since $\det (C\Omega +D)=\chi
^{(2)}(\beta )=-1$ and $\beta $ exchanges $\tau _{1}$ and $\tau _{2}$. So we
may assume that $\gamma =(\gamma _{1},\gamma _{2})\in \Gamma _{1}$ $\times $
$\Gamma _{2}$.

\medskip Our first proof utilizes the determinant formula (\ref{Z2_1bos}) as
follows. For $\gamma _{1}\in \Gamma _{1}$, define $A_{a}^{\prime }(k,l,\tau
_{a},\epsilon )=A_{a}(k,l,\gamma _{1}\tau _{a},\frac{\epsilon }{c_{1}\tau
_{1}+d_{1}})$ following (\ref{GDeps}). We find from Section 4.4 of \cite{MT2}
that
\begin{eqnarray*}
1-A_{1}^{\prime }A_{2}^{\prime } &=&1-A_{1}A_{2}-\kappa \Delta A_{2} \\
&=&(1-\kappa S).(1-A_{1}A_{2}),
\end{eqnarray*}%
where
\begin{eqnarray*}
\Delta (k,l) &=&\delta _{k1}\delta _{l1}, \\
\kappa &=&-\frac{\epsilon }{2\pi i}\frac{c_{1}}{c_{1}\tau _{1}+d_{1}}, \\
S(k,l) &=&\delta _{k1}(A_{2}(1-A_{1}A_{2})^{-1})(1,l).
\end{eqnarray*}%
Since $\det (1-A_{1}A_{2})$ and $\det (1-A_{1}^{\prime }A_{2}^{\prime })$
are convergent on $\mathcal{D}^{\epsilon }$ we find%
\begin{equation*}
\det (1-A_{1}^{\prime }A_{2}^{\prime })=\det (1-\kappa S)\det (1-A_{1}A_{2}).
\end{equation*}%
But
\begin{eqnarray*}
\det (1-\kappa S) &=&\left\vert
\begin{array}{cccc}
1-\kappa S(1,1) & -\kappa S(1,2) & -\kappa S(1,3) & \cdots \\
0 & 1 & 0 & \cdots \\
0 & 0 & 1 & \cdots \\
\vdots & \vdots & \vdots & \ddots%
\end{array}%
\right\vert \\
&=&1-\kappa S(1,1) \\
&=&\frac{c_{1}\Omega _{11}+d_{1}}{c_{1}\tau _{1}+d_{1}},
\end{eqnarray*}%
using (\ref{Om11eps}). Thus
\begin{equation*}
\det (1-A_{1}^{\prime }A_{2}^{\prime })=\frac{c_{1}\Omega _{11}+d_{1}}{%
c_{1}\tau _{1}+d_{1}}.\det (1-A_{1}A_{2}),
\end{equation*}%
which implies (\ref{eq: Z^2identity}) for $\gamma _{1}\in \Gamma _{1}$. A
similar proof applies for $\gamma _{2}\in \Gamma _{2}$. \ \ \ \ \ \

\medskip
The second proof uses Proposition \ref{Prop_Om12_R21expansion} together with
(\ref{eq: omega_12prod}), which tell us that
\begin{equation}
Z_{M^{2},\mathrm{mod}}^{(2)}(\tau _{1},\tau _{2},\epsilon )=\frac{-2\pi
i\Omega _{12}}{\epsilon \eta (\tau _{1})^{2}\eta (\tau _{2})^{2}}\prod_{%
\mathcal{R}^{\prime }}(1-\omega (L))^{-1},  \label{eq: Z^2 omega_12}
\end{equation}%
where $\mathcal{R}^{\prime }=\mathcal{R}\setminus \mathcal{R}_{21}$. Now in
general a term $\omega (L)$ will not be invariant under the action of $%
\gamma $. This is because of the presence of non-modular terms $A_{a}(1,1)$
arising from $E_{2}(\tau _{a})$. But it is clear from (\ref{GDeps}) and the
definition (\ref{Ckldef}) of $C(k,l,\tau )$ together with its
modular-invariance properties that if $L\in \mathcal{R}^{\prime }$ then such
terms are absent and $\omega (L)$ \emph{is} invariant. So the product term
in (\ref{eq: Z^2 omega_12}) is invariant under the action of $\gamma $.

Next, we see from (\ref{GDeps}) that the expression $\epsilon \eta (\tau
_{1})^{2}\eta (\tau _{2})^{2}$ is invariant under the action of $\gamma $ up
to a scalar $\chi (\gamma _{1})\chi (\gamma _{2})=\chi ^{(2)}(\gamma )$.
This reduces the proof of (\ref{eq: Z^2identity}) to showing that
\begin{equation*}
(\gamma _{1},\gamma _{2}):\Omega _{12}\mapsto \Omega _{12}\ \det (C\Omega
+D)^{-1},
\end{equation*}%
and this is implicit in (\ref{eq: modtrans}) upon applying Theorem \ref%
{TheoremGequiv}. This completes the second proof of Theorem \ref%
{Theorem_Z2_G}. $\square $

\begin{remark}
\label{automorphyremark} An unusual feature of the formulas in Theorem \ref%
{Theorem_Z2_G} and Corollary \ref{Cor_Z2_24G} is that the definition of the
automorphy factor $\det (C\Omega +D)$ requires the map $F^{\epsilon
}:D^{\epsilon }\rightarrow \mathbb{H}_{2}$.
\end{remark}

\begin{remark}
\label{Remark_cofact} The reason for the appearance of the factor $-2\pi
i\Omega _{12}/\epsilon $ in (\ref{eq: Z^2 omega_12}) can be understood as
follows. From Theorem \ref{Theorem_period_eps} we see that $-2\pi i\Omega
_{12}/\epsilon =(1-A_{1}A_{2})^{-1}(1,1)$. However
\begin{equation}
(1-A_{1}A_{2})^{-1}(1,1)=\frac{c(1,1)}{\det (1-A_{1}A_{2})},
\label{cofactor}
\end{equation}%
where $c(1,1)$ is the $(1,1)$ cofactor for $1-A_{1}A_{2}$. Comparison with (%
\ref{eq: Z^2 omega_12}) shows that $c(1,1)=\prod_{\mathcal{R}^{\prime
}}(1-\omega (L))^{-1}$ which contains no $A_{a}(1,1)$ terms, and hence is
modular invariant.
\end{remark}

\section{Genus two partition function for lattice theories in the $\protect%
\epsilon$-formalism}

Let $L$ be an even lattice with $V_{L}$ the corresponding lattice theory
vertex operator algebra. The underlying Fock space is
\begin{equation}
V_{L}=M^{l}\otimes C[L]=\oplus _{\alpha \in L}M^{l}\otimes e^{\alpha },
\label{lattice fock space}
\end{equation}%
where $M^{l}$ is the corresponding Heisenberg free boson theory of rank $l=$
dim $L$ based on $H=C\otimes _{Z}L$. We follow Section 4.1 and \cite{MT1}
concerning further notation for lattice theories.

\medskip The general shape of $Z_{V_{L}}^{(2)}(\tau _{1},\tau _{2},\epsilon
) $ is as in (\ref{eq: part func}). Note that the modes of a state $u\otimes
e^{\alpha }$ map $M^{l}\otimes e^{\beta }$ to $M^{l}\otimes e^{\alpha +\beta
}$. Thus if $\alpha \neq 0$ then $Z_{V_{L}}^{(1)}(u\otimes e^{\alpha },\tau
) $ vanishes, and as a result we see that
\begin{eqnarray}
Z_{V_{L}}^{(2)}(\tau _{1},\tau _{2},\epsilon ) &=&\sum_{n\geq 0}\epsilon
^{n}\sum_{u\in M_{[n]}^{l}}\frac{Z_{V_{L}}^{(1)}(u,\tau
_{1})Z_{V_{L}}^{(1)}(u,\tau _{2})}{\langle u,u\rangle _{\mathrm{sq}}}  \notag
\\
&=&\sum_{n\geq 0}\epsilon ^{n}\sum_{u\in M_{[n]}^{l}}\sum_{\alpha ,\beta \in
L}\frac{Z_{M^{l}\otimes e^{\alpha }}^{(1)}(u,\tau _{1})Z_{M^{l}\otimes
e^{\beta }}^{(1)}(u,\tau _{2})}{\langle u,u\rangle _{\mathrm{sq}}}.
\label{eq: lattice part func diag}
\end{eqnarray}%
Here, as in (\ref{eq: part func diag}), $u$ ranges over a diagonal basis for
$M_{[n]}^{l}$; $M^{l}\otimes e^{\alpha }$ should be viewed as a simple
module for $M^{l}$. An explicit formula for $Z_{M^{l}\otimes e^{\alpha
}}^{(1)}(u,\tau )$ was given in \cite{MT2} (Corollary 3 and Theorem 1). We
are going to use these results, much as in the case of free bosons carried
out in Section 6, to elucidate (\ref{eq: lattice part func diag}). Indeed,
we will establish

\begin{theorem}
\label{Theorem_Z2_L_eps}We have
\begin{equation}
Z_{V_{L}}^{(2)}(\tau _{1},\tau _{2},\epsilon )=Z_{M^{l}}^{(2)}(\tau
_{1},\tau _{2},\epsilon )\theta _{L}^{(2)}(\Omega ),
\label{normalizedpartfunc}
\end{equation}%
where $\theta _{L}^{(2)}(\Omega )$ is the (genus two) Siegel theta function
associated to $L$ (e.g. \cite{F})%
\begin{equation}
\theta _{L}^{(2)}(\Omega )=\sum_{\alpha ,\beta \in L}\exp (\pi i((\alpha
,\alpha )\Omega _{11}-2(\alpha ,\beta )\Omega _{12}+(\beta ,\beta )\Omega
_{22})),  \label{eq: g2omega}
\end{equation}
\end{theorem}

\begin{remark}
\label{Rem_Siegel_Theta}Note that $\theta _{L}^{(2)}$ is independent of the
sign of the coefficient of $\Omega _{12}$ in the exponent.
\end{remark}

We will deal first with the case $l=1$, when $M=M^{1}$, and subsequently
show how to modify the argument to handle a general lattice. To this end,
let $a$ be a unit vector of $H$ and use the Fock basis vectors $v=v(\lambda
) $ (cf. (\ref{eq: sq fock vec})) identified with partitions $\lambda
=\{i^{e_{i}}\}$ as in Section 6. Recall that $\lambda $ defines a labelled
set $\Phi _{\lambda }$ with $e_{i}$ nodes labelled $i$. It is useful to
re-state Corollary 3 of \cite{MT1} in the following form:
\begin{equation}
Z_{M\otimes e^{\alpha }}^{(1)}(v,\tau )=Z_{M}^{(1)}(\tau )q^{(\alpha ,\alpha
)/2}\sum_{\phi }\Gamma _{\lambda ,\alpha }(\phi ).
\label{eq: Z^1latticeformula}
\end{equation}%
Here, $\phi $ ranges over the set of involutions
\begin{equation}
\mathrm{Inv}_{1}(\Phi _{\lambda })=\{\phi \in \mathrm{Inv}(\Phi _{\lambda
})|\ p\in \mathrm{Fix}(\phi )\Rightarrow p\mbox{ has label}\ 1\}.
\label{eq: inv}
\end{equation}%
(In words, $\phi $ is an involution in the symmetric group $\Sigma (\Phi
_{\lambda })$ such that all fixed-points of $\phi $ carry the label $1$.
Note that this includes the fixed-point-free involutions, which were the
only involutions which played a role in the case of free bosons. The main
difference between free bosonic and lattice theories is the need to include
additional involutions in the latter case.)  Finally,
\begin{equation}
\Gamma _{\lambda ,\alpha }(\phi ,\tau )=\Gamma _{\lambda ,\alpha }(\phi
)=\prod_{\Xi }\Gamma (\Xi ),  \label{eq: gammavalpha}
\end{equation}%
where $\Xi $ ranges over the orbits (of length $\leq 2$) of $\phi $ acting
on $\Phi _{\lambda }$ and
\begin{equation}
\Gamma (\Xi )=\left\{
\begin{array}{ll}
C(r,s,\tau ), & \mbox{if $\Xi = \{r, s \}$,} \\
(a,\alpha ), & \mbox{if $\Xi = \{1 \}$}.%
\end{array}%
\right.  \label{gamma_xi}
\end{equation}%
 From (\ref{eq: lattice part func diag})-(\ref{eq: gammavalpha}) we get
\begin{gather}
Z_{V_{L}}^{(2)}(\tau _{1},\tau _{2},\epsilon )=  \notag \\
Z_{M}^{(1)}(\tau _{1})Z_{M}^{(1)}(\tau _{2})\sum_{\alpha ,\beta \in
L}\sum_{\lambda =\{i^{e_{i}}\}}\frac{E_{\alpha ,\beta }(\lambda )}{%
\prod_{i}(-i)^{e_{i}}e_{i}!}q_{1}^{(\alpha ,\alpha )/2}q_{2}^{(\beta ,\beta
)/2}\epsilon ^{\sum ie_{i}},  \label{eq: newpartfunc}
\end{gather}%
where
\begin{equation}
E_{\alpha ,\beta }(\lambda )=\sum_{\phi ,\psi \in \mathrm{Inv}_{1}(\Phi
_{\lambda })}\Gamma _{\lambda ,\alpha }(\phi ,\tau _{1})\Gamma _{\lambda
,\beta }(\psi ,\tau _{2}).  \label{eq: Eeqn}
\end{equation}%
(Compare with eqns. (\ref{eq: part func lambda}) - (\ref{eq: expr E}).) Now
we follow the proof of Proposition \ref{Prop_Z2boson_cheq} to obtain an
expression analogous to (\ref{eq: part func D}), namely
\begin{gather}
Z_{V_{L}}^{(2)}(\tau _{1},\tau _{2},\epsilon )=  \notag \\
Z_{M}^{(1)}(\tau _{1})Z_{M}^{(1)}(\tau _{2})\sum_{\alpha ,\beta \in
L}\sum_{D_{\alpha ,\beta }}\frac{\gamma ^{0}(D_{\alpha ,\beta })}{|\mathrm{%
Aut}(D_{\alpha ,\beta })|}q_{1}^{(\alpha ,\alpha )/2}q_{2}^{(\beta ,\beta
)/2},  \label{eq: lattice part func D}
\end{gather}%
the meaning of which we now enlarge upon. Compared to (\ref{eq: part func D}%
), the chequered diagrams $D_{\alpha ,\beta }$ which occur in (\ref{eq:
lattice part func D}) not only depend on lattice elements but are also more
general than before, in that they reflect the fact that the relevant
involutions may now have fixed-points. Thus $D_{\alpha ,\beta }$ is the
union of its (as yet unoriented) connected components which are either
chequered cycles as before or else chequered necklaces as introduced in
Section 3.2. Necklaces arise from orbits of the group $\langle \psi \phi
\rangle $ on $\Phi _{\lambda }$ in which one of the nodes in the orbit is a
fixed-point of $\phi $ or $\psi $. In that case the orbit will generally
contain two such nodes which comprise the end nodes of the necklace. Note
that these end nodes necessarily carry the label $1$ (cf. (\ref{eq: inv})).
There is degeneracy when \emph{both} $\phi $ and $\psi $ fix the node, in
which case the degenerate necklace obtains. To be quite precise, the
necklaces that we are presently dealing with are not quite the same as those
of Section 3.2: not only are they unoriented, but they should be conceived
as having a loop (self-edge) at each end node that carries a label $%
(a,\alpha )$ or $(a,\beta )$ (or the product of the two in case of
degeneracy). For convenience, we call these \emph{chequered* necklaces} and
denote their isomorphism class by $\mathcal{N}_{\alpha ,\beta }^{\ast }$.

\medskip Similarly to (\ref{gammaD}), the term $\gamma ^{0}(D_{\alpha ,\beta
})$ in (\ref{eq: lattice part func D}) is given by
\begin{equation}
\gamma ^{0}(D_{\alpha ,\beta })=\frac{\prod_{\Xi _{1}}\Gamma (\Xi
_{1})\prod_{\Xi _{2}}\Gamma (\Xi _{2})}{\prod_{i}(-i)^{e_{i}}}\epsilon
^{\sum ie_{i}},  \label{eq: lattice gamma(D)}
\end{equation}%
where $\Xi _{1},\Xi _{2}$ range over the orbits of $\phi ,\psi $
respectively on $\Phi _{\lambda }$. As usual the summands in (\ref{eq:
lattice part func D}) are multiplicative over connected components of the
chequered diagram. This applies, in particular, to the chequered cycles
which occur, and these are independent of the lattice elements. As a result,
(\ref{eq: lattice part func D}) factors as a product of two expressions, the
first a sum over diagrams consisting only of chequered cycles and the second
a sum over diagrams consisting only of chequered* necklaces. However, the
first expression corresponds precisely to the genus two partition function
for the free boson (Proposition \ref{Prop_Z2boson_cheq}). We thus obtain
\begin{equation}
\frac{Z_{V_{L}}^{(2)}(\tau _{1},\tau _{2},\epsilon )}{Z_{M}^{(2)}(\tau
_{1},\tau _{2},\epsilon )}=\sum_{\alpha ,\beta \in L}\sum_{D_{\alpha ,\beta
}^{N}}\frac{\gamma ^{0}(D_{\alpha ,\beta }^{N})}{|\mathrm{Aut}(D_{\alpha
,\beta }^{N})|}q_{1}^{(\alpha ,\alpha )/2}q_{2}^{(\beta ,\beta )/2},
\label{eq: spliteqn}
\end{equation}%
where here $D_{\alpha ,\beta }^{N}$ ranges over all chequered diagrams all
of whose connected components are chequered* necklaces in $\mathcal{N}%
_{\alpha ,\beta }^{\ast }$. So (at least if $l=1$) Theorem \ref%
{Theorem_Z2_L_eps} is reduced to establishing

\begin{proposition}
\label{Prop_thetaL_eps}We have
\begin{equation}
\theta _{L}^{(2)}(\Omega )=\sum_{\alpha ,\beta \in L}\sum_{D_{\alpha ,\beta
}^{N}}\frac{\gamma ^{0}(D_{\alpha ,\beta }^{N})}{|\mathrm{Aut}(D_{\alpha
,\beta }^{N})|}q_{1}^{(\alpha ,\alpha )/2}q_{2}^{(\beta ,\beta )/2}.
\label{eq: redomega}
\end{equation}
\end{proposition}

We may apply the argument of (\ref{DLprod}) et. seq. to the inner sum in (%
\ref{eq: redomega}) to write it as an exponential expression
\begin{equation}
\exp \{\pi i((\alpha ,\alpha )\tau _{1}+(\beta ,\beta )\tau
_{2})+\sum_{N^{\ast }\in \mathcal{N}_{\alpha ,\beta }^{\ast }}\frac{\gamma
^{0}(N^{\ast })}{|\mathrm{Aut}(N^{\ast })|}\},  \label{eq: expgammasum}
\end{equation}%
where $N^{\ast }$ ranges over all chequered* necklaces $\mathcal{N}_{\alpha
,\beta }^{\ast }$.

\medskip Recall the isomorphism class $\mathcal{N}_{ab}$ of oriented
chequered necklaces of type $ab$ as displayed in Fig.4 of Section 3.2. We
can similarly consider $N^{\ast }$ a chequered necklaces* of type $ab$. Then
$|\mathrm{Aut}(N^{\ast })|\leq 2$ with equality if, and only if, $a=b$ and
the two orientations of $N^{\ast }$ are isomorphic. Then (\ref{eq:
expgammasum}) can be written as
\begin{equation}
\exp \{\pi i((\alpha ,\alpha )\tau _{1}+(\beta ,\beta )\tau _{2})+\frac{1}{2}%
\sum_{a,b\in \{1,2\}}\sum_{N^{\ast }}\gamma ^{0}(N^{\ast })\}.
\label{eq: orexpgammasum}
\end{equation}%
where here $N^{\ast }$ ranges over \emph{oriented} chequered* necklaces of
type $ab$.

\medskip From (\ref{gamma_xi}) we see that the contribution of the end nodes
to $\gamma ^{0}(N^{\ast })$ is equal to $(a,\alpha )^{2}=(\alpha ,\alpha )$
for a type $11$ necklace, and similarly $(\beta ,\beta ),(\alpha ,\beta )$
and $(\beta ,\alpha )$ for types $22,12$ and $21$ respectively. The
remaining factors of $\gamma ^{0}(N^{\ast })$ have product $\gamma
(N)=\omega (N)$ for some $N\in \mathcal{N}_{ab}$ by Lemma \ref%
{Lemma_om_gamma}. Then (\ref{eq: orexpgammasum}) may be re-expressed as

\begin{equation}
\sum_{\alpha ,\beta \in L}\exp \{\frac{1}{2}(\alpha ,\alpha )(2\pi i\tau
_{1}+\epsilon \omega _{11})+\frac{1}{2}(\beta ,\beta )(2\pi i\tau
_{2}+\epsilon \omega _{22})+(\alpha ,\beta )\epsilon \omega _{12}\},
\label{thetaLexp}
\end{equation}%
recalling $\omega _{ab}=\sum_{N\in \mathcal{N}_{ab}}\omega (N)$ and where $%
\omega _{12}=\omega _{21}$. Note that in passing to $\omega (N)$ the (weight
of) end nodes of each necklace contribute an additional overall factor $%
\epsilon $, as displayed in (\ref{thetaLexp}). (\ref{thetaLexp}) reproduces (%
\ref{eq: g2omega}) on applying Proposition \ref{Propepsperiodgraph}. This
completes the proof of Theorem \ref{Theorem_Z2_L_eps} in the case of a rank $%
1$ lattice.

\bigskip We now consider the general case of an even lattice $L$ of rank $l$%
. We again start with (\ref{eq: lattice part func diag}), using results of
\cite{MT1} to describe $Z_{M^{l}\otimes e^{\alpha }}^{(1)}(u,\tau )$. We
carry over the notation of \cite{MT1} (especially the paragraphs prior to
Theorem 1), in particular $a_{1},...,a_{l}$ is an orthonormal base of $%
H=C\otimes L$ and $u$ is as in eqn.(75) (loc. cit.). Then the case $n=1$ and
$\alpha _{1}=0$ of Theorem 1 (loc. cit.) tells us that
\begin{equation}
Z_{M^{l}\otimes e^{\alpha }}^{(1)}(u,\tau )=Z_{M^{l}}^{(1)}(\tau )q^{(\alpha
,\alpha )/2}\prod_{r=1}^{l}\sum_{\phi \in \mathrm{Inv}_{1}(\Phi _{\lambda
^{r})}}\Gamma ^{r}(\phi ).  \label{genlatticegamma}
\end{equation}%
Here, $\Phi _{\lambda ^{r}}$ is the labelled set corresponding to the $%
r^{th} $ partition $\lambda ^{r}$, which itself is the partition determined
by the occurrences of $a_{r}$ in the representation of the state $u$; $%
\Gamma ^{r}(\phi )$ is as in (\ref{eq: gammavalpha}) with the understanding
that the relevant labelled set is now that determined by $\lambda ^{r}$. As
a result, we obtain an expression for $Z_{V_{L}}^{(2)}(\tau _{1},\tau
_{2},\epsilon )$ analogous to (\ref{eq: newpartfunc}), the only difference
now being that the sum over partitions $\lambda $ must be replaced by a sum
over $l$-tuples of partitions $\lambda ^{r},1\leq r\leq l$. In terms of
chequered diagrams, it is much as before except that the connected
components (either chequered cycles or necklaces) have an additional overall
color $r$ in the range $1\leq r\leq l$. Then the proof of Proposition \ref%
{Prop_Z2boson_cheq} still applies and we obtain the analogues of (\ref{eq:
lattice part func D}) - (\ref{eq: orexpgammasum}) in which we again must sum
also over the $l$ different colors of the pertinent graphs. Explicitly, (\ref%
{eq: orexpgammasum}) now reads
\begin{equation}
\exp \{\pi i((\alpha ,\alpha )\tau _{1}+(\beta ,\beta )\tau _{2})+\frac{1}{2}%
\sum_{r=1}^{l}\sum_{a,b\in \{1,2\}}\sum_{N^{r\ast }}\gamma ^{0}(N^{r\ast
})\}.  \label{eq: genorexpgammasum}
\end{equation}%
Finally, the contribution of the end nodes of $N^{r\ast }$ to $\gamma
^{0}(N^{r\ast })$ is $(a_{r},\alpha )^{2}$ for type 11 necklaces, and
similarly for the other types. This being the only dependence of the
expression $\gamma ^{0}(N^{r\ast })$ on $r$, it follows as before that (\ref%
{eq: genorexpgammasum}) is equal to
\begin{eqnarray*}
&&\exp \{\pi i((\alpha ,\alpha )\tau _{1}+(\beta ,\beta )\tau _{2})+\frac{%
\epsilon }{2}\sum_{r=1}^{l}(a_{r},\alpha )^{2}\omega _{11} \\
&&+\epsilon \sum_{r=1}^{l}(a_{r},\alpha )(a_{r},\beta )\omega _{12}+\frac{%
\epsilon }{2}\sum_{r=1}^{l}(a_{r},\beta )^{2}\omega _{22}\},
\end{eqnarray*}%
that is
\begin{equation*}
\exp \{\frac{1}{2}(\alpha ,\alpha )(2\pi i\tau _{1}+\epsilon \omega _{11})+%
\frac{1}{2}(\beta ,\beta )(2\pi i\tau _{2}+\epsilon \omega _{22})+\epsilon
(\alpha ,\beta )\omega _{12}\}.
\end{equation*}%
The remaining argument is now as before, and completes the proof of Theorem %
\ref{Theorem_Z2_L_eps}. \ \ \ \ \ \ \ $\square $

\bigskip We complete this Section with a brief discussion of the automorphic
properties of $Z_{V_{L}}^{(2)}(\tau _{1},\tau _{2},\epsilon )$. There is
more that one can say here, but a fuller discussion must wait for another
time. The function $\theta _{L}^{(2)}(\Omega )$ is a Siegel modular form of
weight $l/2$ (\cite{F}), in particular it is holomorphic on the Siegel upper
half-space $\mathbb{H}_{2}$. From Theorems \ref{Theorem_period_eps}, \ref%
{Theorem_Z2_boson_eps_hol} and \ref{Theorem_Z2_L_eps}, we deduce

\begin{theorem}
\label{Theorem_Z2_lattice_eps_hol}$Z_{V_{L}}^{(2)}(\tau _{1},\tau
_{2},\epsilon )$ is holomorphic on\ the domain $\mathcal{D}^{\epsilon }$. \
\ \ \ \ $\square $
\end{theorem}

We can obtain the automorphic properties of $Z_{V_{L}}^{(2)}(\tau _{1},\tau
_{2},\epsilon )$ in the same way using that for $\theta _{L}^{(2)}(\Omega )$
together with Theorem \ref{Theorem_Z2_G}. Rather than do this explicitly,
let us introduce a third variation of the $Z$-function, namely the \emph{%
normalized partition function}
\begin{equation}
\hat{Z}_{V_{L},\epsilon }^{(2)}(\tau _{1},\tau _{2},\epsilon )=\frac{%
Z_{V_{L}}^{(2)}(\tau _{1},\tau _{2},\epsilon )}{Z_{M^{l}}^{(2)}(\tau
_{1},\tau _{2},\epsilon )}.  \label{ZV_hat_eps}
\end{equation}%
Bearing in mind the convention (\ref{Omegaconvention}), what (\ref%
{normalizedpartfunc}) says is that there is a commuting diagram of
holomorphic maps
\begin{equation}
\begin{array}{lll}
\ \hspace{0cm}\mathcal{D}^{\epsilon } & \overset{F^{\epsilon }}{%
\longrightarrow } & \hspace{0.7cm}\mathbb{H}_{2} \\
\hat{Z}_{V_{L},\epsilon }^{(2)}\searrow &  & \swarrow \theta _{L}^{(2)} \\
& \hspace{0.2cm}\mathbb{C} &
\end{array}
\label{commhol_diag}
\end{equation}

Furthermore, the $G$-actions on the two functions in question are
compatible. More precisely, if $\gamma \in G$ then we have
\begin{eqnarray}
\hat{Z}_{V_{L},\epsilon }^{(2)}(\tau _{1},\tau _{2},\epsilon )|_{l/2}\
\gamma &=&\hat{Z}_{V_{L},\epsilon }^{(2)}(\gamma (\tau _{1},\tau
_{2},\epsilon ))\det (C\Omega +D)^{-l/2}  \notag \\
&=&\theta _{L}^{(2)}(F^{\epsilon }(\gamma (\tau _{1},\tau _{2},\epsilon
)))\det (C\Omega +D)^{-l/2}\qquad (\text{(\ref{commhol_diag})})  \notag \\
&=&\theta _{L}^{(2)}(\gamma (F^{\epsilon }(\tau _{1},\tau _{2},\epsilon
)))\det (C\Omega +D)^{-l/2}\qquad (\text{Theorem \ref{TheoremGequiv}})
\notag \\
&=&\theta _{L}^{(2)}(\gamma \Omega )\det (C\Omega +D)^{-l/2}\ ]\qquad
\ \ \ \ \ \ \ \ \ \ \ \ \   (\text{%
(\ref{Omegaconvention})})  \notag \\
&=&\theta _{L}^{(2)}(\Omega )|_{l/2}\ \gamma .  \label{wtlaction}
\end{eqnarray}

\medskip For example, if the lattice $L$ is \emph{unimodular} as well as
even then $\theta _{L}^{(2)}$ is a Siegel modular form of weight $l/2$ on
the full group $Sp(4,\mathbb{Z})$. Then (\ref{wtlaction}) informs us that
\begin{equation*}
\hat{Z}_{V_{L},\epsilon }^{(2)}(\tau _{1},\tau _{2},\epsilon )|_{l/2}\gamma =%
\hat{Z}_{V_{L},\epsilon }^{(2)}(\tau _{1},\tau _{2},\epsilon ),\ \gamma \in
G,
\end{equation*}%
i.e. $\hat{Z}_{V_{L},\epsilon }^{(2)}(\tau _{1},\tau _{2},\epsilon )$ is
automorphic of weight $l/2$ with respect to the group $G$.

\section{Genus two partition function for free bosons in the $\protect\rho $%
-formalism}

\subsection{The genus two partition function $Z_{V}^{(2)}(\protect\tau ,w,%
\protect\rho )$}

In this and the next Section we carry out an analysis of the partition
function $Z_{V}^{(2)}(\tau ,w,\rho )$ for Heisenberg free bosonic and
lattice VOAs $V$ in the $\rho $-formalism, which obtains when we consider
sewing a twice-punctured torus to itself as reviewed in Section 2.3. The
main results mirror those obtained in the $\epsilon $-formalism. Recall that
our definition (\ref{Z2_def_rho}) of the partition function in this case is

\begin{equation}
Z_{V,\rho }^{(2)}(\tau ,w,\rho )=Z_{V}^{(2)}(\tau ,w,\rho )=\sum_{n\geq
0}\sum_{u\in V_{[n]}}Z_{V}^{(1)}(u,\bar{u},w,\tau )\rho ^{n}.
\label{eq: part func II}
\end{equation}%
Here, as before, $u$ ranges over a basis of $V_{[n]}$ and $\bar{u}$ is the dual
state with respect to the square-bracket LiZ metric.
$Z^{(1)}(u,v,w,\tau )$ is the
genus one $2$-point function (\ref{Z1_2pt}). Again (\ref{eq: part func II})
is multiplicative over $V$ and independent of the choice of basis, and as
long as we can choose a diagonal basis $\{u\}$ (e.g. the theory $M$ for a
free boson) then
\begin{equation}
Z_{V}^{(2)}(\tau ,w,\rho )=\sum_{n\geq 0}\sum_{u\in V_{[n]}}\frac{%
Z_{V}^{(1)}(u,u,w,\tau )}{\langle u,u\rangle _{\mathrm{sq}}}\rho ^{n}.
\label{eq: part func diag II}
\end{equation}%
For one free boson, we choose the same diagonal basis (\ref{eq: sq fock vec}%
) as before, so that these states $v=v(\lambda )$ may be indexed by
unrestricted partitions $\lambda $. We have seen that $\lambda $ determines
a labelled set $\Phi _{\lambda }$; what is also important in the present
context is another labelled set $\Phi _{\lambda ,2}$ defined to be the
disjoint union of two copies of $\Phi _{\lambda }$, call them $\Phi
_{\lambda }^{(1)}$, $\Phi _{\lambda }^{(2)}$. Fix an identification $\iota
:\Phi _{\lambda }^{(1)}\leftrightarrow \Phi _{\lambda }^{(2)}$.

\medskip The $2$-point function $Z_{M}^{(1)}(v(\lambda ),v(\lambda ),w,\tau
) $ was explicitly described in (\cite{MT1}, Corollary 1), where it was
denoted $F_{M}(v,w_{1},v,w_{2};\tau )$. In what follows we continue to use
the notation $w_{12}=w_{1}-w_{2}=w,w_{21}=w_{2}-w_{1}=-w$. Then we have
(loc. cit.)
\begin{equation}
Z_{M}^{(1)}(v(\lambda ),v(\lambda ),w,\tau )=Z_{M}^{(1)}(\tau )\sum_{\phi
\in F(\Phi _{\lambda ,2})}\Gamma (\phi ),  \label{eq: genus1 2-point func}
\end{equation}%
with
\begin{equation}
\Gamma (\phi ,w,\tau )=\Gamma (\phi )=\prod_{\{r,s\}}\xi (r,s,w,\tau ).
\label{eq: Gamma II}
\end{equation}%
Furthermore, $\phi $ ranges over the elements of $F(\Phi _{\lambda ,2})$
(fixed-point-free involutions in $\Sigma (\Phi _{\lambda ,2})$), $\{r,s\}$
ranges over the orbits of $\phi $ on $\Phi _{\lambda ,2}$, and
\begin{equation*}
\xi (r,s,w,\tau )=\left\{
\begin{array}{ll}
C(r,s,\tau ), & \mbox{if}\ \{r,s\}\subseteq \Phi _{\lambda }^{(i)},i=1\ %
\mbox{or}\ 2, \\
D(r,s,w_{ij},\tau ) & \mbox{if}\ r\in \Phi _{\lambda }^{(i)},s\in \Phi
_{\lambda }^{(j)},i\neq j.%
\end{array}%
\right.
\end{equation*}

\begin{remark}
\label{Rem_Dsym} Note that $\xi $ is well-defined since $D(r,s,w_{ij},\tau
)=D(s,r,w_{ji},\tau )$.
\end{remark}

As usual, the fixed-point-free involution $\phi $ defines a complete
matching $\mu _{\phi }$ on the underlying labelled set $\Phi _{\lambda ,2}$.
In the $\epsilon$-formalism we dealt with a pair of complete matchings based
on $\Phi _{\lambda }$; now we have just a single fixed-point-free involution
acting on two copies of $\Phi _{\lambda }$. However, we can supplement $\mu
_{\phi }$ with the `canonical' matching defined by $\iota $, which can also
be considered as a fixed-point-free involution. Graphically we may represent
$\mu _{\iota }$ by a broken line, so that the analog of figure 5 is
\begin{equation*}
\begin{array}{ccccc}
\overset{r_{1}}{\bullet } & \overset{\phi }{\longrightarrow } & \overset{%
s_{1}}{\bullet } & \overset{\iota }{--\rightarrow } & \overset{s_{1}^{\prime
}}{\bullet } \\
\overset{r_{2}}{\bullet } & \overset{\phi }{\longrightarrow } & \overset{%
s_{2}}{\bullet } & \overset{\iota }{--\rightarrow } & \overset{s_{2}^{\prime
}}{\bullet } \\
\vdots &  & \vdots & \vdots & \vdots \\
\overset{r_{b}}{\bullet } & \overset{\phi }{\longrightarrow } & \overset{%
s_{b}}{\bullet } & \overset{\iota }{--\rightarrow } & \overset{s_{b}^{\prime
}}{\bullet }%
\end{array}%
\end{equation*}

\begin{center}
{\small {Fig. 7} }
\end{center}

\noindent where $\Phi _{\lambda ,2}=\{r_{1},...,r_{b}\}=\{s_{1},...,s_{b}\}$
and $\iota :s\mapsto s^{\prime }$ is the canonical label-preserving
identification of the two copies of $\Phi _{\lambda }$.

\medskip We create analogs of the chequered diagrams of Section 6.1, where
now the cycles correspond to orbits of the cyclic group $\langle \iota \phi
\rangle $ and the edges are not alternately labelled by integers $1,2$ but
rather by solid and broken lines corresponding to the action of the
involutions $\phi $ and $\iota $ respectively. In addition to a positive
integer, each node carries a label $a\in \{1,2\}$ to indicate that it
belongs to $\Phi _{\lambda }^{(a)}$, where the canonical involution
satisfies $\iota :a\mapsto \bar{a}$ with convention (\ref{abar}). We call
such objects \emph{doubly-indexed chequered cycles (diagrams)}. Thus a
doubly-indexed chequered cycle with four nodes looks as follows:

\begin{center}
\begin{picture}(250,60)


\put(100,50){\line(1,0){4}}
\put(108,50){\line(1,0){4}}
\put(116,50){\line(1,0){4}}
\put(85,55){\makebox(0,0){$i_1,a_1$}}
\put(100,50){\circle*{4}}

\put(120,50){\line(0,-1){20}}
\put(135,55){\makebox(0,0){$i_1,\bar{a}_1$}}
\put(120,50){\circle*{4}}

\put(100,30){\line(1,0){4}}\put(108,30){\line(1,0){4}}\put(116,30){\line(1,0){4}}
\put(135,25){\makebox(0,0){$i_2,a_2$}}
\put(120,30){\circle*{4}}

\put(100,50){\line(0,-1){20}}
\put(85,25){\makebox(0,0){$i_2,\bar{a}_2$}}
\put(100,30){\circle*{4}}

\end{picture}

{\small Fig. 8 A Doubly-Indexed Chequered Cycle}
\end{center}

We have no need to consider all permutations of $\Phi _{\lambda ,2}$; the
only ones of relevance are those that commute with $\iota $ and preserve
both $\Phi _{\lambda }^{(1)}$ and $\Phi _{\lambda }^{(2)}$. We denote this
group, which is plainly isomorphic to $\Sigma (\Phi _{\lambda })$, by $%
\Delta _{\lambda }$. By definition, an automorphism of a chequered diagram $%
D $ in the above sense is an element of $\Delta _{\lambda }$ which preserves
edges and node labels.

\medskip For a chequered diagram $D$ corresponding to the partition $\lambda
=\{1^{e_{1}}...p^{e_{p}}\}$ we set
\begin{equation}
\gamma (D)=\frac{\prod_{\{k,l\}}\xi (k,l,w,\tau )}{\prod (-i)^{e_{i}}}\rho
^{\sum ie_{i}}  \label{eq: diagram weight II}
\end{equation}%
where $\{k,l\}$ ranges over the solid edges of $D$. We now have all the
pieces assembled to copy the arguments of Section 6. First use the group $%
\Delta _{\lambda }$ in place of $\Sigma (\Phi _{\lambda })$ to get the
analog of Proposition \ref{Prop_Z2boson_cheq}, namely
\begin{equation}
Z_{M}^{(2)}(\tau ,w,\rho )=Z_{M}^{(1)}(\tau )\sum_{D}\frac{\gamma (D)}{|%
\mathrm{Aut}(D)|},  \label{eq: Z2Mrho}
\end{equation}%
the sum ranging over all chequered diagrams. The further analysis of this
equality again proceeds as before, and we find that
\begin{equation}
Z_{M}^{(2)}(\tau ,w,\rho )=Z_{M}^{(1)}(\tau )\prod_{\mathcal{R}}(1-\gamma
(L))^{-1/2}  \label{eq: prodII}
\end{equation}%
where $L$ ranges over (oriented) rotationless chequered cycles $\mathcal{R}$.

\medskip Now introduce the infinite matrix $R^{\vee }=(R_{ab}^{\vee }(k,l))$%
, $k,l\geq 1,a,b\in \{1,2\}$, given by the block matrix%
\begin{equation*}
R^{\vee }(k,l)=-\frac{\rho ^{(k+l)/2}}{\sqrt{kl}}\left[
\begin{array}{cc}
C(k,l,\tau ) & D(k,l,\tau ,w) \\
D(l,k,\tau ,w) & C(k,l,\tau )%
\end{array}%
\right] .
\end{equation*}%
Compared to (\ref{Rdef}) then, we have $R_{ab}^{\vee }(k,l)=R_{a\bar{b}%
}(k,l) $.

\medskip The analog of the weight function $\omega $ (cf. (\ref{eq: eqomega}%
)) is as follows: for a doubly-indexed chequered diagram $D$ consisting of
solid and broken edges as above we set $\omega ^{\vee }(D)=\prod \omega
^{\vee }(E)$, the product running over all solid edges. Moreover for a solid
edge $E$ of type $\overset{k,a}{\bullet }\longrightarrow \overset{l,b}{%
\bullet }$ with nodes $k,l$ lying in $\Phi _{\lambda }^{(a)},\Phi _{\lambda
}^{(b)}$ respectively, we set
\begin{equation*}
\omega ^{\vee }(E)=R_{a,b}^{\vee }(k,l)=\mbox{$(a,b)$-entry of}\ R^{\vee
}(k,l).
\end{equation*}%
The analog of Lemma \ref{Lemma_om_gamma} is

\begin{lemma}
\label{Lemma_om_gamma_R}$\omega ^{\vee }(L)=\gamma (L)$.
\end{lemma}

\noindent \textbf{Proof. }From (\ref{eq: diagram weight II}) it follows that
for a chequered cycle $L$ we have
\begin{equation}
\gamma (L)=\prod_{\{k,l\}}-\frac{\xi (k,l,w,\tau )\rho ^{(k+l)/2}}{\sqrt{kl}}
,  \label{eq: gammaprod II}
\end{equation}%
the product ranging over solid edges $\{k,l\}$ of $L$. This follows because
in (\ref{eq: diagram weight II}) each part of the relevant partition $%
\lambda $ occurs twice. So to prove the Lemma it suffices to show that if
the nodes $k,l$ lie in $\Phi _{\lambda }^{(a)},\Phi _{\lambda }^{(b)}$
respectively then the $(a,b)$-entry of $R^{\vee }(k,l)$ coincides with the
corresponding factor of (\ref{eq: gammaprod II}). This follows from our
previous discussion together with Remark \ref{Rem_Dsym}. \ \ \ \ \ \ $%
\square $

\medskip From Lemma \ref{Lemma_om_gamma_R} and (\ref{eq: prodII}) we obtain
the analogs of Theorems \ref{Theorem_Z2_boson} and \ref%
{Theorem_Z2_boson_prod}, namely

\begin{theorem}
\label{Theorem_Z2_boson_rho_hat} Let $M$ be the VOA of one free boson. Then
\begin{equation}
Z_{M}^{(2)}(\tau ,w,\rho )= \frac{Z_{M}^{(1)}(\tau )}{\det (1-R^{\vee
})^{1/2}}.  \label{Z2_1bos_rho_hat}
\end{equation}
\end{theorem}

\begin{theorem}
\label{Theorem_Z2_boson_prod_rho_hat}Let $M$ be the VOA of one free boson.
Then
\begin{equation}
Z_{M}^{(2)}(\tau ,w,\rho ) = \frac{Z_{M}^{(1)}(\tau )}{\prod_{\mathcal{R}%
}(1-\omega ^{\vee }(L))^{1/2}}  \label{eq: prodIII}
\end{equation}%
where $L$ ranges over (oriented) rotationless chequered cycles $\mathcal{R}$.
\end{theorem}

We can re-formulate these results in a slightly different way using diagrams
in which the canonical involution $\iota $ and broken edges are dispensed
with, and utilizing a weight function $\omega $ associated to the matrix $R$
rather than $R^{\vee }$. First note that for a doubly-indexed chequered
cycle $L$ we have for suitable indices that
\begin{equation*}
\omega ^{\vee }(L)=R_{a_{1}a_{2}}^{\vee }(k_{1},k_{2})R_{\bar{a}%
_{2}a_{3}}^{\vee }(k_{2},k_{3})R_{\bar{a}_{3}a_{4}}^{\vee
}(k_{3},k_{4})\ldots R_{\bar{a}_{d}\bar{a}_{1}}^{\vee }(k_{d},k_{1}).
\end{equation*}%
Since $R_{ab}^{\vee }(k,l)=R_{a\bar{b}}(k,l)$ then
\begin{equation}
\omega ^{\vee }(L)=R_{a_{1}\bar{a}_{2}}(k_{1},k_{2})R_{\bar{a}_{2}\bar{a}%
_{3}}(k_{2},k_{3})R_{\bar{a}_{3}\bar{a}_{4}}(k_{3},k_{4})\ldots R_{\bar{a}%
_{d}a_{1}}(k_{d},k_{1}).  \label{eq: altomegavee}
\end{equation}%
It is apparent that we can interpret (\ref{eq: altomegavee}) as a weight
function evaluated on the doubly-indexed cycles introduced in Section 3.3
with only regular (solid) edges and nodes indexed as before by a pair of
positive integers $k,a$ with $a\in \{1,2\}$. Thus, a typical doubly-indexed
cycle looks as follows:

\begin{center}
\begin{picture}(300,80)


\put(100,50){\line(1,2){10}}
\put(82,50){\makebox(0,0){$k_1,a_1$}}
\put(100,50){\circle*{4}}

\put(110,70){\line(1,0){20}}
\put(100,78){\makebox(0,0){$k_2,a_2$}}
\put(110,70){\circle*{4}}

\put(130,70){\line(1,-2){10}}
\put(145,78){\makebox(0,0){$k_3,a_3$}}
\put(130,70){\circle*{4}}

\put(140,50){\line(-1,-2){10}}
\put(160,50){\makebox(0,0){$k_4,a_4$}}
\put(140,50){\circle*{4}}

\put(110,30){\line(1,0){20}}
\put(145,20){\makebox(0,0){$k_5,a_5$}}
\put(130,30){\circle*{4}}

\put(100,50){\line(1,-2){10}}
\put(100,20){\makebox(0,0){$k_6,a_6$}}
\put(110,30){\circle*{4}}

\end{picture}

{\small Fig. 9 Doubly-Indexed Cycle}
\end{center}

If a doubly-indexed cycle has $d$ nodes with consecutive edges labeled $%
(k_{1},k_{2})$, $(k_{2},k_{3})$, ...$(k_{d},k_{1})$ and corresponding second
node indices $a_{1},...,a_{d}\in \{1,2\}$, we set
\begin{equation}
\omega (L)=R_{a_{1}a_{2}}(k_{1},k_{2})R_{a_{2}a_{3}}(k_{2},k_{3})\ldots
R_{a_{d}a_{1}}(k_{d},k_{1}).  \label{eq: omegatoomegavee}
\end{equation}%
Our discussion shows that we may now reformulate Theorems \ref%
{Theorem_Z2_boson_rho_hat} and \ref{Theorem_Z2_boson_prod_rho_hat} as
follows:

\begin{theorem}
\label{Theorem_Z2_boson_rho} Let $M$ be the vertex operator algebra of one
free boson. Then
\begin{equation}
Z_{M}^{(2)}(\tau ,w,\rho )=\frac{Z_{M}^{(1)}(\tau )}{\det (1-R)^{1/2}}.
\label{Z2_1bos_rho}
\end{equation}
\end{theorem}

\begin{theorem}
\label{Theorem_Z2_boson_prod_rho}Let $M$ be the vertex operator algebra of
one free boson. Then
\begin{equation}
Z_{M}^{(2)}(\tau ,w,\rho )= \frac{Z_{M}^{(1)}(\tau )}{\prod_{L}(1-\omega
(L))^{1/2}}  \label{eq: prodIV}
\end{equation}%
where $L$ ranges over (oriented) rotationless doubly-indexed cycles. $\ \ \
\ \ \square $
\end{theorem}

\subsection{Holomorphic and modular-invariance properties}

In Section 2.3 we reviewed the genus two $\rho $-sewing formalism and
introduced the domain $\mathcal{D}^{\rho }$ which parameterizes the genus
two surface. An immediate consequence of Theorem \ref{Theorem_R} is the
following.

\begin{theorem}
\label{Theorem_Z2_boson_rho_hol}$Z_{M}^{(2)}(\tau ,w,\rho )$ is holomorphic
in $\mathcal{D}^{\rho }$. \ \ \ \ \ \ $\square $
\end{theorem}

We next consider the invariance properties of the genus two partition
function with respect to the action of the $\mathcal{D}^{\rho }$-preserving
group $\Gamma _{1}$ reviewed in Section 2.3. As in Section 6.2, we again
define a modular genus two partition function for two free bosons by
\begin{eqnarray}
Z_{M^{2},\mathrm{mod}}^{(2)}(\tau ,w,\rho ) &=&Z_{M^{2}}^{(2)}(\tau ,w,\rho
)q^{-1/12}  \notag \\
&=&\frac{1}{\eta (\tau )^{2}\det (1-R)}.  \label{Z2mod_rho}
\end{eqnarray}%
We then have a natural analog of Theorem \ref{Theorem_Z2_G}:

\begin{theorem}
\label{Theorem_Z2_rho_G1}If $\gamma \in \Gamma _{1}$ then
\begin{equation*}
Z_{M^{2},\mathrm{mod}}^{(2)}(\tau ,w,\rho )|_{-1}\gamma =\chi (\gamma
)Z_{M^{2},\mathrm{mod}}^{(2)}(\tau ,w,\rho ).
\end{equation*}
\end{theorem}

\begin{corollary}
\label{Cor_Z2_24_eps_G1}If $\gamma \in \Gamma _{1}$ with $Z_{M^{24},\mathrm{%
mod}}^{(2)}=(Z_{M^{2},\mathrm{mod}}^{(2)})^{12}$ then
\begin{equation*}
Z_{M^{24},\mathrm{mod}}^{(2)}(\tau ,w,\rho )|_{-12}\gamma =Z_{M^{24},\mathrm{%
mod}}^{(2)}(\tau ,w,\rho ).
\end{equation*}
\end{corollary}

\noindent \textbf{Proof. }We have to show that
\begin{equation}
Z_{M^{2},\mathrm{mod}}^{(2)}(\gamma .(\tau ,w,\rho ))\det (C\Omega +D)=\chi
(\gamma )Z_{M^{2},\mathrm{mod}}^{(2)}(\tau ,w,\rho )
\label{eq: Z^2_rho_identity}
\end{equation}%
for $\gamma \in \Gamma _{1}$ where $\det (C\Omega _{11}+D)=c_{1}\Omega
_{11}+d_{1}$. The proof is similar to that of Theorem \ref{Theorem_Z2_G}.
Consider the determinant formula (\ref{Z2_1bos_rho}). For $\gamma \in \Gamma
_{1}$ define
\begin{equation*}
R_{ab}^{\prime }(k,l,\tau ,w,\rho )=R_{ab}(k,l,\frac{a_{1}\tau +b_{1}}{%
c_{1}\tau +d_{1}},\frac{w}{c_{1}\tau +d_{1}},\frac{\rho }{(c_{1}\tau
+d_{1})^{2}})
\end{equation*}%
following (\ref{gam1_rho}). We find from Section 6.3 of \cite{MT2} that
\begin{eqnarray*}
1-R^{\prime } &=&1-R-\kappa \Delta \\
&=&(1-\kappa S).(1-R),
\end{eqnarray*}%
where
\begin{eqnarray*}
\Delta _{ab}(k,l) &=&\delta _{k1}\delta _{l1}, \\
\kappa &=&\frac{\rho }{2\pi i}\frac{c_{1}}{c_{1}\tau +d_{1}}, \\
S_{ab}(k,l) &=&\delta _{k1}\sum_{c\in \{1,2\}}((1-R)^{-1})_{cb}(1,l).
\end{eqnarray*}%
Since $\det (1-R)$ and $\det (1-R^{\prime })$ are convergent on $\mathcal{D}%
^{\rho }$ we find%
\begin{equation*}
\det (1-R^{\prime })=\det (1-\kappa S).\det (1-R).
\end{equation*}%
Indexing the columns and rows by $(a,k)=(1,1),(2,1),\ldots (1,k),(2,k)\ldots
$ and noting that $S_{1b}(k,l)=S_{2b}(k,l)$ we find that
\begin{eqnarray*}
\det (1-\kappa S) &=&\left\vert
\begin{array}{cccc}
1-\kappa S_{11}(1,1) & -\kappa S_{12}(1,1) & -\kappa S_{11}(1,2) & \cdots \\
-\kappa S_{11}(1,1) & 1-\kappa S_{12}(1,1) & -\kappa S_{11}(1,2) & \cdots \\
0 & 0 & 1 & \cdots \\
\vdots & \vdots & \vdots & \ddots%
\end{array}%
\right\vert \\
&=&1-\kappa S_{11}(1,1)-\kappa S_{12}(1,1), \\
&=&1-\kappa \sigma ((1-R)^{-1})(1,1)).
\end{eqnarray*}%
Here and below, $\sigma (M)$ denotes the sum of the entries of a (finite)
matrix $M$. Applying (\ref{Om11rho}), it is clear that
\begin{equation*}
\det (1-\kappa S)=\frac{c_{1}\Omega _{11}+d_{1}}{c_{1}\tau +d_{1}}.
\end{equation*}%
The Theorem then follows from (\ref{Z1modgam}). \ \ \ \ \ \ \ $\square $

\begin{remark}
\label{Rem_Drhohat_mod} $Z_{M^{2},\mathrm{mod}}^{(2)}(\tau ,w,\rho )$ can be
trivially considered as function on the covering space $\mathcal{\hat{D}}%
^{\rho }$ discussed in Section 6.3 of ref. \cite{MT2}. Then $Z_{M^{2},%
\mathrm{mod}}^{(2)}(\tau ,w,\rho )$ is modular with respect to $L=\hat{H}%
\Gamma _{1}$ with trivial invariance under the action of the Heisenberg
group $\hat{H}$ (op. cite.).
\end{remark}

\section{Genus two partition function for lattice theories in the $\protect%
\rho$-formalism}

Let $L$ be an even lattice of rank $l$. The underlying Fock space for the
corresponding lattice theory is given by (\ref{lattice fock space}), and the
general shape of the corresponding partition function is as in (\ref{eq:
part func II}). We utilize the same basis for Fock space as before, namely $%
\{u\otimes e^{\alpha }\}$ where $\alpha $ ranges over $L$ and $u$ ranges
over the usual orthogonal basis for $M^{l}$. From Lemma \ref%
{Lemma_LiZ_lattice} and Corollary \ref{Corollary_cocycle_choice} we see that
\begin{gather}
Z_{V_{L}}^{(2)}(\tau ,w,\rho )=  \notag \\
\sum_{\alpha ,\beta \in L}\sum_{n\geq 0}\sum_{u\in M_{[n]}^{l}}\frac{%
Z_{M^{l}\otimes e^{\beta }}^{(1)}(u\otimes e^{\alpha },u\otimes e^{-\alpha
},w,\tau )}{\langle u,u\rangle _{\mathrm{sq}}}\rho ^{n+(\alpha ,\alpha )/2}.
\label{eq: rk1rholattpartfunc}
\end{gather}%
The general shape of the $2$-point function occurring in (\ref{eq:
rk1rholattpartfunc}) is discussed extensively in \cite{MT1}. By Proposition
1 (loc. cit.) it splits as a product
\begin{gather}
Z_{M^{l}\otimes e^{\beta }}^{(1)}(u\otimes e^{\alpha },u\otimes e^{-\alpha
},w,\tau )=  \notag \\
Q_{M^{l}\otimes e^{\beta }}^{\alpha }(u,u,w,\tau )Z_{M^{l}\otimes e^{\beta
}}^{(1)}(e^{\alpha },e^{-\alpha },w,\tau ),  \label{eq: rhosplitform}
\end{gather}%
where we have identified $e^{\alpha }$ with $\mathbf{1}\otimes e^{\alpha }$,
and where $Q_{M^{l}\otimes e^{\beta }}^{\alpha }$ is a function\footnote{%
Note: in \cite{MT1} the functional dependence on $\alpha $, here denoted by
a superscript, was omitted.} that we will shortly discuss in greater detail.
In \cite{MT1}, Corollary 5 (cf. the Appendix to the present paper) we
established also that
\begin{equation}
Z_{M^{l}\otimes e^{\beta }}^{(1)}(e^{\alpha },e^{-\alpha },w,\tau )=\epsilon
(\alpha ,-\alpha )q^{(\beta ,\beta )/2}\frac{\exp ((\beta ,\alpha )w)}{%
K(w,\tau )^{(\alpha ,\alpha )}}Z_{M^{l}}^{(1)}(\tau ),
\label{eq: latt2ptform}
\end{equation}%
where as usual we are taking $w$ in place of $z_{12}=z_{1}-z_{2}$. With
cocycle choice $\epsilon (\alpha ,-\alpha )=(-1)^{(\alpha ,\alpha )/2}$ (cf.
Appendix) we may then rewrite (\ref{eq: rk1rholattpartfunc}) as
\begin{eqnarray}
&&\hspace{4.5cm}Z_{V_{L}}^{(2)}(\tau ,w,\rho )=  \notag \\
&&Z_{M^{l}}^{(1)}(\tau )\sum_{\alpha ,\beta \in L}\exp \{\pi i[(\beta ,\beta
)\tau +2(\alpha ,\beta )\frac{w}{2\pi i}+\frac{(\alpha ,\alpha )}{2\pi i}%
\log (\frac{-\rho }{K(w,\tau )^{2}})]\}  \notag \\
&&\hspace{3cm}\ \sum_{n\geq 0}\sum_{u\in M_{[n]}^{l}}\frac{Q_{M^{l}\otimes
e^{\beta }}^{\alpha }(u,u,w,\tau )}{\langle u,u\rangle _{\mathrm{sq}}}\rho
^{n}.  \label{eq: newrholattpartfunc}
\end{eqnarray}%
We note that this expression is, as it should be, independent of the choice
of branch for the logarithm function. We are going to establish the \emph{%
precise} analog of Theorem \ref{Theorem_Z2_L_eps}, to wit:

\begin{theorem}
\label{Theorem_Z2_L_rho}We have
\begin{equation*}
Z_{V_{L}}^{(2)}(\tau ,w,\rho )=Z_{M^{l}}^{(2)}(\tau ,w,\rho )\theta
_{L}^{(2)}(\Omega ).
\end{equation*}
\end{theorem}

As in the case of the $\epsilon$-formalism, we first handle the case of
rank $1$ lattices, then consider the general case. The inner double sum in (%
\ref{eq: newrholattpartfunc}) is the object which requires attention, and we
can begin to deal with it along the lines of previous Sections. Namely,
arguments that we have already used several times show that the double sum
may be written in the form
\begin{equation*}
\sum_{D}\frac{\gamma (D)}{|\mathrm{Aut}(D)|}=\exp \left( \frac{1}{2}%
\sum_{N}\gamma (N)\right) .
\end{equation*}%
Here, $D$ ranges over the doubly indexed chequered diagrams of Section 8,
while $N$ ranges over oriented, doubly indexed chequered diagrams which are
connected. Leaving aside the definition of $\gamma (N)$ for now, we
recognize as before that the piece involving only connected diagrams with no
end nodes (aka doubly indexed chequered cycles) splits off as a factor.
Apart from a $Z_{M}^{(1)}(\tau )$ term this factor is, of course, precisely
the expression (\ref{eq: Z2Mrho}) for the free boson. With these
observations, we see from (\ref{eq: newrholattpartfunc}) that the following
holds:
\begin{eqnarray}
&&\hspace{3cm}\frac{Z_{V_{L}}^{(2)}(\tau ,w,\rho )}{Z_{M}^{(2)}(\tau ,w,\rho
)}=  \notag \\
&&\sum_{\alpha ,\beta \in L}\exp \{\pi i[(\beta ,\beta )\tau +2(\alpha
,\beta )\frac{w}{2\pi i}+\frac{(\alpha ,\alpha )}{2\pi i}\log (\frac{-\rho }{%
K(w,\tau )^{2}})  \notag \\
&&\hspace{2.5cm}+\frac{1}{2\pi i}\sum_{N_{\alpha ,\beta }}\gamma (N_{\alpha
,\beta })]\}.  \label{eq: lattrhopartqnt}
\end{eqnarray}%
The expression (\ref{eq: g2omega}) for $\theta _{L}^{(2)}(\Omega )$ is
invariant under interchange of $\alpha $ and $\beta $ as well as replacing $%
(\alpha ,\beta )$ by $-(\alpha ,\beta )$, as noted in Remark \ref%
{Rem_Siegel_Theta}. So to prove Theorem \ref{Theorem_Z2_L_rho}, we see from (%
\ref{eq: lattrhopartqnt}) that it is sufficient to establish that for each
pair of lattice elements $\alpha ,\beta \in L$, we have
\begin{eqnarray}
&&\hspace{1cm}(\beta ,\beta )\Omega _{11}+2(\alpha ,\beta )\Omega
_{12}+(\alpha ,\alpha )\Omega _{22}=  \notag \\
&&(\beta ,\beta )\tau +2(\alpha ,\beta )\frac{w}{2\pi i}+\frac{(\alpha
,\alpha )}{2\pi i}\log (\frac{-\rho }{K(w,\tau )^{2}})  \notag \\
&&\hspace{2.5cm}+\frac{1}{2\pi i}\sum_{N_{\alpha ,\beta }}\gamma (N_{\alpha
,\beta }).  \label{eq: rhoomegas}
\end{eqnarray}

Recall the formulas for the $\Omega _{ij}$  (\ref{Om11rho}-\ref{Om22rho}%
). From now on we fix a pair of lattice elements $\alpha ,\beta $. In order
to reconcile (\ref{eq: rhoomegas}) with the formulas for the $\Omega _{ij}$,
we must carefully consider the expression $\sum_{N_{\alpha ,\beta }}\gamma
(N_{\alpha ,\beta })$. To begin with, the $N_{\alpha ,\beta }$ here
essentially range over oriented chequered necklaces as used in Sections 3
and 8, except that in the present case the integer labelling of the end
nodes is unrestricted (heretofore, it was required to be $1$). The function $%
\gamma $ is essentially (\ref{eq: diagram weight II}), except that we also
get contributions from the end nodes which are now present. Suppose that an
end node has label $k\in \Phi ^{(a)},a\in \{1,2\}$. Then according to
Proposition 1 and display (45) of \cite{MT1} (cf. the Appendix to the
present paper), the contribution of the end node is equal to
\begin{gather}
\xi _{\alpha ,\beta }(k,a,\tau ,w,\rho )=  \notag \\
\left\{
\begin{array}{ll}
\frac{\rho ^{k/2}}{\sqrt{k}}(a,\delta _{k,1}\beta +C(k,0,\tau )\alpha
-D(k,0,w,\tau )\alpha ), & a=1 \\
\frac{\rho ^{k/2}}{\sqrt{k}}(a,\delta _{k,1}\beta -C(k,0,\tau )\alpha
+D(k,0,-w,\tau )\alpha ), & a=2%
\end{array}%
\right.  \label{eq: xirho}
\end{gather}%
together with a contribution arising from the $-1$ in the denominator of (%
\ref{eq: diagram weight II}) (we will come back to this point later). Using
(cf. \cite{MT1}, displays (6), (11) and (12))
\begin{eqnarray*}
D(k,0,-w,\tau ) &=&(-1)^{k+1}P_{k}(-w,\tau )=-P_{k}(w,\tau ), \\
C(k,0,\tau ) &=&(-1)^{k+1}E_{k}(\tau ),
\end{eqnarray*}%
we can combine the two possibilities in (\ref{eq: xirho}) thus (recalling
that $E_{k}=0$ for odd $k$):
\begin{gather}
\xi _{\alpha ,\beta }(k,a,\tau ,w,\rho )=  \notag \\
-\frac{\rho ^{k/2}}{\sqrt{k}}(a,\delta _{k,1}\beta
+(-1)^{(k+1)a}[P_{k}(w,\tau )-E_{k}(\tau )]\alpha ).
\label{eq: xialphabetaform}
\end{gather}

Now consider an oriented chequered necklace $N=N_{\alpha ,\beta }$ with a
pair of end nodes labelled by $k\in \Phi ^{(a)}$ and $l\in \Phi ^{(b)}$
respectively. It looks like
\begin{equation*}
\overset{k,a}{\bullet }\longrightarrow \bullet -\rightarrow
-.....-\rightarrow -\bullet \longrightarrow \overset{l,b}{\bullet }
\end{equation*}

\begin{center}
{\small {Fig. 10} }
\end{center}

\noindent From (\ref{eq: xialphabetaform}) we see that the total
contribution of the end nodes to $\gamma (N)$ is
\begin{eqnarray*}
&&-\xi _{\alpha ,\beta }(k,a,\tau ,w,\rho )\xi _{\alpha ,\beta }(l,b,\tau
,w,\rho )= \\
&&-\frac{\rho ^{(k+l)/2}}{\sqrt{kl}}(a,\delta _{k,1}\beta
+(-1)^{(k+1)a}[P_{k}(w,\tau )-E_{k}(\tau )]\alpha ) \\
&& (a,\delta _{l,1}\beta +(-1)^{(l+1)b}[P_{l}(w,\tau )-E_{l}(\tau )]\alpha ),
\end{eqnarray*}%
where we note that a sign $-1$ arises from each \emph{pair} of nodes, as
follows from (\ref{eq: diagram weight II}). Furthermore, if $N_{\alpha
,\beta }^{\prime }$ denotes the oriented chequered necklace from which the
two end nodes and edges have been \emph{removed} (we refer to these as \emph{%
shortened} necklaces), then we have
\begin{eqnarray}
\gamma (N_{\alpha ,\beta }) &=&-\xi _{\alpha ,\beta }(k,a,\tau ,w,\rho )\xi
_{\alpha ,\beta }(l,b,\tau ,w,\rho )\gamma (N_{\alpha ,\beta }^{\prime })
\notag \\
&=&-\frac{\rho ^{(k+l)/2}}{\sqrt{kl}}\gamma (N_{\alpha ,\beta }^{\prime
})\{\delta _{k,1}\delta _{l,1}(\beta ,\beta )+  \notag \\
&&(-1)^{(k+1)a+(l+1)b}(\alpha ,\alpha )[P_{k}(w,\tau )-E_{k}(\tau
)][P_{l}(w,\tau )-E_{l}(\tau )]+  \notag \\
&&(\alpha ,\beta )[\delta _{k,1}(-1)^{(l+1)b}[P_{l}(w,\tau )-E_{l}(\tau )]
\notag \\
&&+\delta _{l,1}(-1)^{(k+1)a}[P_{k}(w,\tau )-E_{k}(\tau )]\}.
\label{eq: gammaNalphabetaform}
\end{eqnarray}

\medskip We will now consider the terms corresponding to $(\beta ,\beta
),(\alpha ,\beta )$ and $(\alpha ,\alpha )$ separately, and show that they
are precisely the corresponding terms on each side of (\ref{eq: rhoomegas}).
This will complete the proof of Theorem \ref{Theorem_Z2_L_rho} in
the case of rank $1$ lattices. From (\ref{eq: gammaNalphabetaform}), a $%
(\beta ,\beta )$ term arises only if the end node weights $k,l$ are both
equal to $1$, in which case we get
\begin{equation*}
-\rho \gamma (N_{\alpha ,\beta }^{\prime }).
\end{equation*}%
Then the total $(\beta ,\beta )$-contribution to the right-hand-side of (\ref%
{eq: rhoomegas}) is equal to
\begin{equation}
\tau -\frac{\rho }{2\pi i}\sum \gamma (N_{\alpha ,\beta }^{\prime }),
\label{Omega11seen}
\end{equation}%
where the sum ranges over shortened necklaces with end nodes having weight $%
1 $. However, from \cite{MT2} and also from Section 8 (esp. following (\ref%
{eq: prodII})) it follows that this sum is nothing else than $\sigma
((I-R)^{-1}(1,1))$. So we see from (\ref{eq: rhoomegas}) that (\ref%
{Omega11seen}) coincides with $\Omega _{11}$ (cf. (\ref{Om11rho})) as
required.

\medskip Next, from (\ref{eq: gammaNalphabetaform}) we see that an $(\alpha
,\beta )$-contribution arises whenever at least one of the end nodes has
label $1$. If the labels of the end nodes are unequal then the shortened
necklace with the \emph{opposite} orientation makes an equal contribution.
The upshot is that we may assume that the end node to the right of the
shortened necklace has label $l=1\in \Phi ^{(\bar{b})}$, as long as we count
accordingly. Then the contribution to the $(\alpha ,\beta )$-term from (\ref%
{eq: gammaNalphabetaform}) is equal to
\begin{equation*}
-2\rho ^{1/2}\sum_{k\geq 1}\frac{\rho ^{k/2}}{\sqrt{k}}%
(-1)^{(k+1)a+1}[P_{k}(w,\tau )-E_{k}(\tau )]\sum \gamma (N_{\alpha ,\beta
}^{\prime }),
\end{equation*}%
where the inner sum ranges over shortened necklaces with end nodes of weight
$1\in \Phi ^{(\bar{b})}$ and $k\in \Phi ^{(\bar{a})}$. In this case the
discussion in Section 8 shows that
\begin{eqnarray*}
&&\sum (-1)^{(k+1)a+1}\gamma (N_{\alpha ,\beta }^{\prime })=\sum
(-1)^{(k+1)a}R_{\bar{a},a_{1}}^{\vee }(k,k_{1})\ldots R_{a_{d-1},\bar{b}%
}^{\vee }(k_{d-1},1) \\
&=&\sum (-1)^{(k+1)(\bar{a}+1)+1}R_{\bar{a},\bar{a}_{1}}(k,k_{1})\ldots R_{%
\bar{a}_{d-1},b}(k_{d-1},1)
\end{eqnarray*}%
(summed over all indices $k_{1},..,k_{d-1}\geq 1$ and $a,a_{1},...,a_{d-1},b%
\in \{1,2\}$)
\begin{eqnarray*}
&=&\sum_{a,b=1}^{2}(-1)^{(k+1)(a+1)+1}(I-R)_{ab}^{-1}(k,1) \\
&=&\sum_{b=1}^{2}\{-(I-R)_{1b}^{-1}(k,1)+(-1)^{k}(I-R)_{2b}^{-1}(k,1)\} \\
&=&\sigma \{(-1,(-1)^{k})(I-R)^{-1}(k,1)\}.
\end{eqnarray*}%
So the contribution to the $(\alpha ,\beta )$-term from (\ref{eq:
gammaNalphabetaform}) is now seen to be equal to
\begin{eqnarray*}
&&-2\rho ^{1/2}\sum_{k\geq 1}\sigma \{\frac{\rho ^{k/2}}{\sqrt{k}}%
[P_{k}(w,\tau )-E_{k}(\tau )](-1,(-1)^{k})(I-R)^{-1}(k,1)\} \\
&=&-2\rho ^{1/2}\sigma ((b(I-R)^{-1})(1)),
\end{eqnarray*}%
where $b$ is defined by (\ref{bk}). Finally then, the total contribution to
the $(\alpha ,\beta )$ term on the right-hand-side of (\ref{eq: rhoomegas})
is found, using (\ref{Om12rho}), to be
\begin{equation*}
2\frac{w}{2\pi i}-2\frac{\rho ^{1/2}}{2\pi i}\sigma (b(I-R)^{-1}(1))=2\Omega
_{12},
\end{equation*}%
as required.

\medskip It remains to deal with the $(\alpha ,\alpha )$ term, the details
of which are very much along the lines as the case $(\alpha ,\beta )$ just
handled. A similar argument shows that the contribution to the $(\alpha
,\alpha )$-term from (\ref{eq: gammaNalphabetaform}) is equal to
\begin{equation*}
-b(I-R)^{-1}\bar{b}^{T},
\end{equation*}%
so that the total contribution to the $(\alpha ,\alpha )$ term on the
right-hand-side of (\ref{eq: rhoomegas}) is
\begin{equation*}
\frac{1}{2\pi i}\log (\frac{-\rho }{K(w,\tau )^{2}})-\frac{1}{2\pi i}%
b(I-R)^{-1}\bar{b}^{T}=\Omega _{22},
\end{equation*}%
as in (\ref{Om22rho}). This finally completes the proof of Theorem \ref%
{Theorem_Z2_L_rho} in the case of rank $1$ lattices. As for the general case
- we adopt the mercy rule by omitting details! The reader who has progressed
this far will have no difficulty in dealing with the general case, which
follows by combining the calculations in the rank $1$ case just completed
together with the method used in Section 7 to deal with a general lattice in
the $\epsilon $-formalism. \ \ \ \ \ \ \ $\square $

\medskip
Finally, the automorphic properties of $Z_{V_{L}}^{(2)}(\tau ,w,\rho )$ can
be analyzed much as in previous sections to find

\begin{theorem}
\label{Theorem_Z2_lattice_rho_hol}$Z_{V_{L}}^{(2)}(\tau ,w,\rho )$ is
holomorphic on\ the domain $\mathcal{D}^{\rho }$. \ \ \ \ \ $\square $
\end{theorem}

We again  introduce a normalized partition function
\begin{equation*}
\hat{Z}_{V_{L},\rho }^{(2)}(\tau ,w,\rho )=\frac{Z_{V_{L}}^{(2)}(\tau
,w,\rho )}{Z_{M^{l}}^{(2)}(\tau ,w,\rho )},
\end{equation*}%
so that by Theorem \ref{Theorem_Z2_L_rho} there is a commuting
diagram of holomorphic maps
\begin{equation*}
\begin{array}{lll}
\ \hspace{0cm}\mathcal{D}^{\rho } & \overset{F^{\rho }}{\longrightarrow } &
\hspace{0.7cm}\mathbb{H}_{2} \\
\hat{Z}_{V_{L},\rho }^{(2)}\searrow &  & \swarrow \theta _{L}^{(2)} \\
& \hspace{0.2cm}\mathbb{C} &
\end{array}%
\end{equation*}

Furthermore, the $\Gamma _{1}$-actions on the two functions are compatible
with
\begin{equation}
\hat{Z}_{V_{L},\rho }^{(2)}(\tau ,w,\rho )|_{l/2}\ \gamma =\theta
_{L}^{(2)}(\Omega )|_{l/2}\ \gamma ,  \label{wtlaction_rho}
\end{equation}%
for all $\gamma \in \Gamma _{1}$.

\section{Comparison of partition functions in the two formalisms}

In this section we consider the relationship between the genus two boson and
lattice partition functions computed in the $\epsilon $- and $\rho $-
formalisms of the previous Sections. Although, for a given VOA, the
partition functions enjoy many similar properties, we find that neither the
original nor the modular partition functions (cf. Sections 6.2 and 8.2) are
equal in the two formalisms. This result follows from an explicit
computation of the partition functions in the neighborhood of a two tori
degeneration point for two free bosons. It therefore follows that there is
likewise no equality between the partition functions in the $\epsilon $- and
$\rho $- formalisms for a lattice VOA.

\medskip As reviewed in Section 2, there exists a 1-1 map between
appropriate $\Gamma _{1}$-invariant neighborhoods of $\mathcal{D}^{\epsilon
} $ and $\mathcal{D}^{\chi }\cup \mathcal{D}_{0}^{\chi }$ about any two tori
degeneration point (cf. Theorem \ref{Theorem_epsrho_11map}). We will make
particular use of the explicit relationship between $(\tau _{1},\tau
_{2},\epsilon )$ and $(\tau ,w,\chi )$ to $O(w^{4})$ as follows:

\begin{proposition}
\label{Proposition_eps_rho_w4} For $(\tau ,w,\chi )\in \mathcal{D}^{\chi
}\cup \mathcal{D}_{0}^{\chi }$,
\begin{eqnarray*}
2\pi i\tau _{1} &=&2\pi i\tau +\frac{1}{12}(1-4\chi )w^{2}+\frac{1}{144}%
\,E_{2}\left( \tau \right) \left( 1-4\chi \right) ^{2}w^{4}+O(w^{6}), \\
2\pi i\tau _{2} &=&\log (f(\chi ))+\frac{1}{12}\,\left( 1-4\,\chi \right)
^{2}E_{4}\left( \tau \right) w^{4}+O(w^{6}), \\
\epsilon &=&-w\sqrt{1-4\chi }(1+(1-4\chi )E_{2}(\tau )w^{2})+O(w^{5}).
\end{eqnarray*}
\end{proposition}

A proof of Proposition \ref{Proposition_eps_rho_w4} appears in the Appendix.

\medskip Theorem \ref{Theorem_Degen_pt} tells us that $Z_{V,\epsilon
}^{(2)}(\tau _{1},\tau _{2},\epsilon )$ and $Z_{V,\rho }^{(2)}(\tau ,w,\rho
) $ agree at any two-tori degeneration point. In particular, for the VOA $%
M^{2} $ of two free bosons, we have%
\begin{gather}
\lim_{\rho ,w\rightarrow 0}Z_{M^{2},\rho }^{(2)}(\tau ,w,\rho
)=Z_{M^{2}}^{(1)}(q)Z_{M^{2}}^{(1)}(f(\chi ))  \notag \\
=\lim_{\epsilon \rightarrow 0}Z_{M^{2},\epsilon }^{(2)}(\tau _{1},\tau
_{2},\epsilon )=Z_{M^{2}}^{(1)}(q_{1})Z_{M^{2}}^{(1)}(q_{2}),
\label{Zeps_rho_zero}
\end{gather}%
with $Z_{M^{2}}^{(1)}(q)=q^{1/12}/\eta ^{2}(q)$. However, away from a
degeneration point, these partition functions are not equal:

\begin{proposition}
\label{Proposition_Zeps_rho}For two free bosons in the neighborhood of a
two-tori degeneration point we have
\begin{equation*}
\frac{Z_{M^{2},\epsilon }^{(2)}(\tau _{1},\tau _{2},\epsilon )}{%
Z_{M^{2},\rho }^{(2)}(\tau ,w,\rho )}=1+\frac{1}{144}(1-4\chi
)w^{2}+O(w^{4}).
\end{equation*}
\end{proposition}

\noindent \textbf{Proof. }We consider $Z_{M^{2},\epsilon }^{(2)}(\tau
_{1},\tau _{2},\epsilon )$ to $O(\epsilon ^{2})$ and then apply Proposition %
\ref{Proposition_eps_rho_w4} to $O(w^{2})$ to find that%
\begin{eqnarray*}
Z_{M^{2},\epsilon }^{(2)}(\tau _{1},\tau _{2},\epsilon )
&=&Z_{M^{2}}^{(1)}(q_{1})Z_{M^{2}}^{(1)}(q_{2})(1+E_{2}(\tau _{1})E_{2}(\tau
_{2})\epsilon ^{2}+O(\epsilon ^{4})) \\
&=&Z_{M^{2}}^{(1)}\left( q(1+\frac{1}{12}(1-4\chi )w^{2})\right)
Z_{M^{2}}^{(1)}(f(\chi )) \\
&&\left( 1+E_{2}(\tau )E_{2}(q=f(\chi ))(1-4\chi )w^{2}+O(w^{4})\right) .
\end{eqnarray*}%
Recall that for modular form $g_{k}$ of weight $k$ on $SL(2, \mathbb{Z})$,
the `modular derivative'
\begin{equation}
\left(\frac{1}{2\,\pi i\,}\frac{d}{d \tau }+kE_{2}(\tau )\right)g_{k}(\tau)
\label{qdqg}
\end{equation}%
is a modular form of weight $k+2$. Moreover the modular derivative of a
cusp-form is again a cusp-form. In particular, the modular derivative of the
cusp-form $\eta ^{24}(\tau )$ of weight $12$ vanishes since there are no
nonzero cusp-forms of weight $14$. Hence we find (as is well-known) that
\begin{equation}
\frac{1}{2\,\pi i\,}\frac{d}{d \tau }\eta ^{-2}(\tau )=E_{2}(\tau )\eta
^{-2}(\tau ).  \label{qdqeta}
\end{equation}%
Therefore
\begin{equation*}
Z_{M^{2}}^{(1)}(q(1+\frac{1}{12}(1-4\chi )w^{2}))=Z_{M^{2}}^{(1)}(q) \left(
1+\frac{1}{12}(E_{2}(\tau )+\frac{1}{12})(1-4\chi )w^{2}\right) +O(w^{4}).
\end{equation*}%
Altogether, we therefore find
\begin{equation}
\frac{Z_{M^{2},\epsilon }^{(2)}(\tau _{1},\tau _{2},\epsilon )}{%
Z_{M^{2}}^{(1)}(q)Z_{M^{2}}^{(1)}(f(\chi ))}=1+(E_{2}(\tau )G(\chi )+\frac{1%
}{144})(1-4\chi )w^{2}+O(w^{4}),  \label{Z2eps_wsq}
\end{equation}%
with $G(\chi )=E_{2}(q=f(\chi ))-E_{2}(0)=E_{2}(f(\chi ))+\frac{1}{12}$.

We next compute $Z_{M^{2},\rho }^{(2)}(\tau ,w,\rho )$ to $O(w^{2})$ using (%
\ref{Z2_1bos_rho}). We expand $R$ to $O(w^{2})$ as described in Subsection
6.4 of \cite{MT2} and in (\ref{Rexp}) of the Appendix. In particular, (\ref%
{Zeps_rho_zero}) implies
\begin{equation}
\det (I-R^{(0)})^{-1}=Z_{M^{2}}^{(1)}(f(\chi )),  \label{detRzero}
\end{equation}%
with $R^{(0)}$ of (\ref{R0}). Then we find
\begin{eqnarray*}
\det (I-R)^{-1} &=&Z_{M^{2}}^{(1)}(f(\chi )) \det (I+w^{2}\chi E_{2}(\tau
)(I-R^{(0)})^{-1}\Delta )+O(w^{4}) \\
&=&Z_{M^{2}}^{(1)}(f(\chi )) \left( 1+w^{2}\chi E_{2}(\tau )\sigma
((I-R^{(0)})^{-1}(1,1))\right) +O(w^{4}),
\end{eqnarray*}%
with $\Delta _{ab}(k,l)=\delta _{k1}\delta _{l1}$. But applying (\ref{R0G})
it follows that
\begin{equation*}
\det (I-R)^{-1}=Z_{M^{2}}^{(1)}(f(\chi )) \left( 1+E_{2}(\tau )G(\chi
)(1-4\chi )w^{2}\right) +O(w^{4}).
\end{equation*}%
Hence we find
\begin{equation}
\frac{Z_{M^{2},\rho }^{(2)}(\tau ,w,\rho )}{%
Z_{M^{2}}^{(1)}(q)Z_{M^{2}}^{(1)}(f(\chi ))}=1+E_{2}(\tau )G(\chi )(1-4\chi
)w^{2}+O(w^{4}).  \label{Z2rho_wsq}
\end{equation}%
Comparing to (\ref{Z2eps_wsq}) we obtain the result. $\square $

\medskip As described in Subsections 6.2 and 8.2, the bosonic genus two
partition functions in both formalisms have enhanced modular properties on
introducing a standard factor of $(q_{1}q_{2})^{-c/24}$ \ in the $\epsilon $%
-formalism and $q^{-c/24}$ in the $\rho $-formalism. It is therefore natural
to also compare the modular partition functions $Z_{M^{2},\mathrm{mod}%
,\epsilon }^{(2)}(\tau _{1},\tau _{2},\epsilon )$ of (\ref{Z2mod_eps}) and $%
Z_{M^{2},\mathrm{mod},\rho }^{(2)}(\tau ,w,\rho )$ of (\ref{Z2mod_rho}). We
firstly note that $Z_{M^{2},\mathrm{mod},\epsilon }^{(2)}(\tau _{1},\tau
_{2},\epsilon )q_{2}^{1/12}$ and $Z_{M^{2},\mathrm{mod},\rho }^{(2)}(\tau
,w,\rho )$ agree in the two tori degeneration limit. Furthermore, we find
from above that
\begin{equation*}
\frac{Z_{M^{2},\mathrm{mod},\epsilon }^{(2)}(\tau _{1},\tau _{2},\epsilon
)q_{2}^{1/12}}{Z_{M^{2},\mathrm{mod},\rho }^{(2)}(\tau ,w,\rho )}=1+O(w^{4}).
\end{equation*}%
(This result also follows from the fact that the ratio is invariant under
the action of $\Gamma _{1}$ following Theorems \ref{Theorem_Z2_G} and \ref%
{Theorem_Z2_rho_G1}). It follows that
\begin{equation}
\frac{Z_{M^{2},\mathrm{mod},\epsilon }^{(2)}(\tau _{1},\tau _{2},\epsilon )}{%
Z_{M^{2},\mathrm{mod},\rho }^{(2)}(\tau ,w,\rho )}=f(\chi
)^{-1/12}(1+O(w^{4})).  \label{modratio}
\end{equation}%
Thus the modular partition functions do not agree in the two torus
degeneration limit! Eqn. (\ref{modratio}) suggests that a further
modification of the modular partition functions be considered. We consider a
possible further $\Gamma _{1}$-invariant factor of $f(\chi )^{-c/24}$ in the
$\rho $-formalism and find

\begin{proposition}
\label{Proposition_Zeps_rho_mod}For two free bosons in the neighborhood of a
two-tori degeneration point we have
\begin{equation*}
\frac{Z_{M^{2},\mathrm{mod},\epsilon }^{(2)}(\tau _{1},\tau _{2},\epsilon )}{%
Z_{M^{2},\mathrm{mod},\rho }^{(2)}(\tau ,w,\rho )f(\chi )^{-1/12}}%
=1+E_{4}(\tau )(\frac{73}{1440}+39H(\chi ))\left( 1-4\chi \right)
^{2}w^{4}+O(w^{6}),
\end{equation*}%
where $H(\chi )=E_{4}(q=f(\chi ))-E_{4}(q=0)=E_{4}(f(\chi ))-\frac{1}{720}$.
\end{proposition}

\noindent \textbf{Proof.} We give a brief sketch of the proof. First
consider $Z_{M^{2},\mathrm{mod},\epsilon }^{(2)}(\tau _{1},\tau
_{2},\epsilon )$ from (\ref{Z2_1bos}) to $O(\epsilon ^{4})$ to find

\begin{equation}
\frac{1+E_{2}(\tau _{1})E_{2}(\tau _{2})\epsilon ^{2}+\left( E_{2}(\tau
_{1})^{2}E_{2}(\tau _{2})^{2}+54E_{4}(\tau _{1})E_{4}(\tau _{2})\right)
\epsilon ^{4}}{\eta (\tau _{1})^{2}\eta (\tau _{2})^{2}}+O(\epsilon ^{6}).
\label{Zeps4}
\end{equation}%
We expand to $O(w^{4})$ using Proposition \ref{Proposition_eps_rho_w4} and
use (\ref{qdqeta}) and (\ref{qdqE2}) to eventually find that%
\begin{gather*}
Z_{M^{2},\mathrm{mod},\epsilon }^{(2)}(\tau _{1},\tau _{2},\epsilon )=\frac{1%
}{\eta (q)^{2}\eta (f(\chi ))^{2}}[1+E_{2}\left( \tau \right) G(\chi )\left(
1-4\chi \right) w^{2} \\
+\left( E_{4}\left( \tau \right) \left( \frac{73}{1440}+54H(\chi )+\frac{1}{2%
}\,G(\chi )\right) +E_{2}\left( \tau \right) ^{2}G(\chi )^{2}\right) \left(
1-4\,\chi \right) ^{2}w^{4}]+O(w^{6}).
\end{gather*}%
On the other hand, using the expansion of $R$ in (\ref{Rexp}) of the
Appendix and the methods of Subsection 5.2.2 of ref. \cite{MT2}, one
eventually finds that
\begin{gather*}
Z_{M^{2},\mathrm{mod},\rho }^{(2)}(\tau ,w,\rho )f(\chi )^{-1/12}=\frac{1}{%
\eta (q)^{2}\eta (f(\chi ))^{2}}[1+E_{2}\left( \tau \right) G(\chi )\left(
1-4\chi \right) w^{2} \\
+\left( E_{4}\left( \tau \right) \left( 15H(\chi )+\frac{1}{2}\,G(\chi
)\right) +E_{2}\left( \tau \right) ^{2}G(\chi )^{2}\right) \left( 1-4\,\chi
\right) ^{2}w^{4}]+O(w^{6}),
\end{gather*}%
leading to the result. $\square $

\begin{remark}
\label{Remark}Notice that $Z_{M^{2},\mathrm{mod},\epsilon }^{(2)}$ and $%
Z_{M^{2},\mathrm{mod},\rho }^{(2)}f(\chi )^{-1/12}$ differ to $O(w^{4})$ in
the leading $\chi ^{0}$ term. In particular, $R(k,l)=O(\chi )$ so that%
\begin{equation*}
Z_{M^{2},\mathrm{mod},\rho }^{(2)}(\tau ,w,\rho )f(\chi )^{-1/12}=\frac{1}{%
\eta (q)^{2}\eta (f(\chi ))^{2}}+O(\chi ),
\end{equation*}%
whereas (\ref{Zeps4}) gives
\begin{equation*}
Z_{M^{2},\mathrm{mod},\epsilon }^{(2)}(\tau _{1},\tau _{2},\epsilon )=\frac{%
1+\frac{73}{1440}E_{4}\left( \tau \right) w^{4}+O(\chi )}{\eta (q)^{2}\eta
(f(\chi ))^{2}}+O(w^{6}).
\end{equation*}
\end{remark}

\section{Final Remarks}

We have defined and explicitly calculated the genus two partition function
for the Heisenberg free bosonic string and lattice VOAs in two separate
sewing schemes, the $\epsilon$-  and $\rho $-formalisms. It is important to
emphasize that our approach is constructively based only on the algebraic
properties of a VOA. In particular, we have proved the holomorphy of the
genus two partition functions on given domains $\mathcal{D}^{\epsilon }$ and
$\mathcal{D}^{\rho }$ in the $\epsilon $- and $\rho $-formalisms
respectively, and have described modular properties for the modular
partition functions. Lastly, we have shown that although the original (and
modular) partition functions enjoy many properties in common between the two
sewing schemes, for a given VOA they do not agree in the neighborhood of a
two-tori degeneration point.

It is natural to try generalize these results to Riemann surfaces of higher
genus $g$. We can inductively construct such surfaces from tori by sewing
together lower genus surfaces and/or attaching handles \cite{Y}, \cite{MT2}.
In particular, we conjecture that for a general $\epsilon $-sewing scheme
joining together surfaces of genus $g_{1}$ and $g_{2}$, then Theorem \ref%
{Theorem_Z2_boson} generalizes to
\begin{equation*}
Z_{M,\epsilon }^{(g_{1}+g_{2})}=Z_{M}^{(g_{1})}Z_{M}^{(g_{2})}(\det
(1-A_{1}A_{2}))^{-1/2},
\end{equation*}%
where now $A_{1}$ and $A_{2}$ are the moment matrices of Section 4. in ref.
\cite{MT2} and $Z_{M}^{(g_{a})}$ is the genus $g_{a}$ partition function (in
some chosen sewing scheme). Similarly, we conjecture that for a lattice VOA
then Theorem \ref{Theorem_Z2_L_eps} generalizes to
\begin{equation*}
Z_{L,\epsilon }^{(g_{1}+g_{2})}=Z_{M^{l},\epsilon }^{(g_{1}+g_{2})}\theta
_{L}^{(g_{1}+g_{2})}(\Omega ),
\end{equation*}%
where $\theta _{L}^{(g)}(\Omega )$ is the genus $g$ Siegel theta function
associated to $L$ and $\Omega $ is the genus $g$ period matrix. Similarly,
we can consider attaching a handle to a genus $g$ surface in a general $\rho
$-sewing scheme to conjecture the following generalizations of Theorems \ref%
{Theorem_Z2_boson_rho} and \ref{Theorem_Z2_L_rho}:
\begin{eqnarray*}
Z_{M,\rho }^{(g+1)} &=&Z_{M}^{(g)}(\det (1-R))^{-1/2}, \\
Z_{L,\rho }^{(g+1)} &=&Z_{M^{l},\rho }^{(g+1)}\theta _{L}^{(g+1)}(\Omega ),
\end{eqnarray*}%
where $R$ is the moment matrix of Section 5 in \cite{MT2}. Indeed,
these can be confirmed for $g=0$ for the two sewing schemes described in
Section 5.2. Then one obtains, respectively, $Z_{M}^{(1)}(q)=(\det
(1-R))^{-1/2}=q^{1/24}/\eta (q)$ for $R=\mathrm{diag}%
(1-q,1-q,1-q^{2},1-q^{2},\ldots )$ and $Z_{M,\rho }^{(1)}=(\det
(1-R^{(0)}))^{-1/2}=Z_{M}^{(1)}(f(\chi ))$ from (\ref{detRzero}).

\medskip
In conclusion, let us briefly and heuristically sketch how these results
compare to some related ideas in the physics and mathematics literature.
There is a wealth of literature concerning the bosonic string e.g. \cite{GSW}%
, \cite{P}. In particular, the conformal anomaly implies that the physically
defined path integral partition function $Z_{\mathrm{string}}$ cannot be
reduced to an integral over the moduli space $\mathcal{M}_{g}$ of a Riemann
surface of genus $g$ except for the 26 dimensional critical string where the
anomaly vanishes. Furthermore, for the critical string, Belavin and Knizhnik
argue that%
\begin{equation*}
Z_{\mathrm{string}}=\int\limits_{\mathcal{M}_{g}}\left\vert F\right\vert
^{2}d\mu ,
\end{equation*}%
where $d\mu $ denotes a natural volume form on $\mathcal{M}_{g}$ and $F$ is
holomorphic and non-vanishing on $\mathcal{M}_{g}$ \cite{BK}, \cite{Kn}.
They also claim that for $g\geq 2$, $F$ is a global section for the line
bundle $K\otimes \lambda ^{-13}$ (where $K$ is the canonical bundle and $%
\lambda $ the Hodge bundle on $\mathcal{M}_{g}$) which is trivial by
Mumford's theorem in algebraic geometry \cite{Mu2}. In this identification,
the $\lambda ^{-13}$ section  is associated with 26 bosons, the $K$ section
with a $c=-26$ ghost system and the vanishing conformal anomaly to the
vanishing first Chern class  for $K\otimes \lambda ^{-13}$ \cite{N}.
Recently, some of these ideas have also been rigorously proved for a zeta
function regularized determinant of an appropriate Laplacian operator $%
\Delta _{n}$ \cite{McT}. The genus two partition functions $%
Z_{M^{2},\epsilon }^{(2)}(\tau _{1},\tau _{2},\epsilon )$\ and $%
Z_{M^{2},\rho }^{(2)}(\tau ,w,\rho )$ constructed in this paper for a rank $2
$ Heisenberg VOA should correspond in these approaches to a local
description of the holomorphic part of $\left( \frac{\det^{\prime }\Delta
_{1}}{\det N_{1}}\right) ^{-1}$ of refs. \cite{Kn}, $\ $\cite{McT}, giving a
local section of the line bundle $\lambda ^{-1}$. As such, it is therefore
no surprise that $f_{\epsilon ,\rho }=Z_{M^{2},\epsilon
}^{(2)}/Z_{M^{2},\rho }^{(2)}\neq 1$ in the neighborhood of a two-tori
degeneration point where the ratio of the two sections is a non-trivial
transition function $f_{\epsilon ,\rho }$.

In the case of a general rational conformal field theory, the conformal
anomaly continues to obstruct the existence of a global partition function
on moduli space for $g\geq 2$ \footnote{%
This observation appears to be contrary to the main result of \cite{So2}.}.
However, \emph{all} CFTs of a given central charge $c$ are believed to share
the same conformal anomaly e.g. \cite{FS}. Thus, the identification of the
normalized lattice partition functions of (\ref{normpfs}) presumably
reflects the equality of the first Chern class of some bundle associated to
a rank $c$ lattice VOA to that for $\lambda ^{-c}$ with transition function $%
f_{\epsilon ,\rho }^{c/2}$. It is interesting to note that even in the case
of a unimodular lattice VOA with a unique conformal block (\cite{MS}, \cite%
{TUY}) the genus two partition function can therefore only be described
locally. It would obviously be extremely valuable to find a rigorous
description of the relationship between the VOA approach described here and
these related ideas in conformal field theory and algebraic geometry.

\section{Appendix}

\subsection{A product formula}

Here we continue the discussion initiated in Subsection 3.1, with a view to
proving Proposition \ref{Prop_Om12_R21expansion}. Let $\mathcal{M}(I)$ be
the (multiplicative semigroup generated by) the rotationless cycles in the
symbols $x_{i},i\in I$. There is an injection
\begin{equation}
\iota :\bigcup_{n\geq 0}C_{n}\longrightarrow \mathcal{M}(I)
\label{eq: cycleinjection}
\end{equation}%
defined as follows. If $(x)\in C_{n}$ has rotation group of order $r$ then $%
r|n$ and there is a rotationless monomial $y$ such that $x=y^{r}$. We then
map $(x)\mapsto (y)^{r}$. It is readily verified that this is well-defined.
In this way, each cycle is mapped to a power of a rotationless cycle in $%
\mathcal{M}(I)$. A typical element of $\mathcal{M}(I)$ is uniquely
expressible in the form
\begin{equation}
p_{1}^{f_{1}}p_{2}^{f_{2}}...p_{k}^{f_{k}}  \label{eq:RF0}
\end{equation}%
where $p_{1},...,p_{k}$ are distinct rotationless cycles and $%
f_{1},...,f_{k} $ are non-negative integers. We call (\ref{eq:RF0}) the
\textit{reduced form} of an element in $\mathcal{M}(I)$. A general element
of $\mathcal{M}(I)$ is then essentially a labelled graph, each of whose
connected components are rotationless labelled polygons as discussed in
Subsection 3.1.

Now consider a second finite set $T$ together with a map
\begin{equation}
F:T\longrightarrow I.  \label{eq:RF1}
\end{equation}%
Thus elements of $I$ label elements of $T$ via the map $F$. $F$ induces a
natural map
\begin{equation*}
\overline{F}:\Sigma (T)\longrightarrow \mathcal{M}(I)
\end{equation*}%
from the symmetric group $\Sigma (T)$ as follows. For an element $\tau \in
\Sigma (T)$, write $\tau $ as a product of disjoint cycles $\tau =\sigma
_{1}.\sigma _{2}\ldots $. We set $\overline{F}(\tau )=\overline{F}(\sigma
_{1})\overline{F}(\sigma _{2})...$, so it suffices to define $\overline{F}%
(\sigma )$ for a cycle $\sigma =(s_{1}s_{2}...)$ with $s_{1},s_{2},...\in T$%
. In this case we set
\begin{equation*}
\overline{F}(\sigma )=\iota ((x_{F(s_{1})}x_{F(s_{2})}...))
\end{equation*}%
where $\iota $ is as in (\ref{eq: cycleinjection}). When written in the form
(\ref{eq:RF0}), we call $\overline{F}(\tau )$ the \textit{reduced $F$-form}
of $\tau $.

\vspace{0.15in} \noindent For $i\in {I}$, let $s_{i}=|F^{-1}(i)|$ be the
number of elements in $T$ with label $i$. So the number of elements in $T$
is equal to ${\sum }_{i\in I}s_{i}$. We say that two elements $\tau
_{1},\tau _{2}\in \Sigma (T)$ are $F$-equivalent if they have the same
reduced $F$-form, i.e. $\overline{F}(\tau _{1})=\overline{F}(\tau _{2})$. We
will show that each equivalence class contains the same number of elements.
Precisely,

\begin{lemma}
\label{Lemma_Fequiv} Each $F$-equivalence class contains precisely ${\prod }%
_{i\in I}s_{i}$ elements. In particular, the number of $F$-equivalence
classes is $|T|!/{\prod }_{i\in I}s_{i}$.
\end{lemma}

\noindent \textbf{Proof.} An element $\tau \in \Sigma (T)$ may be
represented uniquely as
\begin{equation*}
\left(
\begin{array}{cccc}
0 & 1 & \cdots & M \\
\tau (0) & \tau (1) & \cdots & \tau (M)%
\end{array}%
\right)
\end{equation*}%
so that
\begin{equation*}
\overline{F}(\tau )=\left(
\begin{array}{cccc}
F(0) & F(1) & \cdots & F(M) \\
F(\tau (0)) & F(\tau (1)) & \cdots & F(\tau (M))%
\end{array}%
\right)
\end{equation*}%
with an obvious notation. Exactly $s_{i}$ of the $\tau (j)$ satisfy
\begin{equation*}
\overline{F}(\tau (j))=x_{i}
\end{equation*}%
so that there are $\prod_{i\in I}s_{i}$ choices of $\tau $ which have a
given image under $\overline{F}$. The Lemma follows. ${\square }$

\medskip The next results employs notation introduced in Subsections 3.1 and
3.2.

\begin{lemma}
\label{Lemma_L21_expansion} We have
\begin{equation}
(I-M_{1}M_{2})^{-1}(1,1)=(1-\sum_{L\in \mathcal{L}_{21}}\omega (L))^{-1}.
\label{eq: matrixsum}
\end{equation}
\end{lemma}

\noindent As before, the left-hand-side of (\ref{eq: matrixsum}) means $%
\sum_{n\geq 0}(M_{1}M_{2})^{n}(1,1)$. It is a certain power series with
entries being quasi-modular forms..

\medskip \noindent \textbf{Proof of Lemma.} We have
\begin{equation}
(M_{1}M_{2})^{n}(1,1)=\sum M_{1}(1,k_{1})M_{2}(k_{1},k_{2})\ldots
M_{2}(k_{2n-1},1)  \label{eq: matrixentry}
\end{equation}%
where the sum ranges over all choices of positive integers $%
k_{1},...,k_{2n-1}$. Such a choice corresponds to a (isomorphism class of)
chequered cycle $L$ with $2n$ nodes and with at least one distinguished
node, so that the left-hand-side of (\ref{eq: matrixsum}) is equal to
\begin{equation*}
\sum_{L}\omega (L)
\end{equation*}%
summed over all such $L$. We can formally write $L$ as a product $%
L=L_{1}L_{2}...L_{p}$ where each $L_{i}\in \mathcal{L}_{21}$. This indicates
that $L$ has $p$ distinguished nodes and that the $L_{i}$ are the edges of $%
L $ between consecutive distinguished nodes, which can be naturally thought
of as chequered cycles in $\mathcal{L}_{21}$. Note that in the
representation of $L$ as such a product, the $L_{i}$ do not commute unless
they are equal, moreover $\omega $ is multiplicative. Then
\begin{equation*}
(I-M_{1}M_{2})^{-1}(1,1)=\sum_{L_{i}\in \mathcal{L}_{21}}\omega
(L_{1}...L_{p})=(1-\sum_{L\in \mathcal{L}_{21}}\omega (L))^{-1}
\end{equation*}%
as required. \ \ \ \ \ $\square $

\begin{proposition}
\label{Prop_R21expansion}We have
\begin{equation}
(I-M_{1}M_{2})^{-1}(1,1)=\prod_{L\in \mathcal{R}_{21}}(1-\omega (L))^{-1}
\label{eq: matrixprod}
\end{equation}
\end{proposition}

\noindent \textbf{Proof. }By Lemma \ref{Lemma_L21_expansion} we have
\begin{equation}
(I-M_{1}M_{2})^{-1}(1,1)=\sum m(e_{1},...,e_{k})\omega
(L_{1})^{e_{1}}...\omega (L_{k})^{e_{k}}  \label{eq: mult}
\end{equation}%
where the sum ranges over distinct elements $L_{1},...L_{k}$ of $\mathcal{L}%
_{21}$ and all $k$-tuples of non-negative integers $e_{1},...,e_{k}$, and
where the multiplicity is
\begin{equation*}
m(e_{1},...,e_{k})=\frac{(\sum_{i}e_{i})!}{\prod_{i}(e_{i}!)}.
\end{equation*}%
Let $S$ be the set consisting of $e_{i}$ copies of $L_{i},1\leq i\leq k$,
let $I$ be the integers between $1$ and $k$, and let $F:S\longrightarrow I$
be the obvious labelling map. A reduced $F$-form is then an element of $%
\mathcal{M}(I)$ where the variables $x_{i}$ are now the $L_{i}$. The free
generators of $\mathcal{M}(I)$, i.e. rotationless cycles in the $x_{i}$, are
naturally identified \emph{precisely} with the elements of $\mathcal{R}_{21}$%
, and Lemma \ref{Lemma_Fequiv} implies that each element of $\mathcal{M}(I)$
corresponds to just one term under the summation in (\ref{eq: mult}). Eqn.(%
\ref{eq: matrixprod}) follows immediately from this and the multiplicativity
of $\omega $, and the Proposition is proved. \ \ \ \ $\square $

\subsection{Invertibility about a two tori degeneration point}

Here we prove Proposition \ref{Proposition_eps_rho_w4} describing the
relationship to $O(w^{4})$ between $(\tau _{1},\tau _{2},\epsilon )\in
\mathcal{D}^{\epsilon }$ and $(\tau ,w,\chi )\in \mathcal{D}^{\chi }\cup
\mathcal{D}_{0}^{\chi }$ in the neighborhood of a two tori degeneration
point $\rho ,w\rightarrow 0$ for fixed $\chi =-\rho /w^{2}$. We first expand
$\Omega _{ij}$ to $O(\epsilon ^{4})$ using Theorem \ref{Theorem_period_eps}
to obtain%
\begin{eqnarray*}
2\,\pi i\Omega _{11} &=&2\,\pi i\,\tau _{{1}}+E_{2}\left( \tau _{2}\right) {%
\epsilon }^{2}+E_{2}\left( \tau _{1}\right) E_{2}\left( \tau _{2}\right) ^{2}%
{\epsilon }^{4}+O(\epsilon ^{6}), \\
2\,\pi i\Omega _{12} &=&-{\epsilon }-E_{2}\left( \tau _{1}\right)
E_{2}\left( \tau _{2}\right) {\epsilon }^{3}+O(\epsilon ^{5}), \\
2\,\pi i\Omega _{22} &=&2\,\pi i\,\tau _{{2}}+E_{2}\left( \tau _{1}\right) {%
\epsilon }^{2}+E_{2}\left( \tau _{1}\right) ^{2}E_{2}\left( \tau _{2}\right)
{\epsilon }^{4}+O(\epsilon ^{6}).
\end{eqnarray*}%
Making use of the identity (cf. (\ref{qdqg}))
\begin{equation}
\frac{1}{2\,\pi i\,}\frac{d }{d \tau }E_{2}(\tau )=5E_{4}(\tau )-E_{2}(\tau
)^{2},  \label{qdqE2}
\end{equation}%
it is straightforward to invert $\Omega _{ij}(\tau _{1},\tau _{2},\epsilon )$
to find

\begin{lemma}
\label{Lemma_Omeps}In the neighborhood of a two tori degeneration point $%
\Omega _{12}=0$ of $\Omega \in \mathbb{H}_{2}$ we have
\begin{eqnarray*}
2\,\pi i\,\tau _{{1}} &=&2\,\pi i\Omega _{11}-E_{2}\left( \Omega
_{22}\right) {r}^{2}+5E_{2}\left( \Omega _{11}\right) E_{4}\left( \Omega
_{22}\right) {r}^{4}+O(r^{6}), \\
{\epsilon } &=&-r+E_{2}\left( \Omega _{11}\right) E_{2}\left( \Omega
_{22}\right) {r}^{3}+O(r^{5}), \\
2\,\pi i\,\tau _{2} &=&2\,\pi i\Omega _{22}-E_{2}\left( \Omega _{11}\right) {%
r}^{2}+5E_{2}\left( \Omega _{22}\right) E_{4}\left( \Omega _{11}\right) {r}%
^{4}+O(r^{6}),
\end{eqnarray*}%
where $r=2\,\pi i\Omega _{12}$. $\Box$
\end{lemma}

We next obtain $\Omega _{ij}(\tau ,w,\chi )$ to $O(w^{4})$ in the
neighborhood of a two tori degeneration point:

\begin{proposition}
\label{PropOmrhodegen}For $(\tau ,w,\chi )\in \mathcal{D}^{\chi }\cup
\mathcal{D}_{0}^{\chi }$ we have
\begin{eqnarray}
2\pi i\Omega _{11} &=&2\pi i\tau +(1-4\chi )G(\chi )w^{2}+\left( 1-4\chi
\right) ^{2}G(\chi )^{2}E_{2}\left( \tau \right) w^{4}+O(w^{6}),  \notag \\
&&  \label{Om11_deg} \\
2\pi i\Omega _{12} &=&w\sqrt{1-4\chi }(1+(1-4\chi )G(\chi )E_{2}(\tau
)w^{2})+O(w^{5}),  \label{Om12_deg} \\
2\pi i\Omega _{22} &=&\log f(\chi )+(1-4\chi )E_{2}(\tau )w^{2}  \notag \\
&&+\left( 1-4\,\chi \right) ^{2}\left( G(\chi )E_{2}\left( \tau \right) ^{2}+%
\frac{1}{2}E_{4}\left( \tau \right) \right) w^{4}+O(w^{6}),  \label{Om22_deg}
\end{eqnarray}%
where
\begin{equation}
G(\chi )=\frac{1}{12}+E_{2}(q=f(\chi )),  \label{Gchi}
\end{equation}%
and $f(\chi )$ is the Catalan series (\ref{Catalan}).
\end{proposition}

\noindent \textbf{Proof. }$\Omega _{12}$ is described in Proposition 13 of
\cite{MT2}, whereas $\Omega _{11}$ and $\Omega _{22}$ are described there
only to order $O(w^{2})$. The $O(w^{4})$ terms are calculated similarly, as
follows. From (\ref{P2}) and (\ref{Pkdef}) we have
\begin{eqnarray*}
P_{1}(\tau ,w) &=&\frac{1}{w}(1-E_{2}(\tau )w^{2}-E_{4}(\tau
)w^{4}+O(w^{6})), \\
P_{2}(\tau ,w) &=&\frac{1}{w^{2}}(1+E_{2}(\tau )w^{2}+3E_{4}(\tau
)w^{4}+O(w^{6})), \\
P_{3}(\tau ,w) &=&\frac{1}{w^{3}}(1-3E_{4}(\tau )w^{4}+O(w^{6})), \\
P_{4}(\tau ,w) &=&\frac{1}{w^{4}}(1+E_{4}(\tau )w^{4}+O(w^{6})),
\end{eqnarray*}%
and $P_{n}(\tau ,w)=w^{-n}(1+O(w^{6}))$ for $n>4$. Then $R(k,l)$ of (\ref%
{Rdef}) gives
\begin{equation}
R(k,l)=R^{(0)}(k,l)+R^{(2)}(k,l)w^{2}+R^{(4)}(k,l)w^{4}+O(w^{6}),
\label{Rexp}
\end{equation}%
where
\begin{equation}
R^{(0)}(k,l)=\frac{(-\chi )^{(k+l)/2}}{\sqrt{kl}}\frac{(k+l-1)!}{(k-1)!(l-1)!%
}\left[
\begin{array}{ll}
(-1)^{k} & 0 \\
0 & (-1)^{l}%
\end{array}%
\right] .  \label{R0}
\end{equation}%
$R^{(0)}$ is associated with the self-sewing of a sphere to form a torus
described in Subsections 5.2.2 and 6.4 of \cite{MT2}. $R^{(2)}(k,l)$ is
given by%
\begin{equation}
R^{(2)}(k,l)=\chi E_{2}(\tau )\Delta (k,l),\quad \Delta (k,l)=\delta
_{k1}\delta _{l1}\left[
\begin{array}{ll}
1 & 1 \\
1 & 1%
\end{array}%
\right] ,  \label{R2}
\end{equation}%
whereas $R^{(4)}(k,l)=0$ for $k+l>4$ with non-zero entries
\begin{eqnarray*}
R^{(4)}(1,1) &=&3\chi E_{4}(\tau )\left[
\begin{array}{ll}
1 & 1 \\
1 & 1%
\end{array}%
\right] , \\
R^{(4)}(1,2) &=&-R^{(4)}(2,1)=3\sqrt{2}(-\chi )^{3/2}E_{4}(\tau )\left[
\begin{array}{ll}
1 & 0 \\
0 & -1%
\end{array}%
\right] , \\
R^{(4)}(1,3) &=&R^{(4)}(3,1)=\sqrt{3}\chi ^{2}E_{4}(\tau )\left[
\begin{array}{ll}
1 & 1 \\
1 & 1%
\end{array}%
\right] , \\
R^{(4)}(2,2) &=&3\chi ^{2}E_{4}(\tau )\left[
\begin{array}{ll}
1 & 1 \\
1 & 1%
\end{array}%
\right] .
\end{eqnarray*}%
We find for $b$ of (\ref{bk}) that
\begin{equation*}
b(k)=b^{(0)}(k)+b^{(2)}(k)w^{2}+b^{(4)}(k)w^{4}+O(w^{6}),\
\end{equation*}%
with%
\begin{eqnarray*}
b^{(0)}(k) &=&\frac{(-\chi )^{k/2}}{\sqrt{k}}[-1,(-1)^{k}], \\
b^{(2)}(k) &=&-\delta _{k1}E_{2}(q)b^{(0)}(1),
\end{eqnarray*}%
and $b^{(4)}(k)=0$ for $k>3$ with non-zero entries
\begin{eqnarray*}
b^{(4)}(1) &=&-E_{4}(\tau )b^{(0)}(1), \\
b^{(4)}(2) &=&3E_{4}(\tau )b^{(0)}(2), \\
b^{(4)}(3) &=&-3E_{4}(\tau )b^{(0)}(3).
\end{eqnarray*}%
It is convenient to define
\begin{equation*}
T=(I-R)^{-1}=\sum_{n\geq 0}R^{n}=T^{(0)}+T^{(2)}w^{2}+T^{(4)}w^{4}+O(w^{6}),
\end{equation*}%
with
\begin{eqnarray*}
T^{(0)} &=&(I-R^{(0)})^{-1}, \\
T^{(2)} &=&T^{(0)}R^{(2)}T^{(0)}, \\
T^{(4)} &=&T^{(0)}R^{(2)}T^{(0)}R^{(2)}T^{(0)}+T^{(0)}R^{(4)}T^{(0)}.
\end{eqnarray*}

We next recall from Proposition 10 and (170) of op.cite. that
\begin{equation}
\chi \sigma (T^{(0)}(1,1))=\chi \sigma ((I-R^{(0)})^{-1}(1,1))=(1-4\chi
)G(\chi ).  \label{R0G}
\end{equation}%
We also find that
\begin{equation*}
(b^{(0)}T^{(0)})_{a}(k)=(T^{(0)}\bar{b}^{(0)T})_{a}(k)=(-1)^{k}(-1)^{(k+1)a}.%
\frac{(-\chi )^{k/2}}{\sqrt{k}}\sum\limits_{n\geq 1}S_{n,k}(\chi ),
\end{equation*}%
where $S_{1,k}(\chi )=1$ and
\begin{eqnarray}
S_{n,k}(\chi ) &=&\sum_{k_{n-1},\ldots k_{1}\geq 1}\chi ^{k_{n-1}+\ldots
+k_{1}}\binom{k+k_{n-1}-1}{k_{n-1}}\binom{k_{n-1}+k_{n-2}-1}{k_{n-2}}  \notag
\\
&&\ldots \binom{k_{2}+k_{1}-1}{k_{1}},  \label{Snk}
\end{eqnarray}%
for $n>1$. But (129) of op.cite. states that
\begin{equation*}
\sum\limits_{n\geq 1}S_{n,k}(\chi )=(1+f(\chi ))^{k}.
\end{equation*}%
Hence we find%
\begin{equation}
\frac{(-\chi )^{k/2}}{\sqrt{k}}(b^{(0)}T^{(0)})_{a}(k)=\frac{(-\chi )^{k/2}}{%
\sqrt{k}}(T^{(0)}\bar{b}^{(0)T})_{a}(k)=(-1)^{(k+1)a}.\frac{X(\chi )^{k}}{k},
\label{bR0}
\end{equation}%
where for later convenience we have defined
\begin{equation*}
X=X(\chi )=\chi (1+f(\chi )).
\end{equation*}%
We also note that
\begin{equation}
(1-2X)^{2}=1-4\chi .  \label{Xchi}
\end{equation}

We may then compute $\Omega _{ij}$ to $O(w^{4})$ employing the identities (%
\ref{R0G}), (\ref{bR0}) and (\ref{Xchi}) as follows. We begin with $\Omega
_{11}$ of (\ref{Om11rho}) which to $O(w^{4})$ is given by
\begin{equation*}
2\pi i\Omega _{11}=2\pi i\tau +w^{2}\chi \sigma (T^{(0)}(1,1))+w^{4}\chi
\sigma (T^{(2)}(1,1))+O(w^{6}).
\end{equation*}%
Using (\ref{R2}) we find
\begin{eqnarray*}
\chi \sigma (T^{(2)}(1,1)) &=&E_{2}(\tau )\chi ^{2}\sigma (T^{(0)}\Delta
T^{(0)})(1,1) \\
&=&E_{2}(\tau )(\chi \sigma (T^{(0)}(1,1)))^{2} \\
&=&E_{2}(\tau )(1-4\chi )^{2}G(\chi )^{2},
\end{eqnarray*}%
on applying (\ref{R0G}). Hence (\ref{Om11_deg}) follows.

Expanding $\Omega _{22}$ of (\ref{Om22rho}) to $O(w^{4})$ we find
\begin{eqnarray*}
2\pi i\Omega _{22} &=&\log \chi -b^{(0)}T^{(0)}{}\bar{b}^{(0)T}+\left(
E_{2}(\tau )-b^{(0)}T^{(2)}{}\bar{b}^{(0)T}-2b^{(2)}T^{(0)}{}\bar{b}%
^{(0)T}\right) w^{2} \\
&&+(\frac{1}{2}E_{4}(\tau )-b^{(0)}T^{(4)}{}\bar{b}^{(0)T}-b^{(2)}T^{(0)}{}%
\bar{b}^{(2)T}-2b^{(2)}T^{(2)}{}\bar{b}^{(0)T} \\
&&-2b^{(4)}T^{(0)}{}\bar{b}^{(0)T})w^{4}+O(w^{6}).
\end{eqnarray*}%
Proposition 9 of op.cite. states that
\begin{equation*}
\log f(\chi )=\log \chi -b^{(0)}T^{(0)}{}\bar{b}^{(0)T}.
\end{equation*}%
By applying (\ref{R0G}) and (\ref{bR0}) we obtain%
\begin{eqnarray*}
-b^{(0)}T^{(2)}{}\bar{b}^{(0)T} &=&E_{2}(\tau )\left[ (-\chi )^{1/2}\sigma
(b^{(0)}T^{(0)})(1)\right] ^{2}=4E_{2}(\tau )X^{2}, \\
-2b^{(2)}T^{(2)}{}\bar{b}^{(0)T} &=&-2E_{2}(\tau )(-\chi )^{1/2}\sigma
(b^{(0)}T^{(0)})(1)=-4E_{2}(\tau )X.
\end{eqnarray*}%
Thus the $O(w^{2})$ term\ of $2\pi i\Omega _{22}$ is found using (\ref{Xchi}%
).

Similarly, one finds that%
\begin{eqnarray*}
-b^{(0)}T^{(4)}{}\bar{b}^{(0)T} &=&E_{2}(\tau )^{2} \chi \sigma
(T^{(0)})(1,1) \left[ (-\chi )^{1/2}\sigma (b^{(0)}T^{(0)})(1)\right] ^{2} \\
&&+3E_{4}(\tau ) \sum_{a=1,2}\left[ (-\chi )^{1/2}(b^{(0)}T^{(0)})_{a}(1)%
\right] ^{2} \\
&&-6E_{4}(\tau ) (-\chi )^{1/2}\sigma (b^{(0)}T^{(0)})(1).\frac{(-\chi )^{1}%
}{\sqrt{2}}\sigma (b^{(0)}T^{(0)})(2) \\
&&+6E_{4}(\tau ) (-\chi )^{1/2}\sigma (b^{(0)}T^{(0)})(1) \frac{(-\chi
)^{3/2}}{\sqrt{3}}\sigma (b^{(0)}T^{(0)})(3) \\
&=&4E_{2}(\tau )^{2}(1-4\chi )G(\chi )X^{2}+E_{4}(\tau
)(6X^{2}-12X^{3}+8X^{4})
\end{eqnarray*}%
\begin{eqnarray*}
-b^{(2)}T^{(0)}{}\bar{b}^{(2)T} &=&E_{2}(\tau )^{2} \chi \sigma
(T^{(0)})(1,1)=E_{2}(\tau )^{2}(1-4\chi )G(\chi ), \\
-2b^{(2)}T^{(2)}{}\bar{b}^{(0)T} &=&-2E_{2}(q)^{2} \chi \sigma
(T^{(0)})(1,1) (-\chi )^{1/2}\sigma (b^{(0)}T^{(0)})(1) \\
&=&-4E_{2}(q)^{2}(1-4\chi )G(\chi )X, \\
-2b^{(4)}T^{(0)}{}\bar{b}^{(0)T} &=&-2E_{4}(\tau ) [(-\chi )^{1/2}\sigma
(b^{(0)}T^{(0)})(1)-3\frac{(-\chi )^{1}}{\sqrt{2}}\sigma (b^{(0)}T^{(0)})(2)
\\
&&+3\frac{(-\chi )^{3/2}}{\sqrt{3}}\sigma (b^{(0)}T^{(0)})(3)] \\
&=&E_{4}(\tau )(-4X+6X^{2}-4X^{3}).
\end{eqnarray*}%
Combining, we find the coefficient of $E_{2}\left( \tau \right) ^{2}(1-4\chi
)G(\chi )w^{4}$ in $\Omega _{22}$ is $(1-2X)^{2}=1-4\chi $ whereas the
coefficient of $E_{4}\left( \tau \right) w^{4}$ is
\begin{equation*}
\frac{1}{2}+6X^{2}-12X^{3}+8X^{4}-4X+6X^{2}-4X^{3}=\frac{1}{2}(1-2X)^{4}=%
\frac{1}{2}(1-4\chi )^{2}.
\end{equation*}%
Hence (\ref{Om22_deg}) follows. $\square $

\medskip We next combine Lemma \ref{Lemma_Omeps}\ and Proposition \ref%
{PropOmrhodegen} to prove Proposition \ref{Proposition_eps_rho_w4}. Thus $%
2\,\pi i\,\tau _{{1}}$ is given by
\begin{eqnarray*}
&&2\,\pi i\Omega _{11}-E_{2}\left( \Omega _{22}\right) {(2\pi i\Omega _{12})}%
^{2}+5E_{2}\left( \Omega _{11}\right) E_{4}\left( \Omega _{22}\right) {(2\pi
i\Omega _{12})}^{4}+O({\Omega _{12}^{6})} \\
&=&2\pi i\tau +(1-4\chi )G(\chi )w^{2}+\left( 1-4\chi \right) ^{2}G(\chi
)^{2}E_{2}\left( \tau \right) w^{4} \\
&&-\left[ E_{2}\left( f\right) +(5E_{4}(f)-E_{2}(f)^{2})(1-4\chi )E_{2}(\tau
)w^{2}\right] \\
&&.w^{2}(1-4\chi )(1+2(1-4\chi )G(\chi )E_{2}(\tau )w^{2})+ \\
&&+5E_{2}\left( \tau \right) E_{4}\left( f\right) (1-4\chi )^{2}{w}%
^{4}+O(w^{6}) \\
&=&2\pi i\tau +\left[ G(\chi )-E_{2}\left( f\right) \right] (1-4\chi )w^{2},
\end{eqnarray*}%
where $E_{2}\left( f\right) =E_{2}\left( q=f(\chi )\right) $ and $%
E_{4}\left( f\right) =E_{4}\left( q=f(\chi )\right) $. Applying (\ref{Gchi})
we find%
\begin{equation*}
2\,\pi i\,\tau _{{1}}=2\pi i\tau +\frac{1}{12}(1-4\chi )w^{2}+\frac{1}{144}%
E_{2}(\tau )(1-4\chi )^{2}w^{4}+O(w^{6}).
\end{equation*}%
Similarly
\begin{eqnarray*}
{\epsilon } &=&-{2\pi i\Omega _{12}}+E_{2}\left( \Omega _{11}\right)
E_{2}\left( \Omega _{22}\right) {(2\pi i\Omega _{12})}^{3}+O({\Omega
_{12}^{5})} \\
&=&-w\sqrt{1-4\chi }(1+(1-4\chi )\left[ G(\chi )-E_{2}\left( f\right) \right]
E_{2}(\tau )w^{2})+O(w^{5}).
\end{eqnarray*}%
Finally, $2\,\pi i\,\tau _{{2}}$ is given by
\begin{eqnarray*}
&&2\,\pi i\Omega _{22}-E_{2}\left( \Omega _{11}\right) {(2\pi i\Omega _{12})}%
^{2}+5E_{2}\left( \Omega _{22}\right) E_{4}\left( \Omega _{11}\right) {(2\pi
i\Omega _{12})}^{4}+O({\Omega _{12}^{6}}) \\
&=&\log f(\chi )+(1-4\chi )E_{2}(\tau )w^{2}+\left( 1-4\,\chi \right)
^{2}\left( G(\chi )E_{2}\left( \tau \right) ^{2}+\frac{1}{2}E_{4}\left( \tau
\right) \right) w^{4} \\
&&-\left[ E_{2}\left( \tau \right) +(5E_{4}(\tau )-E_{2}(\tau )^{2})(1-4\chi
)G(\chi )w^{2}\right] \\
&&.w^{2}(1-4\chi )(1+2(1-4\chi )G(\chi )E_{2}(\tau )w^{2}) \\
&&+5E_{2}\left( f\right) E_{4}\left( \tau \right) (1-4\chi
)^{2}w^{4}+O(w^{6}) \\
&=&\log f(\chi )+\left( \frac{1}{2}+5[E_{2}(f)-G(\chi )]\right) E_{4}\left(
\tau \right) (1-4\chi )^{2}w^{4}+O(w^{6}).
\end{eqnarray*}%
Thus Proposition \ref{Proposition_eps_rho_w4} is proved. $\square $

\subsection{Corrections}

We list here some corrections to \cite{MT1} that we needed above. All
references below are to \cite{MT1}.

(a) Display (27) should read
\begin{equation}
\epsilon (\alpha ,-\alpha )=\epsilon (\alpha ,\alpha )=(-1)^{(\alpha ,\alpha
)/2}.  \label{A.1}
\end{equation}

(b) Display (45) should read
\begin{equation}
\gamma (\Xi )=(a,\delta _{r,1}\beta +C(r,0,\tau ))\alpha _{k}+\sum_{l\neq
k}D(r,0,z_{kl},\tau ).  \label{A.2}
\end{equation}

(c) As a result of (a), displays (79) and (80) are modified and now read
\begin{eqnarray}
F_{N}(e^{\alpha },z_{1};e^{-\alpha },z_{2};q) &=&\epsilon (\alpha ,-\alpha )%
\frac{q^{(\beta ,\beta )/2}}{\eta ^{l}(\tau )}\frac{\exp ((\beta ,\alpha
)z_{12})}{K(z_{12},\tau )^{(\alpha ,\alpha )}},  \label{A.3} \\
F_{V_{L}}(e^{\alpha },z_{1};e^{-\alpha },z_{2};q) &=&\epsilon (\alpha
,-\alpha )\frac{1}{\eta ^{l}(\tau )}\frac{\Theta _{\alpha ,L}(\tau
,z_{12}/2\pi i)}{K(z_{12},\tau )^{(\alpha ,\alpha )}}.  \label{A.4}
\end{eqnarray}

\end{document}